\documentclass[11pt]{amsart}

\usepackage{mathrsfs}
\usepackage{amsmath, amscd, amsthm,amssymb, amsfonts, verbatim,subfigure}
\usepackage[mathcal]{eucal}

\usepackage[all]{xy}

\usepackage{mathpple}

\title[Renormalisation and the Batalin-Vilkovisky formalism]{Renormalisation and the Batalin-Vilkovisky formalism}
\author{Kevin Costello}
\address{Department of Mathematics \\ Northwestern University}
\email{costello@math.northwestern.edu}

\date{}

\renewcommand{\Re}{\op{Re} } 

\newcommand{\eps}{\varepsilon}
\renewcommand{\epsilon}{\varepsilon}
\newcommand{\xto}{\xrightarrow}
\newcommand{\E}{\mscr{E}}
\newcommand{\what}{\widehat}

\newcommand{\mscr}{\mathscr}
\newcommand{\br}{\overline}
\newcommand{\iso}{\cong}
\newcommand{\C}{\mathbb C}

\newcommand{\norm}[1]{\left\| #1 \right\|}
\newcommand{\Oo}{\mscr O}
\newcommand{\Z}{\mathbb Z}
\newcommand{\defeq}{\overset{\text{def}}{=}}

\newcommand{\op}{\operatorname}
\newcommand{\mbf}{\mathbf}
\newcommand{\mbb}{\mathbb}
\newcommand{\mf}{\mathfrak}
\newcommand{\mc}{\mathcal}

\newcommand{\ip}[1]{\left\langle #1 \right\rangle}
\newcommand{\abs}[1]{\left| #1 \right|}

\newcommand{\R}{\mbb R}
\renewcommand{\d}{\mathrm{d}}

\newcommand{\liminv}{ \varprojlim }
\newcommand{\limdir}{\varinjlim}

\newcommand{\dbar}{\br{\partial}}
\renewcommand{\Im}{\op{Im}}

\DeclareMathOperator{\dens}{Densities}
\DeclareMathOperator{\Aut}{Aut} \DeclareMathOperator{\End}{End}
 
\DeclareMathOperator{\Sym}{Sym} \DeclareMathOperator{\Hom}{Hom}
 
\DeclareMathOperator{\Diff}{Diff}   
\DeclareMathOperator{\Tr}{Tr} 
\DeclareMathOperator{\Ker}{Ker}

\newtheoremstyle{thm}
  {7pt}
  {7pt}
  {\itshape}
  {}
  {\bf}
  {.}
  {5pt}
  {\thmnumber{#2 }\thmname{#1}\thmnote{ (#3)}}

\newtheoremstyle{def}
  {7pt}
  {10pt}
  {\itshape}
  {}
  {\bf}
  {.}
  {5pt}
  {\thmnumber{#2} \thmname{#1}\thmnote{ (#3)}}

\newtheoremstyle{rem}
  {4pt}
  {7pt}
  {}
  {}
  {\itshape}
  {:}
  {3pt}
  {}

\theoremstyle{thm}

\newtheorem*{utheorem}{Theorem}

\newtheorem*{uproposition}{Proposition}

\newtheorem*{thmA}{Theorem A}
\newtheorem*{thmB}{Theorem B}
\newtheorem*{thmC}{Theorem C}
\newtheorem*{thmD}{Theorem D}

\newtheorem*{ucorollary}{Corollary}
\newtheorem*{ulemma}{Lemma}

\newtheorem{theorem}{Theorem}[subsection]
\newtheorem{thm-def}{Theorem/Definition}[theorem]
\newtheorem{proposition}[theorem]{Proposition}
\newtheorem{lemma}[theorem]{Lemma}

\newtheorem{corollary}[theorem]{Corollary}

\numberwithin{equation}{subsection}

\theoremstyle{def}

\newtheorem{definition}[theorem]{Definition}

\theoremstyle{rem}

\newtheorem*{remark}{Remark}

\newcommand{\cinfty}{C^{\infty}}
\newcommand{\GF}{Q^{GF}}

\begin{document}

\begin{abstract}
This paper gives a way to renormalise certain quantum field theories on compact manifolds.  Examples include Yang-Mills theory (in dimension 4 only), Chern-Simons theory and holomorphic Chern-Simons theory.  The method is within the framework of the Batalin-Vilkovisky formalism.    Chern-Simons theory is renormalised in a way respecting all symmetries (up to homotopy).  This yields an invariant of smooth manifolds:  a certain algebraic structure on the cohomology of the manifold tensored with a Lie algebra, which is a ``higher loop'' enrichment of the natural $L_\infty$ structure.
\end{abstract}

\maketitle
\section{Introduction}

This paper gives a method for renormalising  a class of quantum field theories.  The field theories we are going to consider have space of fields of the form $\E = \Gamma(M,E)$, where $M$ is a compact manifold and $E$ is a super vector bundle on $M$.    We work within the Batalin-Vilkovisky formalism, so that $\E$ is equipped with an odd symplectic pairing $\E \otimes \E \to \C$ satisfying a certain non-degeneracy condition\footnote{Here, and throughout, $\otimes$ refers to the completed projective tensor product, so that $\E \otimes \E = \Gamma (M^2, E \boxtimes E)$.}.  $\E$ will also be equipped with a differential $Q$, which is an odd differential operator $Q : \E \to \E$ which is skew self adjoint and of square zero.  The action functionals in our quantum field theory will be of the form
$$
\frac{1}{2} \ip{e, Q e } + S(e )
$$
where $S$ is a local functional on the space of fields $\E$, which is at least cubic.  

The functional integrals we will renormalise are of the form
$$
Z(S, \hbar, a) = \int_{x \in L} \exp \left( \frac{1}{2  \hbar} \ip{x, Q x } + \frac{1}{\hbar} S(x + a )  \right) \d x
$$
where $a \in \E$ and $L \subset \E$ is an isotropic linear subspace, such that the map $Q : L \to \Im Q$ is an isomorphism.   Such an $L$ is known as a \emph{gauge fixing condition}.  $Z(S,a)$ is a formal functional of the variable $a \in \E$, and can be viewed as a generating function for certain Green's functions of the theory.

A fairly wide class of theories can be put in the form we use, including pure Yang-Mills theory in dimension 4, and Chern-Simons theory in any dimension.

This introduction will give a sketch of the results and of the underlying philosophy, without worrying too much about technical details. 

\subsection{}
Our gauge fixing conditions are always of the form
$$
L = \Im \GF
$$
where $\GF : \E \to \E$ is an odd self adjoint differential operator of square zero, with the property that the super-commutator
$$
H \defeq [Q,\GF]
$$
 is a positive elliptic operator of second order.  The facts that $\GF$ is self-adjoint with respect to the symplectic pairing on $\E$, and that $\GF$ is of square zero, imply that $L = \Im \GF$ is an isotropic subspace.

Only certain theories admit gauge fixing conditions of this form (this is the main restriction on the kind of theories the techniques from this paper can treat).  In many examples, $Q$ is a first order elliptic operator, and $\GF$ is a Hermitian adjoint of $Q$.    

A basic example to bear in mind is Chern-Simons theory in dimension $3$.  If $M$ is a compact oriented $3$-manifold, and $\mf g$ is a Lie algebra with an invariant non-degenerate pairing, then 
$$
\E = \Omega^\ast(M) \otimes \mf g[1].
$$
$[1]$ denotes parity shift.  The symplectic pairing on $\E$ arises from the Poincar\'e pairing on $\Omega^\ast(M)$ and the pairing on $\mf g$.  The operator $Q$ is simply the de Rham differential $\d_{DR}$, and $S$ is the cubic term in the standard Chern-Simons action.  The choice of a metric on $M$ leads to a gauge fixing condition $\GF = \d_{DR}^\ast$.  Further examples, including Chern-Simons theory in other dimensions, will be discussed later.

\subsection{}
Let us write 
$$
P(\eps,T) =  \int_{\eps}^T (\GF \otimes 1) K_t \d t \in \E \otimes \E
$$
where $K_t \in \E \otimes \E$ is the heat kernel for $H  = [Q, \GF]$. Here, $\E \otimes \E$ denotes the completed projective tensor product, 
$$
\E \otimes \E = \cinfty(M\times M, E \boxtimes E ). 
$$
The propagator of our theory is 
$$
P (0,\infty) = \int_{0}^\infty (\GF \otimes 1) K_t \d t. 
$$
This is not an element of $\E \otimes \E$, because of the singularities in the heat kernel at $t = 0$. Instead, $P(0,\infty)$ is an element of a distributional completion of $\E \otimes \E$. 

Let $\Oo(\E)$ denote the algebra of functions on $\E$,
$$
\Oo(\E) = \prod_{n \ge 0} \Hom ( \E^{\otimes n}, \C )_{S_n}
$$
where, as above, the tensor products are completed projective tensor products, and $\Hom$ denotes continuous linear maps. 

Any element $P \in \E \otimes \E$ gives rise to an order two differential operator $\partial_P$ on $\Oo(\E)$ in a standard way.

Up to a constant factor (the determinant of $Q$), one can write our functional integral as
$$
Z(S, \hbar, a) =  \lim_{\eps \to 0} \left( e^{\hbar \partial_{P(\eps,\infty ) } } e^{S / \hbar } \right) (a).
$$
The right hand side is the exponential of a differential operator applied to a function on $\E$, yielding a function on $\E$. This identity is rather formal; in finite dimensions, it is a simple consequence of integration by parts. In infinite dimensions we take it as an attempt at a definition. 

 When $\eps > 0$, the right hand side of this equation is well-defined. However, the $\eps \to 0$ limit is singular, because $P(0,\infty)$ is not an element of $\E \otimes \E$, but has singularities  along the diagonal in $M^2$.

Let us use the notation
$$
\Gamma ( P(\eps,T) , S  )  =  \hbar \log \left ( e^{\hbar \partial_{P(\eps,\infty )} } e^{S / \hbar }\right). 
$$
This is an $\hbar$ dependent function on $\E$, that is, an element of $\Oo(\E) [[\hbar]]$.  We will typically omit the variables $a \in \E$ and $\hbar$ from the notation.

The expression we would like to make sense of is 
$$\hbar \log Z(S, \hbar, a) =  \lim_{\eps \to 0} \Gamma( P(\eps,\infty), S).$$

\subsection{}
 An effective action\footnote{ What I mean by effective action is related to the Wilsonian effective action.  Ignoring for the moment considerations of renormalisation, the Wilsonian effective action is obtained by integrating out all the high-energy fields, i.e.\ all the eigenspaces of $H$ with high eigenvalues.      The effective action considered here is obtained by averaging over all interactions occurring at small scales.  More precisely, the effective action at scale $\eps$ is the sum over all Feynman graphs of the theory, using the propagator $P(0,\eps)$. Using the propagator $P(0,\eps)$ amounts to allowing particles to propagate for a distance of between $0$ and $\eps$.   } at scale $\eps$ is a function $f \in \Oo(\E) [[\hbar]]$ which describes all interactions occurring at a scale below $\eps$.   One can reconstruct the effective action at any other scale using the effective action at scale $\eps$ and the propagator.  The map $f \mapsto \Gamma ( P(\eps,T) , f)$ is the operation taking an effective action at scale $\eps$ to the corresponding effective action at scale $T$.  One can imagine that the scale $T$ effective action $\Gamma(P(\eps,T), f)$ is obtained from the scale $\eps$ effective action $f$ by allowing particles to interact according to $f$, and to propagate a distance between $\eps$ and $T$. 

$\Gamma ( P(\eps,T) , f)$ is the renormalisation group flow from scale $\eps$ to scale $T$ applied to $f$. This is well-defined for all those $f \in \Oo(\E)[[\hbar]]$ which are at least cubic modulo $\hbar$, as long as $0 < \eps < T \le \infty$.     We have the semi-group law
$$
\Gamma ( P(T_2,T_3 ) , \Gamma ( P(T_1,T_2) , f) ) = \Gamma ( P(T_1,T_3) , f).
$$
This equation is  a version of the exact renormalisation group equation.  The operation $f \to \Gamma( P(\eps,T), f)$ is invertible, so that if we know the effective action at any scale, we know it at all other scales. 

The only part of this renormalisation group flow that is problematic is taking an effective action at scale $0$ and turning it into an effective action at any positive scale $\eps$. This part of the procedure needs to be renormalised.   This is to be expected: an effective action at scale $0$ would describe interactions occurring at infinitely high energy.  

One way to phrase the main result we prove is that there is a bijection between functionals $S \in \Oo(\E) [[\hbar]]$, satisfying a locality axiom, and systems $\{ S^{eff}(T) \mid T \in \R_{> 0} \}$ of effective actions at all scales $T > 0$, related by the renormalisation group equation. The effective action $S^{eff}(T)$ must also satisfy a locality condition as $T \to 0$.  The original action $S$ plays the role of the scale $0$ effective action, and the positive scale effective actions $S^{eff}(T)$ are obtained by renormalising the expression $\Gamma ( P( 0,T) , S)$.  Every such system of effective actions $S^{eff}(T)$ arises from a unique local functional $S$ in this way.

\subsection{}
In order to renormalise the expression $\lim_{\eps \to 0} \Gamma ( P(\eps, T ) , S)$, and thus construct the scale $T$ effective action, we will need a way to extract the ``singular part'' of expressions of the form $\Gamma ( P(\eps, T ) , S)$. This will rely on some results about the small $\eps$ asymptotic expansion of $\Gamma ( P(\eps, T ) , S)$.

It's convenient to represent $\Gamma ( P(\eps, T ) , S)$ as 
$$
\Gamma ( P(\eps, T ) , S)= \sum_{i,k \ge 0} \hbar^i \Gamma_{i,k} ( P(\eps, T ) , S)
$$
where $\Gamma_{i,k} ( P(\eps, T ) , S)$ is homogeneous of degree $k$ as a function of $a \in \E$.  This expression is just the Taylor expansion of $\Gamma ( P(\eps, T), S)$ as a function of $\hbar$ and $a \in \E$.   In terms of Feynman graphs, $\Gamma_{i,k}(P(\eps,T),S)$ is the sum over all graphs with first Betti number $i$  and $k$ external edges. 

Let $\op{An} ( (0,\infty) )$ be the algebra of analytic functions on $(0,\infty)$, where $\eps$ is the coordinate on $(0,\infty)$.  
\begin{thmA}
There exists a subalgebra $\mscr{A} \subset \op{An} ( (0,\infty) )$ with a countable basis, such that for all local functionals $S \in \Oo(\E) [[\hbar]]$, there exists 
a small $\epsilon$ asymptotic expansion
$$
\Gamma_{i,k} ( P(\eps, T ) , S) \simeq  \sum f_r(\epsilon) \otimes \Phi_{i,k,r}(T,a) 
$$
where $f_r \in \mscr{A}$, and $\Phi_{i,k,r}$ are certain functions of the variables $T \in (0,\infty)$ and $a \in \E$. 
\end{thmA}
A functional $S \in \Oo(\E) [[\hbar]]$ is \emph{local} if, roughly, its homogeneous components $S_{i,k} \in \Hom(\E^{\otimes k}, \C)_{S_k}$, which are distributions on the vector bundle $E^{\boxtimes k}$ on $M^k$, are supported on the small diagonal, and non-singular in the diagonal directions.  We will denote the space of local functionals by 
$$
\Oo_l(\E) \subset \Oo(\E) . 
$$

Let $\mscr{A}_{\ge 0} \subset \mscr{A}$ be the subspace of those functions whose $\epsilon \to 0$ limit exists.  In order to extract the singular part of functions in $\mscr{A}$, we need to pick a complementary subspace to $\mscr{A}_{\ge 0}$.
\begin{definition}
A \emph{renormalisation scheme} is a subspace $\mscr{A}_{< 0}\subset \mscr{A}$ such that 
$$
\mscr{A}  = \mscr{A}_{\ge 0} \oplus \mscr{A}_{< 0} . 
$$
\end{definition}
Let us fix a renormalisation scheme $\mscr{A}_{< 0}$.  Later we will see that things are independent in a certain sense of the choice of renormalisation scheme.  

\begin{remark}
The functions in $\mscr{A}$ are quite explicit integrals of multi-variable rational functions.   The algebra $\mscr{A}$ only depends on the dimension of the manifold $M$, and not on the details of the particular theory we are working with. 
\end{remark}
\begin{remark}
Instead of using the algebra $\mscr{A}$, one could use any larger algebra of functions on $(0,\infty)$, and obtain the same results.  It is technically easier to use an algebra with a countable basis. 
\end{remark}

\begin{remark}
An alternative regularisation scheme, which Jack Morava suggested to me, would be to use the propagator $\int_0^\infty t^z (\GF \otimes 1)K_t \d t$, where $z$ is a complex parameter.  If we use this propagator, then integrals attached to Feynman graphs converge if $\Re z \gg 0$.  The expressions should admit an analytic continuation to $\C$ with poles on $\tfrac{1}{2} \Z$.  Unfortunately, I wasn't able to prove the existence of the analytic continuation. 
\end{remark}

\subsection{}
The first main result of this paper is the following.
\begin{thmB}
Fix a renormalisation scheme $\mscr{A}_{< 0} \subset \mscr{A}$. 
Then there exists a unique series of counter-terms
$$
S^{CT} (\hbar, \epsilon ,a ) = \sum_{i > 0, k\ge 0} \hbar^i S^{CT}_{i,k} (\epsilon, a )
$$
where 
\begin{enumerate}
\item
each $S^{CT}_{i,k}(\epsilon,a)$ is a formal local functional of $a \in \E$, homogeneous of order $k$ as a function of $a$, with values in $\mscr{A}_{< 0}$ Thus,
$$
S^{CT}_{i,k} \in \Oo_l(\E) \otimes \mscr{A}_{< 0}
$$ 
where $\Oo_l(\E) \subset \Oo(\E)$ is the space of local functionals on $\E$.
\item
the limit 
$$\lim_{\epsilon \to 0} \Gamma ( P(\eps, T ) , S - S^{CT})$$ 
exists. 
\end{enumerate}
These counter-terms are independent of $T$. 
\end{thmB}

Let me give a brief sketch of the (surprisingly simple) proof of this theorem. As before, let us write 
$$
\Gamma ( P(\eps, T ) , S)= \sum_{i,k \ge 0} \hbar^i \Gamma_{i,k} ( P(\eps, T ) , S).
$$
The $\Gamma_{0,k} ( P(\eps, T ) , S)$ all have well-defined $\epsilon \to 0$ limits.  So the first counter-term we need to construct is $S^{CT}_{1,1}$.  Our choice of renormalisation scheme allows us to extract the singular part of any function of $f(\epsilon) \in \mscr{A}$; this singular part is simply the projection of $f$ onto $\mscr{A}_{< 0}$.  We define
$$
S^{CT}_{1,1}(\epsilon,a) = \text{singular part of the small $\epsilon$ expansion of } \Gamma_{1,1} ( P(\eps, T ) , S)\in \mscr{A}_{< 0}. 
$$
It is easy to see that $\Gamma_{1,1} ( P(\eps, T ) , S- \hbar S^{CT}_{1,1})$ is non-singular as $\epsilon \to 0$.  

The next step is to replace $S$ by $S - \hbar S^{CT}_{1,1}$, and use this new action to produce the next counter-term, $S^{CT}_{1,2}$.  That is, we define
$$
S^{CT}_{1,2}= \text{singular part of } \Gamma_{1,2} ( P(\eps, T ) , S - \hbar S_{1,1}^{CT})
$$
where we understand ``singular part'' in the same way as before. We continue like this, to define $S^{CT}_{1,k}$ for all $k \ge 1$.

The next step is to define
$$
S^{CT}_{2,0} = \text{singular part of } \Gamma_{2,0} ( P(\eps, T ) , S - \hbar \sum_{k = 0}^\infty S_{1,k}^{CT})
$$
Continuing in this manner defines all the $S^{CT}_{i,k}$. 

The difficult part of the proof is the verification that the counter-terms $S^{CT}_{i,k}$ are local functionals.  Locality is desirable for many physical and mathematical reasons. More practically, we need $S^{CT}_{i,k}$ to be local in order to apply the procedure at the next inductive step. We only know the existence of a small $\eps$ asymptotic expansion of the $\Gamma_{I,K}(P(\eps,T) , S - \sum \hbar^i S^{CT}_{i,k} )$  when the counter-terms $S^{CT}_{i,k}$ are local. 

The main step in the proof of locality is showing that the $S^{CS}_{i,k}$ are independent of $T$.  Once we know they are independent of $T$, we can take $T \to 0$.  $\Gamma_{i,k}( P(\eps, T) , S)$ is concentrated near the diagonal in $M^k$, if $\eps < T$ and both $\eps, T$ are very small.   Thus, the counter-terms become supported on the diagonal, and local.

This theorem allows one to define unambiguously the renormalised scale $T$ effective action
$$
\Gamma^R ( P(0,T) , S ) =  \lim_{\eps \to 0} \Gamma( P(\eps,T) , S - S^{CT} ) .
$$
Thus, we can define the renormalised functional integral
\begin{multline*}
Z^R =  \exp  ( \Gamma^R ( P(0,\infty), S)  / \hbar ) = \text { renormalisation of }  \\
\int_{x \in L} \exp \left( \frac{1}{2  \hbar} \ip{x, Q x } + \frac{1}{\hbar} S(x + a )  \right) \d x.
\end{multline*}
This renormalised partition function is an element of $\Oo(\E)((\hbar))$, that is, a  non-local formal function on the space $\E$ of fields.

\subsection{Independence of choice of renormalisation scheme}
We should interpret the expression $\Gamma^R ( P(0,T) , S ) $ constructed using Theorem B  as the scale $T$ renormalised effective action.  The renormalisation group equation holds:
$$
\Gamma ( P(T_1,T_2 ) , \Gamma^R ( P(0,T_1) , S) ) = \Gamma^R ( P(0,T_2) , S).
$$
\begin{definition}
A \emph{system of effective actions} on the space of fields $\E$ is given by an effective action
$$
S^{eff}(T) \in \Oo(\E) [[\hbar]] 
$$
for all $T \in \R_{> 0}$, varying smoothly with $T$, 
such that
\begin{enumerate}
\item
Each $S^{eff}(T)$ is at least cubic modulo $\hbar$. 
\item
The renormalisation group equation is satisfied,
$$
S^{eff}(T_2) =  \Gamma ( P(T_1,T_2) , S^{eff}(T_1)).
$$
\item
As $T \to 0$, $S^{eff}(T)$ must become local, in the following sense.  There must exist some $T$-dependent local functional $\Phi(T)$ such that $\lim_{T \to 0} \left( S^{eff}(T) - \Phi(T) \right) = 0$. (The $T \to 0$ limit of $S^{eff}(T)$ itself will generally not exist). 
\end{enumerate}
\end{definition}
The renormalised effective actions $\Gamma^R (  P ( 0,T)  , S)$ constructed from a local functional $S$ satisfy these two axioms. Thus, for any renormalisation scheme $\mscr{A}_{< 0}$, theorem B provides a map
\begin{align*}
\text{local functionals } S \in \Oo_l(\E) [[\hbar]] & \to \text{ systems of effective actions } \\
S & \mapsto \{ \Gamma^R ( P(0,T), S ) \mid T \in \R_{> 0} \} .
\end{align*}
(The local functionals $S$ must be at least cubic modulo $\hbar$, as must the effective actions $S^{eff}(T)$. ) .
\begin{thmC}
For any renormalisation scheme $\mscr{A}_{< 0}$, this map is a bijection.
\end{thmC}
This set of systems of effective actions on $\E$ is a canonical object associated to $(\E,Q, \GF)$, independent of the choice of renormalisation scheme.    Renormalisation and regularisation techniques other than those considered should lead to different ways of parametrising the same set of systems of effective actions.  For instance, if one could make sense of dimensional regularisation on general manifolds, one would hope to get simply a different parametrisation of this set.

From this point of view, the formalism of counter-terms is simply a convenient way to describe this set of systems of effective actions.   The counter-terms themselves, and the original action $S$, are not really meaningful in themselves.   The main point of introducing counter-terms is that it is otherwise difficult to produce systems of effective actions.   \emph{A priori}, it is not obvious that there are any non-zero systems of effective actions at all.

This proposition makes clear in what sense renormalisation is independent of the choice of renormalisation scheme.  To any renormalisation scheme $\mscr{A}_{< 0}$ and local funct\-ional $S \in  \Oo_l(\E) [[\hbar]]$ corresponds a ``theory'', i.e.\  a system of effective actions.  If $\mscr{A}'_{< 0}$ is another renormalisation scheme, there is a unique local functional $S'$ such that $(S', \mscr{A}'_{< 0})$ gives the same theory as $(S, \mscr{A}_{< 0})$.   

This statement can be expressed more formally as follows. Let $\op{RS}$ denote the space of renormalisation schemes, and let $\op{EA}$ denote the set of systems of effective actions.  Theorem C implies gives an isomorphism of fibre bundles on $\op{RS}$
$$
\Oo_l(\E) \times \op{RS} \to \op{EA} \times \op{RS} .
$$
Give the right hand side the trivial flat connection; this pulls back to a non-trivial (and non-linear) flat connection on $\Oo_l(\E)$.   This flat connection is uniquely characterised by the property that for any flat section $S : \op{RS} \to \Oo_l(\E) \times \op{RS}$, the system of effective actions $\{\Gamma^R ( P(0,T) , S (\mscr{A}_{< 0} ) ) \}$ associated to the functional $S(\mscr{A}_{< 0} ) \in \Oo_l(\E)$ is independent of the point $\mscr{A}_{< 0}  \in \op{RS}$.  

We will fix once and for all a renormalisation scheme $\mscr{A}_{< 0}$. This allows us to always talk about local functionals, as a convenient proxy for the set of systems of effective actions. The choice of $\mscr{A}_{< 0}$ is analogous to the choice of a basis in a vector space. 

All statements about a theory should be expressed in terms of the effective actions $\Gamma^R(P(0,T), S)$, and not directly in terms of the local functional $S$.  This will ensure everything is independent of the choice of renormalisation scheme.

\subsection{}

 Bogoliubov and Parasiuk \cite{BogPar57}, Hepp \cite{Hep66} and Zimmerman \cite{Zim69}   have given an algorithm for the renormalisation of certain quantum field theories.   Their algorithm is based on combinatorial manipulations of Feynman graphs.   Recently, Connes and Kreimer \cite{ConKre00} have given a beautiful interpretation of the BPHZ algorithm, in terms of the Birkhoff factorisation of loops in a certain Hopf algebra. 
  
In the approach used in this paper, no graph combinatorics are needed; all we use is a very simple inductive argument, sketched above.  The reason that things become so simple seems to be the particular kind of functional integrals we consider, which are always of the form
$$
Z(S, \hbar, a) = \int_{x \in \Im \GF} \exp \left( \frac{1}{2  \hbar} \ip{x, Q x } + \frac{1}{\hbar} S(x + a )  \right) \d x.
$$
Thus, the moral of this paper is that if we consider functional integrals of this form,
then the problem of renormalisation becomes quite simple, and the counter-terms are unique and automatically local.      As we will see shortly, the particular functional integrals we use also play an important role in the Batalin-Vilkovisky formalism.

Another difference between the approach to renormalisation described in this paper and that of Connes-Kreimer and BPHZ is that we do not give finite values to individual Feynman  graphs, but only to the sum over all graphs with a fixed number of loops and external edges. From the point of view of the effective action,  the expression attached to an individual graph has no meaning.

If we tried to renormalise different classes of functional integrals the procedure would not be so simple.  For example, if we try to simply renormalise the integral 
$$\int_{x \in \Im \GF} \exp \left( \frac{1}{2  \hbar} \ip{x, Q x } + \frac{1}{\hbar} S(x )  \right) \d x,$$
without using the variable $a$, then there are many possible choices of counter-terms.  

If we try to renormalise the functional integral $$\int_{x \in \Im \GF} \exp \left( \frac{1}{2  \hbar} \ip{x, Q x } + \frac{1}{\hbar} S(x) + \frac{1}{\hbar} \ip{x,a}  \right) \d x, $$  then the terms in the Feynman graph expansion don't fit together in the right way to produce local counter-terms.   

This type of functional integral is one which appears more commonly in the literature.  A simple change of variables allows one to express this type of functional integral in terms of the kind considered here, but not conversely.  Indeed, if $a =- Q^{-1} \pi b$, where $\pi$ is the projection onto $\Im \GF$ and $Q^{-1} :  \Im Q \to \Im \GF$ is the inverse to $Q$, we can write
\begin{multline*}
\int_{x \in  \Im \GF } \exp \left( \frac{1}{2  \hbar} \ip{x, Q x } + \frac{1}{\hbar} S(x + a)  \right) \d x \\
\shoveright{ =\int_{x \in  \Im \GF } \exp \left( \frac{1}{2  \hbar} \ip{x- a, Q (x-a) } + \frac{1}{\hbar} S(x )  \right) \d x  \hspace{41pt}} \\
=  e^{- \ip{a, b } /\hbar } \int_{x \in  \Im \GF } \exp \left( \frac{1}{2  \hbar} \ip{x, Q x } + \frac{1}{\hbar} S(x) + \frac{1}{\hbar} \ip{x,b}  \right) \d x 
\end{multline*}

\subsection{Batalin-Vilkovisky formalism}
The Batalin-Vilkovisky quantum master equation is the quantum expression of gauge symmetry.  It takes the form
$$
(Q + \hbar \Delta ) \exp ( S / \hbar ) = 0
$$
where $\Delta$ is a certain order $2$ differential operator acting on the space of functionals on $\E$.  This expression makes perfect sense in the simple situation when the space of fields $\E$ is  finite dimensional (i.e.\ the underlying manifold $M$ is of dimension $0$).  In the situation we work in, however, this expression is infinite. This is because  $\Delta S$ involves the multiplication of singular distributions, and thus has the same kind of singularities as appear in one-loop Feynman diagrams.  

This paper gives a definition of a renormalised quantum master equation that works in the infinite dimensional situation.  There are regularised BV operators $\Delta_t$, for $t > 0$, defined by 
$$
\Delta_t = - \partial_{K_t}.
$$
Recall that $K_t \in \E \otimes \E$ is the heat kernel for $H = [Q, \GF]$, and $\partial_{K_t}$ is the order two differential operator on the algebra $\Oo(\E)$ of functions on $\E$, associated to $K_t$.  The operators $\Delta_t$ are thus differential operators on $\Oo(\E)$, for all $t > 0$.   The ``physical'' BV operator is $\Delta_0$, which is ill-defined.

\begin{lemma}	
A function $f \in \Oo(\E)[[\hbar]]$ satisfies the $\Delta_\eps$ quantum master equation
$$
(Q + \hbar \Delta_\eps)e^{f / \hbar} = 0
$$
 if and only if 
$$
\Gamma ( P(\eps, T ) , f)
$$
satisfies the $\Delta_T$ quantum master equation.  
\end{lemma}
This follows from the fact that
$$
Q P(\eps, T) = -K_T + K_\eps
$$
so that the differential operator $\partial_{P(\eps,T)}$ is a chain homotopy between $\Delta_\eps$ and $\Delta_T$.

This lemma motivates the following definition.
\begin{definition}
A local functional $S \in \Oo_l(\E)[[\hbar]]$ satisfies the renormalised quantum master equation if the renormalised expression 
$$\Gamma^R ( P(0,T) , S)  = \lim_{\eps \to 0} \Gamma (P(\eps,T), S - S^{CT} )$$
satisfies the $\Delta_T$ quantum master equation, for some (or equivalently, all) $T > 0$. 
\end{definition}
Thus, the renormalised QME says that the scale $T$ effective action associated to $S$ satisfies the scale $T$ quantum master equation.

One peculiarity of this renormalised quantum master equation is that, unlike in the finite dimensional situation, it depends on the gauge fixing condition $\GF$.    We will see shortly that this dependence is very weak.  

In fact,  the equation depends also on the renormalisation scheme $\mscr{A}_{< 0}$, but only in an artificial way.  Recall that theorem C says that the choice of a renormalisation scheme sets up a bijection between local functionals $S \in \Oo_l(\E)[[\hbar]]$, and systems of effective actions.       The system of effective actions is given by $\{\Gamma^R(P(0,T), S) \mid T \in \R_{> 0} \}$.   The quantum master equation as a statement about the system of effective actions is independent of the choice of renormalisation scheme.   

We will fix once and for all a renormalisation scheme, and use it to parametrise the set of systems of effective actions.  If we use a different renormalisation scheme, everything works the same, except we are parametrising the set of systems of effective actions in a slightly different way. 

\subsection{Homotopies of solutions to the renormalised quantum master equation}
Let's consider the finite dimensional situation again for a moment, with the further assumption that the finite dimensional space of fields $V$ has trivial $Q$ cohomology.  Then the subspace $L \subset V$ is a Lagrangian, and not just isotropic; we have a direct sum decomposition 
$$
V =  L \oplus \Im Q.
$$
The importance of the Batalin-Vilkovisky quantum master equation comes from the fact that in this situation, if $S$ satisfies the Batalin-Vilkovisky quantum master equation, then the partition function $Z(S,\hbar, a = 0)$  remains unchanged under continuous variations of $L$.

We would like to prove a version of this in the infinite dimensional situation, for the renormalised quantum master equation.  However, things are more delicate in this situation; the renormalised QME itself depends on the choice of a gauge fixing condition.  What we will show is that if $\GF(t)$ is a one-parameter family of gauge fixing conditions, corresponding to the family of isotropic subspaces $\Im \GF(t)$,  then the set of homotopy classes of solutions to the renormalised QME using $\GF(0)$ is isomorphic to the corresponding set using $\GF(1)$.  This result is a corollary of a result about certain simplicial sets of gauge fixing conditions and of solutions to the renormalised QME.

Let me explain this picture in more detail.  One of the axioms we need for our gauge fixing conditions is that there is a direct sum decomposition
$$
\E = \Im \GF \oplus \Ker H \oplus \Im Q
$$ 
where $H = [\GF, Q]$ so that $\Ker H$ is the space of harmonic elements of $\E$.  
This leads to the identification
$$
H^\ast(\E, Q) = \Ker H.
$$
This cohomology group is finite dimensional.

There is a notion of homotopy of solutions of the quantum master equation. Briefly, two solutions $S_0, S_1$ of the quantum master equation are homotopic if there exists $S_t + \d t S'_t$, for $t \in [0,1]$, such that 
$$
\left(Q + \d t \frac{\d}{\d t} + \hbar \Delta \right) e^{(S_t + \d t S'_t )/\hbar } = 0. 
$$
In a similar way, one can define a notion of homotopy of solutions of the renormalised quantum master equation.

\begin{utheorem}
\begin{enumerate}
\item
If $S$ satisfies the renormalised quantum master equation, then the restriction of $\Gamma^R (P(0,\infty), S)$ to $H^\ast(\E, Q)$ satisfies the finite dimensional quantum master equation.   The map $S \to \Gamma^R (P(0,\infty), S)$ respects homotopies. 
\item
Let $\GF(t)  \subset \E$ be a smooth one-parameter family of  gauge fixing conditions, for $t \in [0,1]$.    Then there exists a natural bijection between the set of homotopy classes of solutions of the renormalised quantum master equation using the gauge fixing condition $\GF(0)$  and the corresponding set using $\GF(1)$.
\item
Let $S_0$, $S_1$ be solutions of the renormalised quantum master equations using $\GF(0)$ and $\GF(1)$ respectively, which correspond up to homotopy under the bijection coming from the family $\GF(t)$.  Then $\Gamma^R (P(0,\infty), S_0)$ and $\Gamma^R(P(0,\infty), S_1)$ (defined using $\GF(0)$ and $\GF(1)$ respectively) are homotopic solutions of the quantum master equation on $H^\ast(\E, Q)$.  
\end{enumerate}

\end{utheorem}
In fact, this result will be a corollary of a more abstract result concerning simplicial sets of solutions of the quantum master equation.
\begin{utheorem}
There exist simplicial sets
\begin{enumerate}
\item $\mathbf{GF}(\E,Q)$ of gauge fixing conditions 
\item $\mathbf{BV}( \E, Q )$ of solutions of the renormalised quantum master equation on $\E$ 
\item $\mathbf{BV}( H^\ast(\E, Q) )$ of solutions to the finite-dimensional quantum master equation on $H^\ast(\E,Q)$
\end{enumerate}
which fit  into a diagram of maps of simplicial sets
$$
\xymatrix{ \mathbf{BV}( \E, Q )  \ar[rr]^(.43){\Gamma^R(P(0,\infty), - ) } \ar[d]^{\pi} & &   \mathbf{BV}( H^\ast(\E, Q) )  \\ 
\mathbf{GF}( \E, Q)  &  }
$$
where the vertical arrow  $\pi$ is a fibration of simplicial sets.  
\end{utheorem}
A solution of the renormalised quantum master equation must be with respect to some gauge fixing condition; this defines the vertical arrow $\pi : \mathbf{BV}( \E, Q ) \to \mathbf{GF}(\E,Q)$.   The map $\mathbf{BV}( \E, Q ) \to   \mathbf{BV}( H^\ast(\E, Q) )$ is the simplicial version of the map discussed earlier, which takes a solution $S$ of the renormalised quantum master equation to $\Gamma^R(P(0,\infty), S) \vert_{H^\ast(\E,Q)}$.  This is a solution of the quantum master equation on $H^\ast(\E,Q)$.

One can deduce the previous more concrete result from this abstract statement about simplicial sets using the simplicial analogs of path and homotopy lifting properties for fibrations.  

\subsection{}
We have seen that the set of systems of effective actions is a canonical object associated to $(\E,Q, \GF)$, defined with out reference to a renormalisation scheme.  In a similar way,  we could say that the simplicial set $\mbf{BV}(\E,Q)$ is a canonical object associated to $(\E,Q)$ without reference to the choice of gauge fixing condition.  A choice of a renormalisation scheme and a gauge fixing condition gives us a convenient parametrisation of this simplicial set, as a set of local functionals  satisfying the renormalised quantum master equation.  However, if we choose a different renormalisation scheme, we get something canonically isomorphic; and if we choose a different gauge fixing condition we something canonically homotopy equivalent. (At least, this is true as long as the space of natural gauge fixing conditions is contractible, which it always seems to be in examples). 

If we have a classical action $S_0$ which solves the classical master equation, there is a simplicial set $\mbf{BV}(\E,Q, S_0)$ of quantisations of $S_0$, i.e.\ solutions to the renormalised quantum master of the form $S_0+ \hbar S_1 + \cdots$.    If the space of natural gauge fixing conditions is contractible, then this simplicial set is canonically associated to the classical theory $(\E,Q,S_0)$, up to canonical homotopy equivalences. 

Thus, there is a homotopy action of the group of symmetries of the classical theory on the simplicial set $\mbf{BV}(\E,Q, S_0)$ of quantisations.  One would like to quantise a given classical theory in a way preserving as many symmetries as possible.

\subsection{Quantisation of Chern-Simons theory}
The results of this paper allow one to renormalise a wide variety of quantum field theories, for example Chern-Simons theory on a compact oriented manifold, or pure Yang-Mills theory on a compact $4$-dimen\-sional manifold with a conformal class of metrics.  By ``renormalisation'' I simply mean the construction of counter-terms.

I would like to distinguish between renormalisation and quantisation.  By \emph{quantisation} I mean the replacement of an action $S_0$ by an action $S = S_0 + \sum_{i > 0} \hbar^i S_i$ which satisfies the renormalised quantum master equation.   The only non-trivial theory we succeed in quantising in this paper is Chern-Simons theory.  In fact, we only quantise Chern-Simons theory modulo the constant term (constant as a function on the space of fields). 

The quantisation of Chern-Simons theory is based on a general local-to-global principle, which allows one to construct global solutions to the renormalised QME from local ones.  This local-to-global result is stated and proved for a general class of theories in the body of the paper.  To keep things simple, I'll only discuss Chern-Simons theory in this introduction. 

Chern-Simons theory is a perturbative gauge theory associated to a compact oriented manifold $M$ and a flat bundle $\mf g$ of super Lie algebras on $M$, with an invariant pairing of parity opposite to that of $\dim M$. 

In this situation,  let
$$
\E = \Omega^\ast(M , \mf g ) [1].
$$
$\E$ is the Batalin-Vilkovisky odd symplectic manifold associated to the Chern-Simons gauge theory of connections with values in $\mf g$.  

A gauge fixing condition on $\E$ is given by a choice of a metric on $M$.  The space of metrics is of course contractible,.   As we have seen above, the spaces of homotopy classes of solutions to the renormalised QME for different gauge fixing conditions are homotopy equivalent.  Thus, one can speak about homotopy classes of solutions to the renormalised QME without reference to a metric. 
\begin{utheorem}
Let $M$, $\mf g$ be as above.  Then there is a canonical (up to a contractible choice) quantisation $S$ of the Chern-Simons action $S_0$ on $\E$, modulo constant terms. 

That is, there is a canonical up to homotopy functional $S = S_0 + \sum_{i > 0} \hbar^i  S_i$ on $\E$, where each $S_i$ is defined modulo constants (as a function on $\E$), which satisfies the renormalised quantum master equation. 
\end{utheorem}
I should emphasise that the specific form the $S_i$ will take is dependent on both the metric and the renormalisation scheme we choose to work with. If we use a different renormalisation scheme, then the $S_i$ will change, but the effective action $\Gamma^R ( P (0,T) , \sum \hbar^i S_i)$ remains unchanged.  If we use a different metric, then this effective action changes by a homotopy. 

\begin{ucorollary}
There is a canonical, up to homotopy and modulo constants, function $\Phi^{CS}$ on 
$$H^\ast(\E) = H^\ast(M, \mf g ) [1]$$
which satisfies the quantum master equation.     
\end{ucorollary}
As $\Phi^{CS}$ is well-defined modulo constants, $\Phi^{CS}$ is an element of 
$$
\Phi^{CS}  \in \left( \prod_{k > 0}  \Sym^k H^\ast(\E)^{\vee}  \right)[[\hbar]] .
$$
The quantum master equation is
$$
\{ \Phi^{CS}, \Phi^{CS} \} + \hbar \Delta \Phi^{CS} = 0
$$
which must hold modulo constants, that is, modulo $\C[[\hbar]]$.  

One can write this identity in more explicit terms.  The Hamiltonian vector field associated to $\Phi^{CS}$ has Taylor components which are linear maps
$$
\sum \hbar^i l_{i,k} : H^\ast(M, \mf g)^{\otimes k} \to H^\ast(M, \mf g) [[\hbar]]. 
$$
where $l_{i,k}$ is independent of $\hbar$.   The $l_{i,k}$ are defined for all $i, k \ge 0$.  The $l_{0,k}$ are the usual $L_\infty$ operations, which arise via the homological perturbation lemma.  The $l_{i,k}$ for $i  > 0$ are the new structure.  The quantum master equation is encoded in a sequence of identities of the form
\begin{multline*}
\sum_{\substack{i_1+i_2 = i \\ k_1 + k_2 = k-1  }} \pm l_{i_1,k_1}  ( x_1,\ldots, x_j,  l_{i_2,k_2}  
(x_{j+1}, \ldots, x_{j+i_1}   )  , x_{j + i_1 + 1}, \ldots, x_{i}       )  \\
+ \sum \pm  l_{i-1,k+2} ( x_1, \ldots, x_{j'}, \delta', x_{j'+1}, \ldots, x_{j''}, \delta'', x_{j''+1}, \ldots, x_i  ) = 0. 
\end{multline*}
for all $i,k$. 
Here, the $x_a \in H^\ast(M , \mf g)$, and $\sum \delta' \otimes \delta'' \in H^\ast(M, \mf g)^{\otimes 2}$ is the tensor inverse to the pairing.    The first term in this identity is the expression appearing in the usual $L_\infty$ equation.  This algebraic structure is sometimes called a ``quantum $L_\infty$ algebra''. 

As we have seen, the $l_{0,k}$ give the usual $L_\infty$ structure on $H^\ast(M, \mf g)$. These $L_\infty$ algebras, for varying $\mf g$,  encode a great deal of the rational homotopy type of $M$. Thus, the structure defined by the  the $l_{i,k}$ for $i > 0$ can be viewed as a kind of quantisation of the rational homotopy type. 

In the case when $H^\ast(\E)  = 0$,  Kontsevich \cite{Kon94} and Axelrod-Singer \cite{AxeSin92,AxeSin94} (when $\dim M = 3$) have already constructed the perturbative Chern-Simons invariants.  In some sense, their construction is orthogonal to the construction in this paper.   Because we work modulo constants, the construction in this paper doesn't give anything in the case when $H^\ast(\E) = 0$.     On the other hand, their constructions don't apply in the situations where our construction gives something non-trivial.    

There seems to be no fundamental reason why  a generalisation of the construction to this paper, including the constant term,  should not exist.  Such a generalisation would also generalise the results of Kontsevich and Axelrod-Singer.    However, the problem of constructing such a generalisation does not seem to be amenable to the  techniques used in this paper.

A theory related to the Chern-Simons theory considered here has been treated in the very interesting recent paper \cite{Mne06}. The BF theory used by Mnev is the same as Chern-Simons theory with a special kind of Lie algebra, of the following form.  Let $\mf g$ be any finite dimensional Lie algebra.  Then $\mf g \oplus \mf g^\vee$ is a Lie algebra with an even invariant pairing, and  $\mf g \oplus \mf g^\vee[1]$ is a Lie algebra with an odd invariant pairing.  The Lie algebra structure arises from the fact that $\mf g^\vee$ carries a $\mf g$ action.  Chern-Simons theory with Lie algebra $\mf g \oplus \mf g^\vee$ (when $\dim M$ is odd) or $\mf g \oplus \mf g^\vee[1]$  (when $\dim M$ is even) is the same as the BF theory considered by Mnev. 

\subsection{Construction of the quantisation of Chern-Simons theory}

Let me sketch the proof of the theorem on quantisation of Chern-Simons theory.  The proof uses the homotopical algebra of simplicial presheaves to glue together local solutions to the renormalised quantum master equation to find a global solution.

Recall that a simplicial presheaf $G$ on $M$ is a presheaf of simplicial sets on $M$; thus, for each open set $U \subset M$ and each integer $n$, we have a set $G(U,n)$ of $n$-simplices of the simplicial set $G(U)$. 

Let $\mbf{Met}$ be the simplicial presheaf such that $\mbf{Met}(U,n)$ is the set of smooth families of Riemannian metrics on $U$, parametrised by $\Delta^n$.  Let $\mbf{FMet} \subset \mbf{Met}$ be the sub-simplicial presheaf given by families of flat metrics. 

 It turns out that whether or not an action functional $S$ satisfies the renormalised quantum master equation is a local property.   Further,  the statement that $S$ satisfies the renormalised QME on an open set $U$ depends only on the metric $g$ on $U$. This is true also in families, parametrised by simplices.  Thus, we can arrange the solutions of the renormalised QME into a simplicial presheaf on $M$.
 
\begin{utheorem}
There is a simplicial presheaf $\mbf{BV}$ on $M$, with a map $\mbf{BV} \to \mbf{Met}$, whose  $0$-simplices $\mbf{BV}(U,0)$ are given by
\begin{enumerate}
\item
metrics $g$ on $U$ 
\item
a quantisation $S$ of the Chern-Simons action on $U$, modulo constants.  That is, $S$ is a solution of the renormalised quantum master equation (modulo constants) associated to $g$, which modulo $\hbar$ is the Chern-Simons action $S_0$.
\end{enumerate}
The one-simplices $\mbf{BV} (U,1)$ are homotopies of this data, and so on.
\end{utheorem}

Our ultimate aim is to construct a canonical (up to contractible choice) point of $\Gamma(M, \mbf{BV})$.   Such a point will consist of a metric $g$ on $M$  and a quantisation of the Chern-Simons action on $\E$ to a solution of the renormalised quantum master equation associated to $g$, as always modulo constants. 

We can always find a solution to the renormalised QME locally. 
\begin{uproposition}
Suppose an open subset $U\subset M$ is equipped with a flat Riemannian metric.  Then the original Chern-Simons action satisfies the renormalised quantum master equation.
\end{uproposition}
\begin{remark}
This proposition is the only result in this subsection which is really special to Chern-Simons theory.  The proof of this proposition relies heavily on the work of Kontsevich \cite{Kon94, Kon03a}.    In particular, we use the compactifications of configuration spaces used in these papers.  The quantum master equation is deduced from a theorem about the vanishing of certain integrals on configuration spaces proved by Kontsevich in \cite{Kon94} when $\dim M > 2$ and in \cite{Kon03a} when $\dim M = 2$.     
\end{remark}

The last proposition shows that we have a map 
$$
\mbf{FMet} \to \mbf{BV}
$$
of simplicial presheaves on $M$.

If $G$ is a simplicial presheaf on $M$, one can construct its derived global sections $\R \Gamma(M,G)$, which is a simplicial set.  We use a \v{C}ech definition of $\R \Gamma(M,G)$.  If $G_1 \to G_2$ is a map of simplicial presheaves which induces a  a weak equivalence of simplicial sets $G_1(U) \to G_2(U)$, for sufficiently small open balls $U$ in $M$, then the map $\R \Gamma(M,G_1) \to \R \Gamma(M,G_2)$ is a weak equivalence.  
\begin{ulemma}
For sufficiently small open balls $U$ in $M$, $\mbf{FMet}(U)$ is contractible.
\end{ulemma}
If $U$ is a ball, then $\mbf{FMet}(U)$ is the simplicial set associated to the space of flat metrics on $U$, which one can show is contractible using a simple rescaling argument.

It follows that $\R \Gamma(M,\mbf{FMet})$ is contractible. 
\begin{thmD}
The map
$$
\Gamma(M, \mbf{BV}) \to \R \Gamma(M,\mbf{BV})
$$
is a weak equivalence.
\end{thmD}
This theorem is the heart of the ``local-to-global'' principle; it allows one to glue together local solutions to the renormalised QME to give global ones.  This result is stated and proved for a general class of theories in the body of the paper.

We have a diagram
$$
\R \Gamma(M,\mbf{FMet}) \to \R \Gamma(M,\mbf{BV}) \xleftarrow{\simeq} \Gamma(M, \mbf{BV}) 
$$
where the first arrow comes from the map of simplicial presheaves $\mbf{FMet} \to \mbf{Met}$ constructed earlier.  The space on the left is contractible, and theorem D says that the arrow on the right is a weak equivalence.  Thus, we get the required point of $\Gamma(M, \mbf{BV}) $ up to contractible choice.

\subsection{Acknowledgements}
This paper benefited a great deal from conversations with many people.  I'd like to thank Dennis Sullivan for inviting me to present these results in his seminar, and for many helpful conversations.   Ezra Getzler's help with heat kernels and simplicial methods was invaluable, and made the paper much clearer.    A conversation with Jack Morava was very helpful early on.   Arthur Greenspoon's comments greatly improved the text.    I'm also grateful to Mohammed Abouzaid, John Baez, Alberto Cattaneo,  Paul Goerss, Dmitry Kaledin, Pavel Mnev, David Nadler, Vasily Pestun, Jared Wunsch and Eric Zaslow, for their help with various aspects of this work.

\section{A crash course in  the Batalin-Vilkovisky formalism}
\label{section intro bv}

This section should probably be skipped by experts; it consists of an informal introduction to the Batalin-Vilkovisky approach to quantising gauge theories.    The only thing which may not be standard is a discussion of a version of the BV formalism where one integrates over isotropic instead of Lagrangian subspaces.  

In this section, our vector spaces will \emph{always} be finite dimensional.  Of course, none of the difficulties of renormalisation are present in this simple case.  Many of the expressions we write in the finite dimensional case are ill-defined in the infinite dimensional case. 

Let us suppose we have a finite dimensional vector space $V$ of fields, a group $G$ acting on $V$ in a possibly non-linear way, and a $G$-invariant function $f$ on $V$ such that $0$ is a critical point of $f$.

One is interested in making sense of functional integrals of the form
$$\int_{ V / G } e^{f/ \hbar }$$
over the quotient space $V / G$.   The starting point in the Batalin-Vilkovisky formalism is the BRST construction, which says one should try to interpret this quotient in a homological fashion.  This means we should consider the  supermanifold
$$
\mf g [ 1 ] \oplus V
$$
where $[1]$ refers to a change of degree, so $\mf g$ is in degree $-1$. The space of functions on this super-manifold is 
$$
\Oo ( \mf g [ 1 ] \oplus V ) = \wedge^\ast \mf {g}^\vee \otimes \Oo (V)
$$
which is the super vector space underlying the Chevalley-Eilenberg Lie algebra co\-chain complex for $\mf g$ with coefficients in the $\mf g$-module $\Oo(V)$.  The Chevalley-Eilenberg differential gives an odd derivation of $\Oo(\mf g [1] \oplus V)$, which can be thought of as an odd vector field on $\mf g [1] \oplus V$.  Let us denote this odd vector field by $X$. 

 Recall that this Lie algebra cochain complex computes the homotopy invariants for the action of $\mf g$ on $\Oo(V)$.  Thus, we can view $\Oo(\mf g[1] \oplus V)$, with this differential, as a ``derived'' version of the algebra of functions on the quotient of $V$ by the action of $G$ 
(at least, in a formal neighbourhood of the origin, which is all we really care about).   The BRST construction says one should replace the integral over $V / G$ by an integral of the form
$$\int_{ \mf {g} [1 ] \oplus  V } e^{f / \hbar }.$$

This leaves us in a better situation than before, as we are in a linear space, and we can attempt to make sense of the integral perturbatively.   However,  we still have problems; the quadratic part of the functional $f$ is highly degenerate on $\mf {g} [1 ] \oplus  V$.  Indeed, $f$ is independent of $\mf g[1]$ and is constant on $G$-orbits on $V$.  Thus,  we cannot compute the integral above by a perturbation expansion around the critical points of the quadratic part of $f$. 

This is where the Batalin-Vilkovisky formalism comes in.    Let $E$ denote the odd cotangent bundle of $\mf{g} [1] \oplus V$, so that
$$
E = \mf{g} [1] \oplus V \oplus V^\vee [-1] \oplus \mf{g}^\vee [-2].
$$
The various summands of $E$ are usually given the following names: $\mf{g} [1]$ is the space of ghosts, $V$ is the space of fields,  $V^\vee [-1]$ is the space of antifields and $\mf{g}^\vee [-2]$ is the space of antighosts. 

The function $f$ on $\mf{g} [1] \oplus V$ pulls back to a function on $E$, via the projection $E \to \mf{g} [1] \oplus V$; we continue to call this function $f$.  By naturality, the vector field $X$ on $\mf{g} [1] \oplus V$ induces one on $E$, which we continue to call $X$.  As $[X, X] = 0$ on $\mf{g} [1] \oplus V$, the same identity holds on $E$.  As $X$ preserves $f$ on $\mf{g} [1] \oplus V$, it continues to preserve $f$ on $E$.

$E$ is an odd symplectic manifold, and $X$ is an odd vector field preserving the symplectic form.  Thus, there exists a unique function $h_X$ on $E$ whose Hamiltonian vector field is $X$, and which vanishes at zero.  As $X$ is odd, $h_X$ is an even function.

As $E$ is odd symplectic, the space of functions on $E$ has an odd Poisson bracket.  The statement that $[X,X] = 0$ translates into the equation $\{h_X, h_X\} = 0$. The statement that $X f = 0$ becomes $\{h_X, f\} = 0$.   And, as $f$ is pulled back from $\mf{g} [1] \oplus V$, it automatically satisfies $\{f,f\} = 0$.  These identities together tell us that the function $f + h_X$ satisfies the \emph{Batalin-Vilkovisky classical master equation}, 
$$
\{ f + h_X, f + h_X \} = 0 . 
$$

Let us write 
$$
f (e)  + h_X ( e ) = \frac{1}{2} \ip{e, Q e }  + S(e)
$$
where $Q : V \to V$ is an odd linear map, skew self adjoint for the pairing $\ip{\quad}$, and $S$ is a function which is at least cubic.  The fact that $f + h_X$ satisfies the classical master equation implies that $Q^2 = 0$.  Also, the identity 
$$
Q S + \frac{1}{2} \{ S , S \} = 0
$$ 
holds as a consequence of the classical master equation for $f + h _X$.

Let $L \subset E$ be a small, generic, Lagrangian perturbation of the zero section $\mf {g} [1] \oplus V \subset E$.  The Batalin-Vilkovisky formalism tells us to consider the functional integral
\begin{equation*}
\int_{e \in L} \exp  (  f(e) / \hbar + h_X (e)/ \hbar  )  = \int_{e \in L} \exp  ( \frac{1}{2 \hbar} \ip{e, Q e }  + \frac{1}{\hbar} S(e) ) 
\end{equation*}
As $L$ is generic, the pairing $\ip{e, Q e}$ will have very little degeneracy on $L$.  In fact, if the complex $(E,Q)$ has zero cohomology, then the pairing $\ip{e, Q e}$ is non-degenerate on a generic Lagrangian $L$.     This means we can perform the above integral perturbatively, around the critical point $0 \in L$. 

\subsection{Quantum master equation}
Let us now turn to a more general situation, where $E$ is a finite dimensional vector space with an odd symplectic pairing, and $Q : E \to E$ is an odd operator of square zero which is skew self adjoint for the pairing.   $E$ is not necessarily of the form constructed above.

Let $x_i, \eta_i$ be a Darboux basis for $E$, so that $x_i$ are even, $\eta_i$ are odd, and $\ip{x_i, \eta_i} = 1$.  Let $\Delta$ be the order two differential operator on $E$ given by the formula
$$
\Delta = \sum \partial_{x_i} \partial_{\eta_i}.
$$
This operator is independent of the choice of basis of $E$. 

Let $S \in \Oo(E ) [[\hbar]]$ be an $\hbar$-dependent function on $E$, which modulo $\hbar$ is at least cubic. The function $S$  satisfies the \emph{quantum master equation} if 
$$
(Q + \hbar \Delta) e^{S/  \hbar }  = 0 . 
$$
This equation is equivalent to the equation
$$
Q S +  \frac{1}{2} \{ S , S \} + \hbar \Delta S = 0 . 
$$

The key lemma in the Batalin-Vilkovisky formalism is the following.
\begin{lemma}
Let  $L \subset E$ be a Lagrangian on which the pairing $\ip{e,Q e }$ is non-degenerate. (Such a Lagrangian exists if and only if $H^\ast(E, Q ) = 0$).   Suppose that $S$ satisfies the quantum master equation. Then the integral
$$
\int_{e \in L} \exp  ( \frac{1}{2 \hbar} \ip{e, Q e }  + \frac{1}{\hbar} S(e) ) 
$$
is unchanged under deformations of $L$.
\end{lemma}
The non-degeneracy of the inner product on $L$, and the fact that $S$ is at least cubic modulo $\hbar$,   means that one can compute this integral perturbatively.

Suppose  $E, Q, \ip{\  , \ }, S$ are obtained as before from a gauge theory.  Then $S$ automatically satisfies the classical master equation $Q S + \frac{1}{2} \{S, S \} = 0$.  If, in addition, $\Delta S = 0$, then $S$ satisfies the quantum master equation.  Thus, we see that we can quantise the gauge theory in a way independent of the choice of $L$ as long as $S$ satisfies the equation $\Delta S = 0$. When $S$ does not satisfy this equation, one looks to replace $S$ by a series $S ' = S + \sum_{ i > 0} \hbar^i S_i$ which does satisfy the quantum master equation $Q S' +  \frac{1}{2} \{ S' , S' \} + \hbar \Delta S' = 0 $.

\subsection{Geometric interpretation of the quantum master equation}
The quantum master equation has a geometric interpretation, first described by Albert Schwarz \cite{Sch93}.  I will give a very brief summary; the reader should refer to this paper for more details.

Let $\mu$ denote the unique up to scale translation invariant ``measure'' on $E$, that is, section of the Berezinian.    The operator $\Delta$ can be interpreted as a kind of divergence associated to the measure $\mu$, as follows.  As $E$ is an odd symplectic manifold,  the algebra $\Oo(E)$ has an odd Poisson bracket.  Every function $S \in \Oo(E)$ has an associated vector field $X_S$, defined by the formula $X_S f = \{ S, f \}$.  

The operator $\Delta$ satisfies the identity
$$
\mc L_{X_S} \mu = (\Delta S) \mu
$$
where $\mc L_{X_S}$ refers to the Lie derivative.  In other words, $\Delta S$ is the infinitesimal change in volume associated to the vector field $X_S$.  

Thus, the two equations
\begin{align*}
\{S,S\} &= 0 \\
\Delta S &= 0
\end{align*}
say that the vector field $X_S$ has square zero and is measure preserving.

This gives an interpretation of the quantum master equation in the case when $S$ is independent of $\hbar$.  When $S$ depends on $\hbar$, the two terms of the quantum master equation do not necessarily hold independently.  In this situation, we can interpret the quantum master equation as follows.  Let $\mu_S$ be the measure on $E$ defined by the formula
$$
\mu_S = e^{S / \hbar} \mu.
$$
We can define an operator $\Delta_S$ on $\Oo(E)$ by the formula
$$
\mc L_{X_f} \mu_S = (\Delta_S f ) \mu_S.
$$
This is the divergence operator associated to the measure $\mu_S$, in the same way that $\Delta$ is the divergence operator associated to the translation invariant measure $\mu$.  

Then, a slightly weaker version of the quantum master equation is  equivalent to the statement
$$
\Delta_S^2 = 0.
$$
Indeed, one can compute that
$$
\hbar \Delta_S f = \{S, f \} + \hbar \Delta f 
$$
so that 
$$
\hbar^2 \Delta_S^2 f = \tfrac{1}{2} \{  \{S, S\} , f \} +  \hbar \{ \Delta S , f \} . 
$$
Thus, $\Delta_S^2 = 0$ if and only if $\tfrac{1}{2} \{S,S\} +\hbar \Delta S$ is in the centre of the Poisson bracket, that is, is constant. 

This discussion shows that the quantum master equation is the statement that the measure $e^{S / \hbar } \mu$ is compatible in a certain sense with the odd symplectic structure on $E$. 

\begin{remark}
In fact it is better to use half-densities rather than densities in this picture.  A solution of the quantum master equation is then given by a half-density which is compatible in a certain sense with the odd symplectic form. As all of our odd symplectic manifolds are linear, we can ignore this subtlety. 
\end{remark}

 \subsection{Integrating over isotropic subspaces}
 As we have described it, the BV formalism only has a chance to work when $H^\ast(E, Q ) = 0$.   This is because one cannot make sense of the relevant integrals perturbatively otherwise.   However, there is a generalisation of the BV formalism which works when $H^\ast(E,Q ) \neq 0$.  In this situation, let $L \subset E$ be an isotropic subspace such that $Q : L \to \Im Q$ is an isomorphism.    Let $\op{Ann}(L) \subset E$ be the set of vectors which pair to zero with any element of $L$. Then we can identify 
$$
H^\ast(E, Q) = \op{Ann}(L) \cap \Ker Q . 
$$
We thus have a direct sum decomposition 
$$
E = L \oplus H^\ast(E,Q) \oplus \Im Q . 
$$
Note that $H^\ast(E,Q)$ acquires an odd symplectic pairing from that on $E$. Thus, there is a BV operator $\Delta_{H^\ast(E,Q) }$ on functions on $H^\ast(E,Q)$.  We say a function $f$ on $H^\ast(E,Q)$ satisfies the quantum master equation if $\Delta e^{f / \hbar  } = 0$. 

The analog of the ``key lemma'' of the Batalin-Vilkovisky formalism is the following.  This lemma is well known to experts in the area.  
\begin{lemma}
Let $S \in \Oo(E) [[\hbar]] $ be an $\hbar$-dependent function on $E$ which satisfies the quantum master equation.  Then the function on $H^\ast(E,Q)$ defined by
$$
a \mapsto \hbar \log \left( \int_{e \in L} \exp  ( \frac{1}{2 \hbar} \ip{e, Q e }  + \frac{1}{\hbar} S(e+  a ) )  \right) 
$$
(where we think of $H^\ast(E,Q)$ as a subspace of $E$) satisfies the quantum master equation. 

Further, if we perturb the isotropic subspace $L$ a small amount, then this solution of the QME on $H^\ast( E, Q)$ is changed to a homotopic solution of the QME. 
\end{lemma}
Note that since $Q a = 0$, we can write the exponential in the integrand in the equivalent way $( \frac{1}{2 \hbar} \ip{e+a, Q (e+a) }  + \frac{1}{\hbar} S(e+  a ) )  $.  Thus, there is no real need here to separate out quadratic and higher terms. However, in the infinite dimensional situation we will discuss later, it will be essential to write the integrand as $( \frac{1}{2 \hbar} \ip{e, Q e }  + \frac{1}{\hbar} S(e+  a ) )  $, because we will take a field $a$ which is not closed.

This integral is an explicit way of writing the homological perturbation lemma for BV algebras, which transfers a solution of the quantum master equation at chain level to a corresponding solution on cohomology.    From this observation it's clear (at least philosophically) why the lemma should be true; the choice of the Lagrangian $L$ is essentially the same as the choice of symplectic homotopy equivalence between $E$ and its cohomology.  I'll omit a formal proof for now, as a proof of a more general statement will be given later. 

Let me explain what a ``homotopy'' of a solution of the QME is.  There is a general concept of homotopy equivalence of algebraic objects, which I learned from the work of Deligne, Griffiths, Morgan and Sullivan \cite{DelGriMor75}.  Two algebraic objects are homotopic if they are connected by a family of such objects parametrised by the commutative differential algebra $\Omega^\ast ([0,1])$.     In our context, this means that two solutions $f_0, f_1$  to the quantum master equation on $H^\ast(E,Q)$ are homotopic if there exists an element $F \in \Oo(H^\ast(E,Q) ) \otimes \Omega^\ast([0,1] ) [[\hbar]]$ which satisfies the quantum master equation
$$
 (  \d_{DR} + \hbar \Delta_{H^\ast(E,Q) } ) e^{F / \hbar } = 0
$$
and which restricts to $f_0$ and $f_1$ when we evaluate at $0$ and $1$.    Here $\d _{DR}$ refers to the de Rham differential on $\Omega^\ast([0,1])$. 

The quantum master equation imposed on $F$ is equivalent to the equation
$$
\d_{DR} F + \frac{1}{2} \{ F , F \} + \hbar \Delta F = 0 . 
$$
If we write $F(t, \d t) =   A(t ) + \d t B(t)$, then the QME imposed on $F$ becomes the system of equations
\begin{align*}
\frac{1}{2} \{ A(t) , A(t) \} + \hbar \Delta A (t) &= 0 \\
\frac{\d}  { \d t} A(t) + \{ A(t) , B(t) \} +  \hbar \Delta B(t) &= 0 . 
\end{align*}
The first equation says that $A(t)$ satisfies the ordinary QME for all $t$, and the second says that the family $A(t)$ is tangent at every point to an orbit of a certain ``gauge group'' acting on the space of solutions to the QME.  

\section{Example : Chern-Simons theory}
I want to discuss a class of quantum field theories in the Batalin-Vilkovisky formalism.  The general definition of the kind of quantum field theory I want to consider is a little technical, so I will start by discussing a simple example in detail, namely Chern-Simons theory on a $3$-manifold.  (Later we will discuss Chern-Simons theory on a manifold of any dimension).    Most of the features of more general theories are already evident in this example.  

Let $M$ be a compact oriented $3$-manifold.  Let $\mf g$ be a flat bundle of complex\footnote{Everything works if we take a real Lie algebra and work over $\R$ throughout.} Lie algebras with a complex valued invariant pairing  on $M$.  For example, $\mf g$ could be the adjoint bundle associated to a flat principal $SL(n,\C)$ bundle, equipped with the Killing form. 

We want to quantise the Chern-Simons gauge theory associated to $\mf g$. This is the theory whose space of fields $\mscr{V}$ is the space $\Omega^1(M, \mf g)$, which we think of as being the space of $\mf g$-valued connections.   The Lie algebra of the gauge group is
$$\mscr{G} = \Omega^0 (M, \mf g ).$$    This Lie algebra acts on the space of fields in the usual way that infinitesimal bundle automorphisms act on connections.  

The Chern-Simons action functional is the function on $\mscr{V}$ given by
$$
\frac{1}{2}\ip{v, \d v}  + \frac{1}{3} \ip{ v, [v,v] } 
$$
where $\d : \Omega^1 (\mf g) \to \Omega^2 (\mf g)$ is the operator obtained by coupling the flat connection on $\mf g$ to the de Rham differential.

Applying the Batalin-Vilkovisky construction, as described in section \ref{section intro bv}, yields an odd symplectic vector space 
\begin{align*}
\E = \mscr{G} [1] \oplus \mscr{V} \oplus \mscr{V}^\vee [-1] \oplus \mscr{G}^\vee [-2].
\end{align*}
If we interpret the duals of the infinite dimensional vector spaces appropriately, we find that 
$$
\E = \Omega^\ast (M, \mf g) [1] . 
$$
This is simply because $\Omega^2 (M,\mf{g})$ has a natural non-degenerate pairing with $ \mscr{V} = \Omega^1 (M,\mf{g} )$, and $\Omega^3 (M,\mf{g})$ has a natural non-degenerate pairing with $\Omega^0(M,\mf{g}) = \mscr{G}$.  

The differential $Q : \E \to \E$ (constructed as in section \ref{section intro bv}) is simply the de Rham differential, coupled to the flat connection on $\mf{g}$.     The odd pairing on $\E$ is given  by
$$
\ip{ \alpha \otimes X, \alpha' \otimes X'  } = (-1)^{\abs{\alpha} }\int_M \alpha\wedge \alpha' \ip{X, X' }_{\mf g}
$$
where $\alpha,\alpha' \in \Omega^\ast(M)$ and $X,X' \in \mf{g}$.  The notation $\ip{\ , \ }_{\mf g}$ refers to the invariant pairing on $\mf g$. 

The action functional $S$ is given by the formula
$$
S( \sum_i \alpha_i \otimes X_i ) = \sum_{i,j,k} (-1)^{(\abs{\alpha_i} + 1) (\abs{\alpha_j} + 1 )} \int_M \alpha_i \wedge \alpha_j \wedge \alpha_k \ip{ X_i, [X_j, X_k ]  }.
$$

\subsection{Gauge fixing and functional integrals}
As we saw in section \ref{section intro bv}, the next step in the Batalin-Vilkovisky formalism is the choice of a gauge fixing condition.  This is an isotropic subspace $L \subset \E$ such that $Q : L \to \Im Q$ is an isomorphism.   

Picking a metric on $M$ gives an operator $\d^\ast  : \Omega^i (M) \to \Omega^{i-1} (M)$.  This extends to an operator $\GF : \E \to \E$, defined (locally) by
$$
\GF ( \omega \otimes A)  = (\d^\ast \omega) \otimes A
$$
where $A$ is a flat section of the bundle $\mf g$.  Our gauge fixing condition will be
$$
L = \Im \GF . 
$$
The functional integral we will construct is
$$
 \int_{e \in \Im \GF} \exp  ( \frac{1}{2 \hbar} \ip{e, Q e }  + \frac{1}{\hbar} S(e+  a ) ) .
$$

Thus, the Chern-Simons gauge theory in the BV formalism is encoded in the following data.
\begin{enumerate}
\item
A vector space $\E$, which is the space of global sections of a super vector bundle on a compact manifold $M$.
\item
An odd anti-symmetric  pairing on $\E$, satisfying a certain non-degeneracy condition.
\item
A skew self adjoint operator $Q : \E \to \E$, of square zero.
\item
An even functional $S : \E \to \C$ which is ``local'', meaning roughly that it has a Taylor expansion in terms of continuous linear maps $\E^{\otimes k} \to \C$ which are distributions on $M^k$, supported on the small diagonal, and non-singular in the diagonal directions.   (A precise definition will be given soon).
\item
An auxiliary operator $\GF : \E \to \E$, which is a self adjoint elliptic operator of square zero.
\item
The super-commutator $[Q,\GF]$ must be an elliptic operator of order $2$ with some positivity conditions.
\end{enumerate}

\section{Batalin-Vilkovisky formalism in infinite dimensions}
\label{section infinite bv}

Now we will define the type of quantum field theory we will consider.  The definition consists essentially of abstracting the data we encountered above in the case of Chern-Simons theory.

\subsection{Functionals}
 If $M,N$ are smooth manifolds and $F,G$ are super vector bundles on $M,N$ respectively, we will use the notation
$$
\Gamma(M,F) \otimes \Gamma(N,G)
$$
to denote the space $\Gamma(M \times N, F \boxtimes G)$ of smooth sections of $F \boxtimes G$.  In other words, $\otimes$ refers to the completed projective tensor product where appropriate.  

Let $E$ be a super vector bundle over $\C$ on a compact manifold $M$.  Let $\E$ denote the space of global sections of $E$.   
\begin{definition}
Let $\Oo(\E)$ be the algebra 
$$
\Oo(\E) = \prod_{n \ge 0} \Hom ( \E^{\otimes n}, \C ) _{S_n}
$$
where $\Hom$ denotes the space of continuous linear maps, and the subscript $S_n$ denotes taking $S_n$ coinvariants.  

Direct product of distributions makes $\Oo(\E)$ into an algebra.
\end{definition}
We can view  $\Oo(\E)$ as the algebra  of formal functions at $0 \in \E$ which have Taylor expansions of the form
$$
f(e) = \sum f_i(e^{\otimes i} ) 
$$
where each 
$$
f_i : \E^{\otimes i} \to \C
$$
is a continuous linear map (i.e. a distribution).  
\begin{definition}
Let $X$ be an auxiliary manifold.  Define
$$
\Oo(\E, \cinfty(X)) = \prod_{n \ge 0} \Hom ( \E^{\otimes n}, \cinfty(X) ) _{S_n}
$$
where, as before, we take the space of continuous linear maps.  
\end{definition}
Thus, $\Oo(\E,\cinfty(X))$ is a certain completed tensor product of $\Oo(\E)$ and $\cinfty(X)$. 
Let $\hbar$ be a formal parameter.  Let 
$$
\Oo(\E, \C[[\hbar]] ) = \liminv \Oo(\E) \otimes \C[\hbar]/\hbar^n.
$$
Similarly, let $\Oo(\E, \cinfty(X) \otimes \C[[\hbar]])$ be the inverse limit $\liminv \Oo(\E, \cinfty(X) ) \otimes \C[\hbar] / \hbar^n$. 

\subsection{Local functionals}
\begin{definition}

Let $\op{Diff}(E,E')$ denote the infinite rank vector bundle on $M$ of differential operators between two vector bundles $E$ and $E'$ on $M$.  Let 
$$
\op{Poly Diff} ( \Gamma(E)^{\otimes n}, \Gamma(E') ) = \Gamma ( M, \op{Diff}(E, \C)^{\otimes n} \otimes E'  ) \subset \op{Hom} (\Gamma(E)^{\otimes n} , \Gamma(E'))
$$
where $\C$ denotes the trivial vector bundle of rank $1$.  All tensor products in this expression are fibrewise tensor products of vector bundles on $M$. 
\end{definition}
It is clear that $\op{Poly Diff}( \Gamma(E)^{\otimes n}, \Gamma(E'))$ is the space of sections of an infinite rank vector bundle on $M$.     If $F$ is a super vector bundle on another manifold $N$, let 
$$
\op{PolyDiff}  ( \Gamma(E)^{\otimes n}, \Gamma(E') ) \otimes \Gamma(F)
$$
denote the completed projective tensor product, as usual, so that 
$$
\op{PolyDiff}  ( \Gamma(E)^{\otimes n}, \Gamma(E') ) \otimes \Gamma(F)  =  \Gamma ( M \times N,  (\op{Diff}(E, \C)^{\otimes n} \otimes E') \boxtimes F  ) .
$$
If $X$ is a manifold, we can think of $\op{PolyDiff}  ( \Gamma(E)^{\otimes n}, \Gamma(E') ) \otimes \cinfty(X)$ as the space of smooth families of polydifferential operators parametrised by $X$. 

One can give an equivalent definition of polydifferential operators in terms of local trivialisations $\{e_i\}, \{e'_j\}$ of $E$ and $E'$,  and local coordinates $y_1,\ldots, y_l$ on $M$.  
A  map  $\Gamma(E)^{\otimes n} \to \Gamma(E')$ is a polydifferential operator if, locally, it is a finite sum of operators of the form 
$$
f_1 e_{i_1} \otimes \cdots \otimes f_n e_{i_n} \mapsto \sum _{j} e'_j  \Phi^{j}_{i_1 \ldots i_n} (y_1,\ldots,y_l)  ( D_{I_{i_1}}  f_1 ) \cdots ( D_{I_{i_n}} f_n ) 
$$
where $\Phi^{j}_{i_1 \ldots i_n} (y_1,\ldots,y_l)$ are smooth functions of the $y_i$, the $I_k$ are multi-indices, and the operators $D_{I_k}$ are the corresponding partial derivatives with respect to the $y_i$.

\begin{definition}
Let 
$$\Oo_l(\E) \subset \Oo(\E) =  \prod_{n \ge 0} \Hom ( \E^{\otimes n}, \C ) _{S_n}$$ be the space of formal functionals on $\E$, each of whose Taylor series terms $f_n : \E^{\otimes n} \to \C$ factors as 
$$
\E^{\otimes n} \to \dens(M) \xto{\int_M} \C
$$
where the first map is a polydifferential operator.    

Elements of $\Oo_l(\E)$ will be called \emph{local functionals} on $\E$.

If $X$ is another manifold (possibly with corners), and $F$ is a super vector bundle on $X$, let 
$$
\Oo_l( \E, \Gamma(X,F) ) \subset \Oo(\E, \Gamma(X,F)) = \prod_{n \ge 0} \Hom ( \E^{\otimes n}, \Gamma(X,F) ) _{S_n}
$$
be the subspace of those functions each of whose terms factors as 
$$
\E^{\otimes n} \to \dens(M) \otimes \Gamma(X,F) \xto{\int_M} \Gamma(X,F)
$$
where the first map is in $\op{Poly Diff} ( \E^{\otimes n}, \dens(M) ) \otimes \Gamma(X,F)$.
\end{definition}
Elements of $\Oo_l(\E, \cinfty(X))$ are smooth families of local functionals parametrised by a manifold $X$.    We will use a similar notation when we want to take functionals with values in $\C[[\hbar]]$.  Let
$$
\Oo_l(\E, \C[[\hbar]] ) = \liminv \Oo_l(\E) \otimes \C[\hbar]/\hbar^n.
$$
We can define $\Oo_l(\E, \cinfty(X) \otimes \C[[\hbar]])$ as an inverse limit in the same way. 

\begin{remark}
The space $\Oo_l(\E)$ is \emph{not} an algebra; the product of two local functionals is not local.
\end{remark}

If $f \in \Oo_l(\E)$ is a local functional, then the $i$-th term of its Taylor expansion $f_i : \E^{\otimes i} \to \C$ is assumed to factor through a map $\E^{\otimes i} \to \dens(M)$.
Note that there will in general exist many such factorisations. However, we have
\begin{lemma}
If 
$$\Phi : \E^{\otimes n} \to \C$$  is a map which factorises through a polydifferential operator $\E^{\otimes n} \to \dens(M)$,  there is a unique polydifferential operator
$$
\Psi : \E^{\otimes n-1} \to \Gamma(E^\vee) \otimes_{\cinfty(M)} \dens(M)
$$
such that 
$$
\Phi ( e_1 \otimes \cdots \otimes e_n ) = \int_M \ip{ \Psi ( e_1 \otimes \cdots \otimes e_{n-1} ) , e_n }.
$$
\label{lemma local functional vector field}
\end{lemma}
\begin{proof}
This is an easy calculation in local coordinates.
\end{proof}

\subsection{Batalin-Vilkovisky formalism}
\begin{definition}
An odd symplectic structure on $\E$ is a map
$$
\ip{\ , \ }_M : E \otimes E \to \dens_M
$$
of graded vector bundles on $M$, which is odd, antisymmetric, and non-degenerate, in the sense that the associated map $E \to E^{\vee} \otimes \dens_M$ is an isomorphism.  
\end{definition}
Such an odd symplectic structure  induces an odd pairing 
$$
\ip{\ , \ } : \E \otimes_\C \E \to \C
$$
on $\E$.

\begin{lemma}
An odd symplectic structure on $\E$ induces an odd Poisson bracket on the space $\Oo_l(\E)$ of local functionals.  Further, if $f \in \Oo_l(\E)$ is local and $g \in \Oo(\E)$ is possibly non-local, the Poisson bracket $\{f,g\}$ is well-defined.
\end{lemma}
\begin{proof}
This is an immediate corollary of Lemma \ref{lemma local functional vector field}.  Indeed, this lemma shows that any local functional on $\E$ can be replaced by something which has a Taylor expansion  in terms of polydifferential operators $\E^{\otimes n} \to \E$.   Here we are using the identification of $\E$ with $\Gamma(E^\vee )\otimes_{\cinfty(M)} \dens(M)$ given by the odd symplectic form.  

Something with a Taylor expansion in terms of maps $\E^{\otimes n} \to \E$ can be considered to be a formal vector field on $\E$.    This yields the Hamiltonian vector field on $\E$ associated to a local functional.  The usual formula for the action of vector fields on functions defines the Poisson bracket $\{f,g\}$ if $f \in \Oo_l(\E)$ and $g \in \Oo(\E)$.   Since this formula amounts to inserting the output of the map $\E^{\otimes n} \to \E$ into the input of a map $\E^{\otimes k} \to \C$, everything is well-defined.
\end{proof} 
In general, the Poisson bracket of two non-local functionals may be ill-defined.

Let us fix, for the rest of the paper, the following data:
\begin{enumerate}
\item
A compact manifold $M$ with a complex super vector bundle $E$ on $M$, whose algebra of global sections is denoted (as above) by $\E$.
\item
An odd linear elliptic differential operator $Q : \E \to \E$, which is self adjoint and satisfies $[Q,Q]  = 0$. $Q$ induces a differential on all the spaces associated to $\E$, such as $\Oo_l(\E)$ and $\Oo(\E)$.
\end{enumerate}
A  local functional $S \in \Oo_l(\E)$, which is at least cubic,  satisfies the classical master equation if
$$
Q S + \frac{1}{2}\{S,S\} = 0.
$$
We will try to quantise, within the Batalin-Vilkovisky formalism, field theories whose action functionals are of the form
$$
\frac{1}{2} \ip{e, Q e } + S(e)
$$
where $S \in \Oo_l(\E )$ satisfies the classical master equation. A wide class of quantum field theories of physical and mathematical interest, including for instance pure Yang-Mills theory, appear in this way.

\begin{definition}
Let $(M,\E,Q,\ip{\ , \ })$ be as above.  A gauge fixing condition is an odd linear differential operator 
$$
\GF :  \E \to \E
$$
such that
\begin{enumerate}
\item
$\GF$ is skew self adjoint with respect to the odd symplectic pairing on $\E$.
\item
$\GF$ is an operator of order $\le 1$. 
\item
$[\GF,\GF] = 0.$
\item
the operator 
$$
H := [Q,\GF]
$$
is a second order elliptic operator, which is a generalised Laplacian in the sense of \cite{BerGetVer92}. This means that the symbol 
$$
\sigma(H) \in \Gamma( M, \Sym^2 TM \otimes \End E) 
$$
is a positive definite metric on $T^\ast M$, tensored with the identity map $E\to E$.
\item
There is a direct sum decomposition
$$
\E = \Im Q \oplus \Im \GF \oplus \Ker H . 
$$
\end{enumerate}
\end{definition}
The fourth condition is a little restrictive. Even so, many interesting quantum field theories admit a natural collection of such gauge fixing conditions.  This condition is imposed because I need to use the asymptotic expansion of the heat kernel of $H$ proved in  \cite{BerGetVer92}.  

\section{Examples}
\label{section examples}
\subsection{$\phi^4$ theory}
The $\phi^4$ theory is a standard simple example of a quantum field theory. The Batalin-Vilkovisky formalism does not say anything interesting about the $\phi^4$ theory, because it is a simple bosonic field theory with no gauge symmetries.  

In our setting, the $\phi^4$ theory works as follows.  We take $M$ to be a compact Riemannian manifold.  Let 
$$
\E = \E^0 \oplus  \E^1 =  \cinfty(M,\C) \oplus \Pi \cinfty(M,\C)
$$ 
where $\Pi$ refers to parity reversal.   Let $e_0 \in \E^0$, $e_1 \in \E^1$ be the elements which correspond to the identity function on $M$.   The odd symplectic pairing on $\E$ is defined by 
$$
\ip{ \phi e_0, \psi e_1   } = \int_M \phi \psi 
$$
where $\phi, \psi \in \cinfty(M,\C)$ and the integration is taken using the measure associated to the Riemannian metric.

The operator $Q$ is defined by 
\begin{align*}
Q (\phi e_0 ) &= ( \Delta \phi + m^2 \phi ) e_1 \\
Q ( \phi e_1 ) &= 0 
\end{align*}
where $m \in \R_{> 0}$ is a ``mass'' parameter, and $\Delta$ is the Laplacian on $\cinfty(M)$ associated to the Riemannian metric.  

The action $S$ is defined by
$$
S ( \phi e_0 + \psi e_1 ) = \int_M \phi^4 .
$$
Then 
$$
\frac{1}{2} \ip{ \phi e_0 , Q \phi e_0  }  + S ( \phi ) = \int_M \frac{1}{2} \phi \Delta \phi + \frac{1}{2} m^2 \phi^2 + \phi^4 
$$
which is the usual action for the $\phi^4$ theory.

The gauge fixing operator $\GF$ on $\E$ is defined by
\begin{align*}
\GF (\phi e_1 ) &= \phi e_0 \\
\GF ( \phi e_0 ) &= 0 .
\end{align*}
The operator $H = [\GF, Q ]$ is given by
$$
H ( \phi e_i) = ( \Delta \phi + m^2 \phi  ) e_i
$$
for $i = 0,1$.   As required, this is a generalised Laplacian.  Also, 
$$\E = \Ker H \oplus \Im \GF \oplus \Im Q . $$  
Note that these axioms would not be satisfied if $m = 0$, because $\Ker H \cap \Im \GF$ would not be zero. 

The functional integral we will construct is 
$$
\int_{x \in \cinfty(M)} \exp \left(  \hbar^{-1}  \int_M \frac{1}{2} \phi \Delta \phi + \frac{1}{2} m^2 \phi^2 + \phi^4     \right).
$$
Of course, this construction works if we use any polynomial of $\phi$ which is at least cubic in place of $\phi^4$.  

Because the space of fields is concentrated in degrees $0$ and $1$, any local functional which is homogeneous of degree $0$ satisfies the renormalised quantum master equation.   

\subsection{Yang-Mills theory} 
We will use a certain first order formulation of Yang-Mills theory, well-known in the physics literature \cite{MarZen96, Wit04}.  This was studied in the Batalin-Vilkovisky formalism in \cite{Cos06a}.  I don't see how to fit the standard formulation of Yang-Mills theory into the framework of this paper; the problem is that I don't see how to construct gauge fixing conditions of the type we need.     I will refer to \cite{Cos06a} for details of this construction. 

Let $M$ be a compact oriented $4$ manifold, and let $V$ be a vector bundle on $M$ with a connection $A$ satisfying the Yang-Mills equation
$$
\d_A F(A)_+ = 0. 
$$
In \cite{Cos06a} a certain differential graded algebra $\mscr{B}$ was constructed from  $(M,V,A)$, which looks like
$$
\xymatrix{   \Omega^{0} (\End(V)) \ar[dd]_{\d_{A}} \ar[ddrr]^{-[F(A)_{+},\quad]} &  &  &  \text{ degree } 0  \\ 
 & & & \\
 \Omega^{1}(\End(V)) \ar[dd]_{\d_{A+}}   \ar'[dr]|{-[F(A)_{+}, \quad]} [ddrr] & \oplus &  \Omega^{2}_{+} (\End(V)) \ar[dd]^{- \d_{A}} \ar[ddll]|(.3){\op{Id}} &  \text{ degree } 1    \\
 & & & \\
  \Omega^{2}_{+}(\End(V))  \ar[ddrr]^{-[F(A)_{+},\quad]} & \oplus & \Omega^{3} (\End(V)) \ar[dd]^{-\d_{A}} &  \text{ degree } 2  \\
  & & & \\
   &  & \Omega^{4} (\End(V)) &  \text{ degree } 3  }
$$
The arrows describe the differential $Q$ in $\mscr{B}$.   The algebra $\mscr{B}$ also has a trace of degree $-3$, which comes from trace over $V$ and integration over $M$.  This trace induces an odd symmetric pairing $\Tr ab$ on $\mscr{B}$. 

Let $\E = \mscr{B}[-1]$ be $\mscr{B}$ shifted by one, with differential $Q$ inherited from $\mscr{B}$.  The pairing on $\mscr{B}$ induces a symplectic pairing on $\E$.   Define a functional $S$ on $\E$ by 
$$
S(e) = \Tr a^3.
$$
It was observed in section 3.3 of \cite{Cos06a} that $\E$, with this action functional,  is the Batalin-Vilkovisky odd symplectic space associated with the first order formulation of pure Yang-Mills gauge theory, as studied in for instance \cite{MarZen96,Wit04}.  Indeed, $\E^{-1}= \mscr{B}^0$ is the space of ghosts, $\E^0 = \mscr{B}^1$ is the space of fields, $\E^1 = \mscr{B}^2$ is the space of antifields, and $\E^2 = \mscr{B}^3$ is the space of antighosts.  

The choice of a metric on $M$ and a Hermitian metric on the vector bundle $V$ leads to a Hermitian metric on $\mscr{B}$.  The operator $\GF$ is the Hermitian adjoint to $Q$. This operator satisfies the axioms of a gauge fixing condition.   

\subsection{Holomorphic Chern-Simons}

Let $M$ be a Calabi-Yau manifold of odd complex dimension. Fix a holomorphic  volume form $\Omega_{Vol}$ on $M$.  Let $\mf g$ be a complex Lie algebra with a non-degenerate , symmetric, invariant pairing $\ip{\ , \ }$.   Let
$$
\E = \Omega^{0,\ast}(M) \otimes_\C \mf g. 
$$
$\E$ has an odd symplectic pairing, defined by
$$
\ip{\alpha \otimes X, \alpha' \otimes X' } = (-1)^{\abs{\alpha} }\int_M \Omega_{Vol} \wedge \alpha \wedge \alpha' \ip{X, X' } 
$$
where $\alpha \in \Omega^{0,\ast}(M)$ and $X \in \mf g$. 

Define $S \in \Oo_l(\E)$ by
$$
S( \sum_i \alpha_i \otimes X_i ) = \sum_{i,j,k} (-1)^{(\abs{\alpha_i} + 1) (\abs{\alpha_j} + 1 )} \int_M \Omega_{Vol} \wedge \alpha_i \wedge \alpha_j \wedge \alpha_k \ip{ X_i, [X_j, X_k ]  }.
$$
This is the holomorphic Chern-Simons action functional.

Picking a metric on $M$ leads to an operator $\dbar^\star$ on $\Omega^{0,\ast}$, and so on $\E$.  This is the gauge fixing operator.  The operator $H= [\dbar,\dbar^\star]$ is a generalised Laplacian.

The functional integral we will renormalise is
$$
\int_{x \in \Im \dbar^\star} \exp\left({\frac{1}{2\hbar} \ip{ x, \dbar x }  + \frac{1}{\hbar} S(x + a)}\right) .
$$

This theory can be generalised in several directions.  Firstly, we can take $\mf g$ to be a non-trivial holomorphic vector bundle on $M$ with a Lie algebra structure (i.e.\ a sheaf of Lie algebras over the sheaf of holomorphic functions on $M$).  In that case, to get a gauge fixing condition would require the choice of both a metric on $M$ and a Hermitian metric on the vector bundle.  

If $M$ is of even complex dimension, one can take $\mf g$ to be a super-Lie algebra with an odd  invariant pairing (see \cite{AleKonSch97}).  Of course, there is a version of this where we take a sheaf of such super-Lie algebras.

Yet another direction is to use $L_\infty$ rather than Lie algebras. In that case the action $S$ is no longer purely cubic, but has higher-order contributions from the higher $L_\infty$ operations.  Such theories are also discussed in \cite{AleKonSch97}.

\subsection{Higher dimensional Chern-Simons}
We have already discussed Chern-Simons theory in dimension three.  One can readily generalise this theory to manifolds of arbitrary dimension, just as with holomorphic Chern-Simons theory.  

For odd-dimensional manifolds, one  takes a complex Lie algebra $\mf g$ with an invariant pairing, and sets $\E = \Omega^\ast(M)\otimes \mf g [1]$. The Chern-Simons action is given by essentially the same formula as in the $3$-dimensional case:   
$$
S( \sum_i \alpha_i \otimes X_i ) = \sum_{i,j,k} (-1)^{(-1)^{(\abs{\alpha_i} + 1) (\abs{\alpha_j} + 1 )}} \int_M\alpha_i \wedge \alpha_j \wedge \alpha_k \ip{ X_i, [X_j, X_k ]  }
$$
where $\alpha_i \in \Omega^\ast(M)$ and $X_i \in \mf g$.  The operator $Q$ is defined to be the de Rham differential $\d$.  If we pick a metric on $M$ then we can define $\GF$ to be the adjoint operator $\d^\star$. 

As with holomorphic Chern-Simons theory, this can be generalised in various directions.  Firstly, we  could replace $\mf g$ by a locally trivial sheaf of Lie algebras.  Or, we could use a sheaf of $L_\infty$ algebras instead Lie algebras.

In the case when $M$ is even-dimensional, we can take $\mf g$ to be an $L_\infty$ algebra with an odd invariant pairing, or a locally trivial sheaf of such.

A particularly interesting case occurs when $\dim M = 2$.  Two dimensional Chern-Simons theory (introduced in \cite{AleKonSch97}) is essentially the same as the Poisson sigma model.   This theory is used in the proof of Kontsevich's formality theorem \cite{Kon03a,CatFel01}.  

If $\mf g$ is an $L_\infty$ algebra with an odd invariant pairing, then we can identify $\mf g^{even} = (\mf g^{odd})^{\vee} [1]$. Thus, if we let $V = \mf g^{odd} [1]$, so that $V$ is a purely even vector space, then 
$$
\mf g [1]  = \Pi T^\ast V . 
$$
The $L_\infty$ structure on $\mf g$ is described by a function $f$ on $\mf g[1] = \Pi T^\ast V$, satisfying the classical master equation $\{f,f\} = 0$. 

We can identify functions on $\Pi T^\ast V$ with polyvector fields on $V$, that is, 
$$
\Oo(\Pi T^\ast V) = \oplus_i \wedge^i T V\otimes \Oo(V). 
$$
The Poisson bracket on $\Oo(\Pi T^\ast V)$ corresponds to the Schouten bracket on polyvector fields. 

Thus, the $L_\infty$ structure on $\mf g$ defines a generalised Poisson structure on $V$.  If this $L_\infty$ structure is purely quadratic in the $\mf g^{ev}[1]$ variables, then we find a Poisson structure in the classical sense on $V$.

The space
$$
\E = \Omega^\ast( M ) \otimes \mf g [ 1]
$$
can be identified with the  supermanifold of maps $\Pi T M \to \Pi T^\ast V$, leading to the sigma model interpretation of the theory.

\section{Differential operators}
\label{section differential operators}

\subsection{Functor of points approach to formal super geometry}

In order to clarify sign issues, we will use the  ``functor of points'' approach to formal super geometry. A formal super-space is viewed as a functor  from the category of nilpotent local commutative super-algebras to the category of sets.    For example, the formal super-space $\mbb{A}^{m,n}$ sends $R \mapsto \left( R^{ev} \right)^{\oplus m}  \oplus  \left( R^{odd} \right)^{\oplus n}$.   

If 
$$
V = V^{ev} \oplus V^{odd}
$$
is a $\Z/2$-graded vector space,  the formal super-space $V$ sends $R \mapsto (V \otimes R)^{ev}$.    

An even formal function on $V$ is then a natural transformation of functors $V \to \mbb{A}^{1,0}$,  and an odd formal function on $V$ is a natural transformation to $V \to \mbb{A}^{0,1}$.   More explicitly, an even  formal function $f$ on $V$ is a collection of maps 
$$f_R  : (V \otimes R)^{ev} \to R^{ev},$$
natural with respect to morphisms $R \to R'$ of nilpotent local commutative super-algebras.   

Let us denote by $\Oo'(V)$ the $\Z/2$-graded algebra of all formal functions on $V$. We can identify this algebra with the inverse limit 
$$\Oo'(V) = \liminv \Sym^\ast V^\vee / \Sym^{\ge n} V^\vee. $$

In particular, if $\E$ as above is the global sections of some super vector bundle $E$ on $M$, we can think of $\E$ as an infinite dimensional formal super-space.  The algebra $\Oo(\E)$ is a subalgebra of the algebra $\Oo'(\E)$.  An element of $\Oo'(\E)$ is in $\Oo(\E)$ whenever its Taylor expansion is in terms of continuous linear maps on $\E^{\otimes n}$.  

In the case when $M$ is a point, so that $\E$ is simply a finite dimensional super vector space, then $\Oo(\E)$ and $\Oo'(\E)$ coincide.  

\subsection{Derivations}
The functor of points way of looking at things is useful (for instance) when thinking about derivations of the algebra $\Oo(\E)$.  Everything we do  could also be written more explicitly, but the advantage of the functor of points formalism is that it allows one to essentially forget about signs.  

Let $\epsilon$ be an odd parameter.  One can define the derivation $Q$ of $\Oo(\E)$ by the formula
$$
f ( x + \epsilon Q x  ) = f( x ) + \epsilon (Q f ) (x) . 
$$
The left hand side is defined using the language of the functor of points, as follows.  For each  nilpotent local commutative superalgebra $R$ as above, $f$ gives a map 
$$f_{R \otimes \C[\epsilon]} :  (\E \otimes R \otimes \C[\epsilon])^{ev} \to
\begin{cases}
 (R \otimes \C[\epsilon] )^{ev} \text{ if } f \text{ is even} \\ 
  (R \otimes \C[\epsilon] )^{odd} \text{ if } f \text{ is odd} . 
\end{cases}
$$
If $x \in (\E \otimes R)^{ev}$, then $x + \eps Q x$ is in $( \E \otimes R \otimes \C[\eps] )^{ev}$.  Thus,  $f ( x + \eps Q x ) \in (R \otimes \C[\eps])$.  The coefficient of $\epsilon$ in $f( x + \epsilon Q x ) $ gives a map $(\E \otimes R)^{ev}$ to $R^{ev}$ (if $f$ is odd), or to $R^{odd}$ (if $f$ is even), and so a formal functional on $\E$ of parity opposite to that of $f$.    This is defined to be $Q f$.

A priori, this construction only yields a derivation of $\Oo'(\E)$, but one can check easily that this derivation preserves the subalgebra $\Oo(\E)$.   In more explicit terms, if we identify
$$
\Oo(\E) = \prod \Hom ( \E^{\otimes n}, \C ) _{S_n} ,
$$
then the derivation $Q$ preserves each factor $\Hom(\E^{\otimes n} , \C)_{S_n}$.  On each such factor, $Q$ is simply the usual tensor product differential on each $\E^{\otimes n}$, which we then transfer to the dual space and to the space of $S_n$ coinvariants.

 To compute the commutator of derivations in terms of the functor of points is very simple.  Let  $D_1,D_2$ be derivations of $\Oo'(V)$, where $D_i$ is given by an automorphism of $\Oo'(V) \otimes \C[\epsilon_i] / (\epsilon_i^2)$, which modulo $\epsilon_i$ is the identity. The parameter $\epsilon_i$ has the same parity as $D_i$.   We can compute the commutator by 
$$
(1 + \epsilon_1 D_1 ) ( 1 + \epsilon_2 D_2 ) ( 1 - \epsilon_1 D_1 )(1 - \epsilon_2 D_2 )  = 1 +  \epsilon_1 \epsilon_2 [D_1,D_2] .
$$

Let $e \in \E$.  Differentiation with respect to $e$ is an operator $\partial_e  : \Oo'(\E) \to \Oo'(\E)$ defined by the formula 
$$
f ( x + \epsilon e )  = f ( x ) + \epsilon \partial_e f 
$$
where $\epsilon$ is an auxiliary parameter of the same parity as $e$, and of square $0$.  

It is easy to see that this derivation preserves $\Oo(\E)$.  Indeed, if we identify $\Oo(\E) = \prod \Hom ( \E^{\otimes n}, \C ) _{S_n}$, then the derivation $\partial_e$ is (up to sign) the map 
$$ \Hom(\E^{\otimes n}, \C)_{S_n} \to \Hom(\E^{\otimes n-1}, \C)_{S_{n-1}}$$ 
given by contraction with $e$.

One can compute easily that 
$$
[ \partial_e , Q] = \partial_{Q e} .
$$

\subsection{Convolution operators}

If $K \in \E \otimes \E$, define a convolution map $K \star : \E \to \E$ by the formula
$$
K \star e = (-1)^{\abs e} \sum K' \otimes \ip{K'', e}
$$
where $K = \sum K' \otimes K''$.    The reason for the choice of sign in the definition of the convolution is that 
$$
( Q  K ) \star e  =   [Q, K\star ]  e.
$$
The first $Q$ in this formula denotes the tensor product differential on $\E \otimes \E$.

\begin{lemma}
An element $K \in \E \otimes \E$ is symmetric if and only if the map $K \star $ is self adjoint, and is antisymmetric if and only if $K \star $ is skew self adjoint. 
\end{lemma}

\subsection{Order two differential operators}
Let
$$
\phi = \sum \phi' \otimes \phi'' \in \Sym^2 \E.
$$
Associated to $\phi$ we have an order two differential operator $\partial_{\phi}$ on $\Oo(\E)$, given by the formula
$$
\partial_{\phi} =  \frac{1}{2} \sum \partial_{\phi''} \partial_{\phi'} =  \frac{1}{2} \sum  (-1)^{\abs{\phi'} \abs{\phi'' } }\partial_{\phi'} \partial_{\phi''}.
$$
Even though the sum may be infinite, this expression is well defined; up to sign, $\partial_\phi$ comes from the map 
$$ \Hom(\E^{\otimes n}, \C)_{S_n} \to \Hom(\E^{\otimes n-2}, \C)_{S_{n-2}}$$ 
given by contraction with $\phi$.  

The differential operators $\partial_\phi$ mutually super-commute, and 
$$
[\partial_\phi, Q] = \partial_{Q \phi} . 
$$

If $e \in \E$, let $e^\vee : \E \to \C$ be the continuous linear functional given by $e^\vee( e' ) = \ip{e',e}$.   Note that $e^\vee \in \Oo(\E)$, so there is an associated left multiplication map $\Oo(\E) \to \Oo(\E)$.
The following lemma will be useful later.
\begin{lemma}
If $\phi \in \Sym^2 \E$ is even, then
$$
[\partial_\phi, e^\vee ] =   -\partial_{\phi \star e } .
$$
\end{lemma}
\begin{proof}
In the language we are using for formal functions, if $R$ is an auxiliary nilpotent ring, and $r \in R$, $\eta \in \E$ are elements of the same parity, then by definition
$$
e^\vee (r \eta ) = \ip{r \eta , e } = r \ip{\eta , e} .
$$
Therefore, in particular, 
$$
\partial_{e'} e^\vee = \ip{ e', e } 
$$
for all $e' \in \E$.    

We can assume
$$
\phi = \phi' \otimes \phi'' + (-1)^{\abs{\phi'} \abs{\phi''}  } \phi'' \otimes \phi'.
$$
Then
$$
\partial_\phi = \partial_{\phi''} \partial_{\phi'}.
$$
Then
$$
[\partial_\phi, e^\vee ]  = \ip{\phi',e}   \partial_{\phi''}  + (-1)^{\abs{\phi'} (\abs{e}+1) } \ip{\phi'',e}   \partial_{\phi''} 
$$
whereas
$$
\phi \star e = (-1)^{\abs{e} } \ip{\phi'', e}  \phi' + (-1)^{\abs{e} + \abs{\phi''} \abs{\phi' }} \ip{\phi'', e}  \phi'.
$$
The lemma follows from these two equations. 
\end{proof}

\subsection{Heat kernels}
A heat kernel for the operator $H : \E \to \E$ is an element 
$$K_t \in\E^{\otimes 2} = \Gamma(M^2,E \boxtimes E),$$ defined for $t \in \R_{> 0}$, such that 
$$
K_t \star e = e^{- t H} e 
$$
for all $e \in \E$.  Because $H$ is a generalised Laplacian, the results of \cite{BerGetVer92} imply that it admits a unique heat kernel $K_t$, which is also smooth as a function of $t$.

Let 
$$L_t = (\GF \otimes 1) K_t$$ 
so that $L_t$ is the kernel representing the operator $\GF e^{-t H}$. This is a smoothing operator for $0< t \le \infty$, so that $L_t \in \E \otimes \E$ for such $t$.

Observe that 
\begin{align*}
(Q \otimes 1 + 1 \otimes Q) K_t &= 0 \\
\frac { \d } { \d t } K_t  + ( Q \otimes 1 + 1 \otimes Q) L_t &= 0 .
\end{align*}
These formulae together say that the expression
$$
K_t + \d t L_t 
$$
is a closed element of $\E \otimes \E \otimes \Omega^\ast(\R_{> 0})$, where we give this space the tensor product differential.  

Let
$$
P(\eps,T)= \int_{\epsilon}^T L_t \d t = \int_{\eps}^{T} (\GF \otimes 1)K_t \d t .
$$
Note that for $0 < \epsilon < T \le \infty$, $P(\epsilon,T)$ is in $\E \otimes \E$.    This is clear if $0 < \epsilon < T < \infty$.  If $T = \infty$, one needs to check that 
$$
\alpha \mapsto \int_\epsilon^\infty \GF e^{-t H} \alpha \d t
$$
is a smoothing operator for any $t > 0$. The only problems that could occur would be on the $0$ eigenspace of $H$, but $\GF$ annihilates this eigenspace.  This is one of the axioms of a gauge fixing condition. 

The reason the kernels $P(\epsilon,T)$ are important comes from the following lemma.
\begin{lemma}
The operator
\begin{align*}
\Im Q &\to \Im \GF \\
\alpha &\mapsto \int_0^\infty \GF e^{-t H} \alpha \d t
\end{align*}
is the inverse to the isomorphism $Q : \Im \GF \to \Im Q$.
\end{lemma}
Thus, the singular kernel $P(0,\infty)$ represents $Q^{-1}$.
\begin{proof}
Indeed,
$$
\int_0^\infty \GF e^{-t H} Q \alpha \d t= - Q \int_0^\infty \GF e^{-t H} \alpha \d t + \int_0^\infty H e^{-t H} \alpha \d t . 
$$
Since $(\GF)^2 = 0$, and $\alpha \in \Im \GF$, the first term on the right hand side is zero.  Thus, 
$$
\int_0^\infty \GF e^{-t H} Q \alpha \d t =  - \int_0^\infty \frac{\d } {\d t }  e^{-t H} \alpha \d t  = \alpha - \pi \alpha
$$ 
where $\pi : \E \to \Ker H$ is the projection onto the zero eigenspace of $H$.  

Since $\Ker H \cap \Im \GF = 0$, and $\Im \GF$ is a direct sum of $H$ eigenspaces, we see that $\pi \alpha = 0$, so that
$$
\int_0^\infty \GF e^{-t H} Q \alpha \d t = \alpha 
$$
as desired.  
\end{proof} 

\subsection{Functional integrals in terms of differential operators}
When $M$ is a point, we will be able to express certain integrals over the finite dimensional vector space $\E$ in terms of the differential operators on $\Oo(\E)$ introduced above.  When $\dim M > 0$, so that $\E$ is infinite dimensional, we will attempt to take this as a definition of the functional integral.

\begin{definition}
Let  $\Phi \in \Sym^2 \E$ be any even element, and  $S \in \Oo(\E, \C [[\hbar]])$ be an even element, which modulo $\hbar$ is at least cubic.  Define
$$
\Gamma( \Phi, S ) = \hbar \log\left( \exp  (\hbar \partial_\Phi ) \exp (  S / \hbar )   \right) \in \Oo(\E, \C [[\hbar ]] ) . 
$$
\end{definition}
It is easy to check that this expression is well defined. However, if $S$ was not at least cubic, but contained a quadratic term, this expression would contain some non-convergent infinite sums. 

Let us assume for a moment that $M$ is a point.
Let $\left(\Im \GF\right)_\R$ be a real slice of $\Im \GF \subset \E$, with the property that the quadratic form $\ip{x,Qx}$ is negative definite on $\left(\Im \GF\right)_\R$.  Then the integral
$$
\int_{x \in \left(\Im \GF\right)_\R}  \exp\left( \tfrac{1}{2} \ip{x,Q x} / \hbar + S(x + a) / \hbar \right) \d \mu
$$
is well defined as a formal series in the variables $\hbar$ and  $a \in \E$.  We use an $\hbar$-dependent Lebesgue measure $\d \mu$ normalised so that 
$$
\int_{x \in \left(\Im \GF\right)_\R}  \exp\left( \tfrac{1}{2} \ip{x,Q x}   / \hbar \right) \d \mu = 1 . 
$$
\begin{lemma} If $M$ is a point,
\label{lemma integral differential}
$$  \Gamma(P(0,\infty), S) (a)  =  \hbar \log \int_{x \in \left(\Im \GF\right)_\R}  \exp\left( \tfrac{1}{2} \ip{x,Q x} / \hbar + S(x + a) / \hbar \right) \d \mu .  $$
\label{lemma integral differential general}
\end{lemma}
\begin{proof}
It suffices to show that for any polynomial  $f \in \Oo(\E)$
$$
 \exp ( \hbar \partial_{P(0,\infty)} ) f =  \int_{\left(\Im \GF\right)_\R}  \exp \left( \tfrac{1}{2} \ip{x,Q x} /\hbar \right)  f (x+a).
$$
Both sides are functions of $a \in \E$.  

Let $v \in \left(\Im \GF\right)_\R \subset \E$. Define a linear function $l_v = (Q v)^\vee$ on $\left(\Im \GF\right)_\R$, so that
$$
l_v(x)  = \ip{x,Q v}.
$$
Then
$$
\partial_v \ip{x, Q x} = 2l_v (x). 
$$
We can integrate by parts, to find
\begin{multline*}
\int_{\left(\Im \GF\right)_\R}  \exp \left( \tfrac{1}{2} \ip{x,Q x} /\hbar \right) l_v (x + a) f(x+a) \\  =   \hbar \int_{\left(\Im \GF\right)_\R} \left(\partial_v  \exp \left( \tfrac{1}{2} \ip{x,Q x} /\hbar \right)  \right) f (x+a) 
\\ \shoveright{ +  l_v(a) \int_{\left(\Im \GF\right)_\R}  \exp \left( \tfrac{1}{2} \ip{x,Q x} /\hbar \right)   f (x+a ))  } \\
= - \hbar \int_{\left(\Im \GF\right)_\R}  \exp \left( \tfrac{1}{2} \ip{x,Q x} /\hbar \right) \partial_v f (x+a)\\
 +   l_v(a) \int_{\left(\Im \GF\right)_\R}  \exp \left( \tfrac{1}{2} \ip{x,Q x} /\hbar \right)   f (x+a ))  
\end{multline*}
A similar identity holds for $\exp(\hbar \partial_{P(0,\infty)}) f$,  namely
$$
 \exp( \hbar \partial_{P(0,\infty)}) l_v  f = - \hbar   \exp (\hbar \partial_{P(0,\infty)}) \partial_v f + l_v(a)  \exp( \hbar \partial_{P(0,\infty)})  f .
$$
This follows from the equation
$$
[\partial_{P(0,\infty)}, l_v ] = [\partial_{P(0,\infty)}, (Q v)^\vee ] =   - \partial_{ P(0,\infty) \star (Q v)  }  = - \partial_{ v } .
$$
The last equation holds because $P(0,\infty) \star$ is the operator $Q^{-1} : \Im Q \to \Im \GF$.  

These identities allow us to use induction to reduce to the case when $f$ is constant.  The normalisation in the measure on  $\left(\Im \GF\right)_\R$ takes care of this initial case.

\end{proof} 

\section{Regularisation}

\subsection{Regularisation}

One can write the formal equality
$$
 \hbar \log \int_{e \in (\Im \GF)_\R} \exp(  \tfrac{1}{2} \ip{e,Qe}  / \hbar +   S ( e + a) / \hbar  ) =  \lim_{ \epsilon\to 0}  \Gamma  (P(\epsilon,\infty), S)
$$
where as above,
$$
 \Gamma  (P(\epsilon,\infty), S) = \hbar \log \left( \exp (\hbar \partial_{P(\epsilon,\infty)}) 
  \exp ( S / \hbar) \right)  . 
$$
When $\dim M = 0$, this equality was proved in Lemma \ref{lemma integral differential}.  When $\dim M > 0$, we will take this equality as an attempt to define the functional integral over $(\Im \GF)_\R$.  

When $\dim M > 0$,  although the expression $\Gamma(P(\epsilon,\infty), S)$ is well defined for all $\epsilon > 0$, the limit  $\lim_{\epsilon \to 0} \Gamma(P(\epsilon,\infty), S)$ is divergent.   This is because $P(0,\infty)$ is a distributional section of the vector bundle $E \boxtimes E$ on $M^2$,  with singularities on the diagonal.  The expression $\exp (\hbar \partial_{P(0,\infty)} )  \exp ( S / \hbar)$ involves multiplication of the distribution $P(0,\infty)$ with distributions with support along the diagonal.    In other words, it involves integrals over products of $M$ where the integrand has singularities on the diagonal.      This is the problem of ultraviolet singularities.

As I explained in the introduction, we will renormalise the limit by subtracting certain counter-terms from the action $S$.  In order to do this, we need some control over the small $\epsilon$ asymptotics of $\Gamma(P(\eps,\infty), S)$. 

Let us write
$$
\Gamma(P(\epsilon,T),S) = \sum_{i\ge 0,k \ge 0} \hbar^i \Gamma_{i,k} (P(\epsilon,T),S)
$$
where $\Gamma_{i,k}(P(\eps,T),S)$ is homogeneous of order $k$ as a formal functional of $e \in \E$.   This expression is the Taylor expansion of $\Gamma$ in the variables $\hbar$ and $e \in \E$. 

We will construct a small $\epsilon$ asymptotic expansion for each $\Gamma_{i,k}(P(\epsilon,T), S)$.   This expansion will take values in a certain subalgebra $\mscr{A}$ of the algebra of analytic functions on $\epsilon \in (0,\infty)$.   
\begin{definition}
Let $\mscr{A} \subset \cinfty( (0,\infty) )$ be the subalgebra spanned over $\C$ by  functions of $\epsilon$ of the form
$$
f(\epsilon )  = \int_{U \subset (\epsilon,1)^n } \frac{  F(t_1,\ldots,t_n)^{1/2} } { G(t_1,\cdots, t_n)^{1/2} } \d t_1 \ldots \d t_n 
$$
and functions of the form
$$
f(\epsilon )  = \int_{U \subset (\epsilon,1)^{n-1} } \frac{  F(t_1,\ldots,t_n = \eps)^{1/2} } { G(t_1,\ldots, t_n = \eps)^{1/2} } \d t_1 \cdots \d t_{n-1}
$$
where 
\begin{enumerate}
\item
$F, G \in \Z_{\ge 0} [t_1,\ldots, t_n] \setminus \{0\}$; $n$ can take on any value. 
\item
the region of integration $U$ is an open subset cut out by finitely many equations of the form $t_i^l > t_j$, for $l \in \Z$.  
\end{enumerate}
\end{definition}
``Spanned'' means in the non-topological sense; we only allow finite sums.  Thus, $\mscr{A}$ has a countable basis, and every element is written as a finite sum of basis elements. We give $\mscr{A}$ the trivial topology, where all subspaces are closed. 

The details in the definition of $\mscr{A}$ aren't all that important; we could always use a larger algebra containing $\mscr{A}$. 

\begin{thmA}
\begin{enumerate}
\item
For each $i$ and $k$, there exist functions $f_r \in \mscr{A}$, $r \in \Z_{\ge 0}$, such that there is a small $\eps$ asymptotic expansion of the form 
$$
\Gamma_{i,k}(P(\epsilon,T), S(\epsilon, \hbar) ) (e) \simeq  \sum_{r \ge 0} f_r (\epsilon) \Phi_r (T,e) 
$$
where $e \in \E$, and each $\Phi_r$ is in $\Oo(\E, \cinfty((0,\infty)))$. 
\item
The $\Phi_r(T,e)$ have a small $T$ asymptotic expansion
$$
\Phi_r(T,e) = \sum g_q(T) \Psi_{q,r} (e)
$$ where the $\Psi_{q,r} \in \Oo_l(\E)$ are local functionals of $e$, and $g_q(T)$ are certain smooth functions of $T$. 
\item
If $k > 0$, so $\Psi_{q,r}$ is a non-constant function on $\E$,  then we can speak of the germ of each $\Psi_{q,r}$ near a point $x \in M$.  This germ only depends on the germ of the data $(E, \ip{\ , \ }, Q, \GF, S )$ near $x$.   
\end{enumerate}
\label{theorem_analytic_continuation} 
\end{thmA}

\begin{remark}
The last point is a little delicate.  For instance, for any $t > 0$, the germ of the heat kernel $K_t$ near a point $(x,x)$ in $M^2$ depends on the \emph{global} behaviour of the elliptic operator $H$.  Only the $t \to 0$ asymptotics of the germ of $K_t$ depends on the germ of $H$.
\end{remark}
The appendix contains a proof of a more general version of this theorem, as well 
as a precise statement of what I mean by small $\eps$ asymptotic expansion.

\section{Renormalisation}
The next step in constructing the quantum field theory is \emph{renormalisation}.  This amounts to replacing our action functional $S \in \Oo_l(\E, \C[[\hbar]])$ by a series
$$
S^R(\hbar,\epsilon) = S(\hbar) - S^{CT}(\hbar,\epsilon) = S - \sum_{i > 0, k \ge 0} \hbar^i S^{CT}_{i,k}(\epsilon).
$$
The $S_{i,k}^{CT}$ are known as counter-terms; they are $\eps$-dependent local functionals on $\E$, homogeneous of degree $k$.   This renormalised action functional $S^R$ is required to be such that the $\epsilon \to 0$ limit of $\Gamma( P(\eps,T), S^R (\eps) )$  exists.  Recall that 
$$
P(\eps,T)= \int_{\eps}^T  (\GF \otimes 1) K_t \d t
$$
so that $\lim_{\eps \to 0}\Gamma(P(\eps,T), S^R(\eps) )$ will be our renormalised effective action.

In order to perform this renormalisation, we need to pick a renormalisation scheme. Let $\mscr{A}_{\ge 0} \subset \mscr{A}$ be the subspace of those functions $f$ such that $\lim_{\eps \to 0} f$ exists. 
\begin{definition}
A \emph{renormalisation scheme} is a subspace $\mscr{A}_{< 0}$ to $\mscr{A}_{\ge 0}$ such that 
$$
\mscr{A} = \mscr{A}_{< 0} \oplus \mscr{A}_{\ge 0} . 
$$
\end{definition}
Note that a renormalisation scheme exists, simply because $\mscr{A}$ is a vector space with a countable basis.

The choice of a renormalisation scheme allows us to extract the singular part of functions $f \in \mscr{A}$.  The singular part of $f$ is simply the projection of $f$ onto $\mscr{A}_{< 0}$.

\begin{thmB}
Let $S \in \Oo_l(\E, \C[[\hbar]])$ be an even functional which, modulo $\hbar$, is at least cubic.   Then there exists a unique series
$$S^{CT}(\eps,\hbar) =   \sum_{i \ge 1, k \ge 0} \hbar^i S_{i,k}^{CT}(s)  $$
such that 
\begin{enumerate}
\item
each $S_{i,k}^{CT}(\eps)$ is an element of $\Oo_l(\E, \mscr{A}_{< 0})$
which is homogeneous of degree $k$ as a function on $\E$; 
\item
the limit $\lim_{\epsilon \to 0} \Gamma(P(\eps,T),S - S^{CT})$ exists.  
\end{enumerate}

If $k > 0$, the germ of the counter-term $S_{i,k}^{CT}$ at $x \in M$ depends only on the germ of the data  $(E, \ip{\ , \ }, Q, \GF,  S )$ near $x$.  
\label{theorem renormalisation}
\end{thmB}
We use the notation 
$$\Gamma^R(P(0,T), S) =  \lim_{\eps \to 0} \Gamma(P(\eps,T),S - S^{CT}). $$ $\Gamma^R(P(0,\infty), S) $ is a renormalised functional integral.

\begin{proof}

Before I begin the proof, let us introduce some notation.   Let us give the set $\Z_{\ge 0} \times \Z_{\ge 0}$ the lexicographic ordering, so that $(i,k) < (i',k')$ if $i < i'$, or if $i = i'$ and $k < k'$. 

Let $A \in \Oo(\E, \C[[\hbar]])$ be any (not necessarily local) functional.    As before, we will write 
$$
A = \sum \hbar^i A_{i,k} 
$$
where $A_{i,k} \in \Oo(\E)$ is a homogeneous functional of degree $k$.  Let us write 
$$
A_{\le (I,K)} = \sum_{(i,k) \le (I,K) } \hbar^i A_{i , k}.
$$

The proof of this theorem is actually very simple; no graph combinatorics are required.  All we need to do is to construct the $S_{i,k}^{CT}$ inductively using the lexicographic ordering on $\Z_{\ge 0} \times \Z_{\ge 0}$.

So, let us suppose, by induction, that we have constructed 
$$S_{i,k}^{CT}(\eps, T) \in \Oo_l(\E, \mscr{A}_{< 0} \otimes \cinfty(\R_{> 0}))$$ 
for all $(i,k) < (I,K)$,  which are homogeneous of degree $k$ as a function of $\alpha \in \E$.   Here, $\R_{> 0}$ has coordinate $T$.     

Let 
$$
S^{CT}_{<(I,K)}  = \sum_{(i,k) < (I,K)} \hbar^i S^{CT}_{i,k}
$$
We can write
$$
\Gamma ( P(\eps,T),  S - S^{CT}_{<(I,K)} )= \sum_{i,k \ge 0} \hbar^i  \Gamma_{i,k}( P(\eps,T), S - S^{CT}_{<(I,K)} )
$$
Let us make the following further induction assumptions : 
\begin{enumerate}
\item
the limit  $\lim_{\eps \to 0} \Gamma_{i,k}( P(\eps,T), S - S^{CT}_{<(I,K)} )$ exists, for all $(i,k) < (I,K)$. 
\item
each $S^{CT}_{i,k}$ is independent of $T$ if $(i,k) <(I,K)$. 
\end{enumerate}

Now simply let
$$
S^{CT}_{I,K} = \text{Singular part of } \Gamma_{I,K}( P(\eps,T), S - S^{CT}_{<(I,K)} ).
$$
By ``singular part'' I mean the following.  We take the small $\eps$ expansion of the right hand side, of the form $\sum f_r(\eps) \Phi( T, a)$ where $f_r(\eps) \in \mscr{A}$.  The singular part is obtained by replacing each $f_r$ by its projection onto $\mscr{A}_{< 0}$. 

Note that 
$$
 \Gamma_{i,k}( P(\eps,T), S -  S^{CT}_{<(I,K)} - S^{CT}_{I,K})  =  \Gamma_{i,k}( P(\eps,T), S -  S^{CT}_{<(I,K)} )  - \delta_{i,I} \delta_{k,K} S^{CT}_{I,K} 
$$
if $(i,k) \le (I,K)$.   It follows that  $\Gamma_{i,k}( P(\eps,T), S - S^{CT}_{<(I,K)} - S^{CT}_{I,K}) $ is non-singular (i.e.\ has a well-defined $\eps \to 0$ limit)  if $(i,k) \le (I,K)$. 

To prove the result, it remains to prove the following. 
\begin{enumerate}
\item
$S^{CT}_{I,K}$ is independent of $T$.
\item
$S^{CT}_{I,K}$ is a local functional.  
\item
$S^{CT}_{0,K} = 0$ (that is, there are no tree-level counter-terms).
\end{enumerate}

\begin{lemma}
$S^{CT}_{I,K}$ is independent of $T$.  
\label{lemma counterterms independent T}
\end{lemma}
\begin{proof}
Observe that 
$$
P(\eps,T') - P(\eps,T)= P(T,T')   =  \int_{T}^{T'}  ( \GF \otimes 1) K_t  \d t
$$
is in $\E \otimes \E$ (that is, it has no singularities).    Let 
$$
A(\eps,\hbar) \in \Oo(\E, \mscr{A}_{\ge 0} \otimes \C[[\hbar]])
$$
be any (not necessarily local) functional.  Non-singularity of $P(T,T')$ implies that 
$$
\Gamma_{\le (I,K)} ( P(T,T'), A(\eps,\hbar))
$$
is non-singular (i.e., has a well-defined $\eps \to 0$ limit). 

We know by induction that $ \Gamma_{\le (I,K)}(P(\eps,T), S  - S^{CT}_{\le(I,K)}(\eps,T) ) $ is non-singular.  It follows that
$$
\Gamma_{\le(I,K)} \left(  P(T,T') ,  \Gamma_{\le (I,K)}\left( P(\eps,T), S - S^{CT}_{\le(I,K)}(\eps,T) \right) \right) 
$$
is non-singular.   But,
\begin{multline*}
\Gamma_{\le(I,K)} \left( P(T,T') ,  \Gamma_{\le(I,K)}\left( P(\eps,T), S - S^{CT}_{\le(I,K)}(\eps,T) \right) \right)   \\ 
\begin{split}
& =   \Gamma_{\le(I,K)}\left(P(\eps,T'), S- S^{CT}_{\le(I,K)} (\eps,T)\right) \\ 
& =  \Gamma_{\le(I,K)}\left(P(\eps,T'), S - S^{CT}_{<(I,K)} (\eps,T') \right) - S^{CT}_{(I,K)} (\eps,T)
 \end{split}
\end{multline*}
(where we are using the induction assumption that $S^{CT}_{<(I,K)} (\eps,T') = S^{CT}_{<(I,K)} (\eps,T)$).

This makes it clear that
\begin{align*}
S^{CT}_{(I,K)} (\eps,T) &= 
\text{Singular part of } \Gamma_{\le(I,K)}\left(P(\eps,T'), S - S^{CT}_{<(I,K)} (\eps,T') \right)\\
&= S^{CT}_{(I,K)} (\eps,T')
\end{align*}
as desired. 
\end{proof}

\begin{lemma}
$S^{CT}_{I,K}$ is local, so that 
$$S^{CT}_{I,K} \in \Oo_l ( \E, \mscr{A}_{< 0} ).$$  
\end{lemma}
\begin{proof}
This follows immediately from the fact that $S^{CT}_{I,K}$ is independent of $T$, and from the $T \to 0$ asymptotic expansion of the singular part of $\Gamma_{I,K}( P(\eps,T), S - S^{CT}_{<(I,K)} )$ proved in theorem A.
\end{proof}

\begin{lemma}
$S_{0,k}^{CT} = 0$.
\end{lemma}
\begin{proof}
All we need to show is that  $\Gamma_{0,k} (P(\eps,T), S ) $ is regular, for all $k$. This is an easy corollary of Lemma \ref{lemma local functional vector field}.  Namely, if we write $S = \sum \hbar^i S_{i,k}$ as usual, then Lemma \ref{lemma local functional vector field} implies that the functionals $S_{0,k}$ are of the form
$$
S_{0,k}(\alpha) = \ip { \Psi_k ( \alpha^{\otimes k-1}) ,\alpha } 
$$
for some polydifferential operator
$$
\Psi_k : \E^{\otimes k-1} \to \E
$$
and $\alpha \in \E$.  The tree-level terms $\Gamma_{0,k}(P(\eps,T),S)$ of $\Gamma(P(\eps,T),S)$ are obtained by composing these polydifferential operators $\Psi_k$ with each other and with the operator 
$$
P(\eps,T) : \E \to \E
$$
constructed by convolution with the kernel $P(\eps,T)$.  Note that 
$$
P(\eps,T)\star \alpha = \int_{\eps}^T  \GF e^{-t H} \alpha 
$$
This operator is a map $\E \to \E$,  even when $\eps = 0$. That is, $\int_0^T e^{-t H}$ takes smooth sections of the vector bundle $E$ on $M$ to smooth sections of $E$. This implies that all tree-level operators are non-singular, as desired.  
\end{proof}
This completes the proof. 
\end{proof}

\section{Renormalisation group flow and converse to theorem B}
Let $S \in \Oo_l(\E)[[\hbar]]$ be a local functional. The expression $\Gamma^R ( P(0,T) , S ) $ constructed using theorem B should be interpreted as the scale $T$ renormalised effective action.  
\begin{lemma}
The renormalisation group equation holds:
$$
\Gamma^R (   P(0,T') , S  )  = \Gamma ( P(T,T') , \Gamma^R ( P(0,T) , S )  ) .
$$ 
\end{lemma}
\begin{proof}
This follows from the fact that the counter-terms $S_{i,k}^{CT}$ are independent of $T$, and the identity
$$
\Gamma ( P(T,T') , \Gamma ( P(\eps,T) , S  - S^{CT} )  )  = \Gamma ( P(\eps,T'), S  -  S^{CT} ) .
$$
\end{proof} 

Thus we have seen that to any local functional $S \in \Oo(\E)[[\hbar]]$, we can associate a system of effective actions $\Gamma^R ( P (0,T), S)$ satisfying the renormalisation group equation.  In a moment we will see that a converse holds: all such systems of effective actions satisfying a certain locality condition arise in this way.

\begin{definition}
A \emph{system of effective actions} on the space of fields $\E$ is given by an effective action
$$
S^{eff} (T)\in \Oo(\E , \C[[\hbar]]  )
$$
for each $T \in \R_{> 0}$, which are all at least cubic modulo $\hbar$, and
such that
\begin{enumerate}
\item
The renormalisation group equation is satisfied,
$$
S^{eff}(T_2) =  \Gamma ( P(T_1,T_2) , S^{eff}(T_1)).
$$
\item
As $T \to 0$, $S^{eff}(T)$ must become local, in the following sense.  There must exist some $T$-dependent local functional 
$$\Phi \in \Oo_l( \E , \cinfty(0,\infty) \otimes \C[[\hbar]]  ) $$ 
where $T$ is the coordinate on $(0,\infty)$,
such that 
$$\lim_{T \to 0} \left( S^{eff}(T) - \Phi(T) \right)= 0.$$ 
(The $T \to 0$ limit of $S^{eff}(T)$ itself will generally not exist). 
\end{enumerate}
\end{definition}
The effective actions defined by $\Gamma^R ( P(0,T) , S)$ for $S \in \Oo_l(\E, \C[[\hbar]] )$ satisfy these axioms.

Fix any renormalisation scheme $\mscr{A}_{< 0}$. Then theorem B provides a map
\begin{align*}
S \in \Oo_l(\E) [[\hbar]] \text{ which is at least cubic modulo }\hbar& \to \text{ systems of effective actions } \\
S & \mapsto \{ \Gamma^R ( P(0,T), S ) \mid T \in \R_{> 0} \} 
\end{align*}
\begin{thmC}
\label{proposition rge local functional}
For any renormalisation scheme $\mscr{A}_{< 0}$, this map is a bijection.
\end{thmC}
This set of systems of effective actions is a canonical object associated to $(\E,Q, \GF)$, independent of the choice of renormalisation scheme.    Renormalisation and regularisation techniques other than those considered should lead to different ways of parametrising the same set of systems of effective actions.  For instance, if one could make sense of dimensional regularisation on general manifolds, one would hope to get simply a different parametrisation of this set.

From this point of view, the formalism of counter-terms is simply a convenient way to describe this set of systems of effective actions.   The counter-terms themselves, and the original action $S$, are not really meaningful in themselves.   

\begin{proof}[Proof of theorem C]
The proof is very simple.  Let $\{S^{eff}(T)\}$ be a system of effective actions.  We want to construct a corresponding local functional $S \in \Oo_l(\E) [[\hbar]]$ such that 
$$
S^{eff} (T) = \Gamma^R ( P(0,T) , S ) .
$$
Suppose, by induction, that we have constructed $S_{<(I,K)}$ such that 
$$
S_{(i,k)} ^{eff} (T) = \Gamma_{(i,k)} ^R ( P(0,T) , S_{<(I,K)} ) 
$$
for all $(i,k) < (I,K)$.    The initial case is when $i = 0$ and $k = 2$. Then,  $S^{eff}_{(0,2)} (T) = 0 = S_{(0,2)}$,  so the identity automatically holds. 

Assume $(I,K) \ge (0,3)$.  The renormalisation group equation says that 
\begin{align*}
S_{(I,K) }^{eff} (T) &= \Gamma_{(I,K)}(P(\eps, T) , S^{eff}_{\le (I,K)} (\eps )) \\
&=  S_{(I,K)}^{eff}(\eps) +   \Gamma_{(I,K)} ( P(\eps, T) , S_{<(I,K)}^{eff} (\eps ) ) \\
&=  S_{(I,K)}^{eff}(\eps) +   \Gamma_{(I,K)} ( P(\eps, T) , \Gamma_{<(I,K)}^R (P(0,\eps ), S_{<(I,K)} ) ) \\
&= S_{(I,K)}^{eff}(\eps) +   \Gamma_{(I,K)} ( P(\eps, T) , \Gamma_{\le (I,K)}^R (P(0,\eps ), S_{<(I,K)} ) )  -   \Gamma^R_{(I,K)} ( P(0, \eps) , S_{<(I,K)} ) ) \\
&=  S_{(I,K)}^{eff}(\eps) +   \Gamma^R_{(I,K)} ( P(0, T) , S_{<(I,K)} ) )  -  \Gamma^R_{(I,K)} ( P(0, \eps) , S_{<(I,K)} ) ).
\end{align*}
The left hand side  is independent of $\eps$, so the right hand side is also.  Thus, 
$$
 S_{(I,K)}^{eff}(\eps)  -  \Gamma^R_{(I,K)} ( P(0, \eps) , S_{<(I,K)} ) )
$$
is independent of $\eps$.  This allows us to define
$$
S_{(I,K)} =  S_{(I,K)}^{eff}(\eps)  -  \Gamma^R_{(I,K)} ( P(0, \eps) , S_{<(I,K)} ) ).
$$
With this definition, we automatically have
$$
S_{(I,K)}^{eff}(T) = \Gamma^R  ( P(0,T) , S_{\le (I,K)}   ) .
$$
It remains to show that $S_{(I,K)}$ is local.   This is an immediate consequence of the locality axiom for $S^{eff}(\eps)$.  There exists $\Phi(\eps)$, an $\eps$ dependent local functional, such that $\lim_{\eps \to 0} S^{eff}(\eps ) - \Phi(\eps ) = 0$.  Thus, 
\begin{align*}
S_{(I,K)} &=  S_{(I,K)}^{eff}(\eps)  -  \Gamma^R_{(I,K)} ( P(0, \eps) , S_{<(I,K)} ) ) \\
& =  S_{(I,K)}^{eff}(\eps) -  \Phi(\eps) + \Phi(\eps) -   \Gamma^R_{(I,K)} ( P(0, \eps) , S_{<(I,K)} ) )\\
&= \lim_{\eps \to 0}  \left( S_{(I,K)}^{eff}(\eps) -  \Phi(\eps) \right) + \lim_{\eps \to 0} \left( \Phi(\eps) -   \Gamma^R_{(I,K)} ( P(0, \eps) , S_{<(I,K)} ) ) \right) \\
&=  \lim_{\eps \to 0} \left( \Phi(\eps) -   \Gamma^R_{(I,K)} ( P(0, \eps) , S_{<(I,K)} ) ) \right) .
\end{align*}
This last quantity is local, because $\Phi(\eps)$ is local and $\Gamma^R_{(I,K)} ( P(0,\eps), S_{< (I,K)} )$ has a small $\eps$ asymptotic expansion in terms of local functionals. 

\end{proof}

As I explained in the introduction, this result elucidates how the renormalisation procedure depends on the choice of renormalisation scheme.    For the rest of the paper, we will fix one renormalisation scheme $\mscr{A}_{< 0}$.  When we talk about local functionals, we are in some sense really talking about systems of effective actions, but using our fixed renormalisation scheme to identify these two sets.  Thus, as long as all the statements we make are about the effective action $\Gamma^R(P(0,T), S)$, and not directly about the local functional $S$, everything we do is independent of the choice of renormalisation scheme. 

\section{Quantum master equation}

\subsection{Quantum master equation in finite dimensions}
Recall that to each element $\Phi \in \Sym^2 \E$ we have associated (in section \ref{section differential operators}) an order two differential operator $\partial_{\Phi}$ on $\Oo(\E)$.  Let $$\Delta_T =  -\partial_{K_T}.$$    When $T > 0$, this is a differential operator on $\Oo(\E)$. When $T = 0$, because $K_0$ is no longer an element of $\Sym^2 \E$ but of some distributional completion, the operator $\Delta_0$ is  ill-defined.   The operators $\Delta_T$ are odd, second-order differential operators on $\Oo(\E)$, which satisfy the equation
$$
\Delta_T^2 = 0.
$$ 
Thus, each $\Delta_T$ endows the algebra $\Oo(\E)$ with the structure of a Batalin-Vilkovisky algebra.  The ill-defined operator $\Delta_0$ is the ``physical'' Batalin-Vilkovisky operator.

Let us assume for a moment that $M$ is a point, so that $\E$ is simply a finite dimensional super vector space.  Then, the operator $\Delta_0$ is perfectly well defined.  In this situation, $\Delta_0$ is the Batalin-Vilkovisky odd Laplacian, as described in section \ref{section intro bv}.   The operator $K_0$ is simply the kernel representing the identity map on $\E$. In terms of a Darboux basis $x_i, \xi_i$ of $V$, where $x_i$ are even, $\xi_i$ are odd, and $\ip{x_i,\xi_i} = \delta_{i j}$, then
$$
K_0 =  -  \sum_i  \left(  x_i \otimes \xi_i + \xi_i \otimes x_i  \right)
$$
and
$$
\Delta_0 = - \partial_{K_0} = \sum \partial_{x_i} \partial_{\xi_i}.
$$

Still in the situation when $\dim M = 0$, the Poisson bracket on the algebra $\Oo(\E)$ can be written in terms of $\Delta_0$, using the equation
$$
 \{ f , g \} = \Delta_0(f g) - \Delta_0 (f) g - (-1)^{\abs f} f \Delta_0 (g) .
$$
When $\dim M > 0$, this Poisson bracket is ill-defined on $\Oo(\E)$. However, as we have seen, the Poisson bracket of an element of $\Oo_l(\E)$ with an element of $\Oo(\E)$ is well-defined. 

A functional $S \in \Oo(\E, \C[[\hbar]])$ satisfies the BV quantum master equation if 
$$
(Q + \hbar \Delta_0 ) \exp ( S / \hbar ) = 0.
$$
Again, this makes sense when $\dim M = 0$, but the left hand side of this equation is ill-defined  when $\dim M > 0$.  This equation can be re-expressed as 
$$
Q S + \tfrac{1}{2} \{ S, S \} + \hbar \Delta_0 S = 0.
$$

\subsection{The quantum master equation and the renormalisation group flow}
Let us return to the situation where $\dim M \ge 0$.  

Recall that 
$$
P(\eps,T)= \int_{\epsilon}^T L_t \d t = \int_{\eps}^T (\GF \otimes 1) K_t \d t.
$$
\begin{lemma}
$$
Q P _\epsilon^T = - K_T + K _\epsilon. 
$$
where $Q$ denotes the tensor product differential on $\E \otimes \E$.
\end{lemma}
\begin{proof}
Indeed, our sign conventions are such that $Q P(\epsilon,T)$ is the kernel for the operator 
\begin{align*}
\alpha &\mapsto Q \int_\epsilon^T \GF e^{-t H }  \alpha + \int_\epsilon^T \GF e^{-t H }  Q \alpha \\
&= \int_\epsilon^T H e^{-t H } \alpha \\
&= -e^{-T H }\alpha + e^{- \epsilon H } \alpha .
\end{align*}
Thus, the operator associated to $Q P(\epsilon,T)$ is the same as that associated to $-K_T + K_\epsilon$.
\end{proof} 

\begin{lemma}
Let $S \in \Oo(\E, \C[[\hbar]])$.  Then $S$ satisfies the $\Delta_{T_1}$ quantum master equation
$$
(Q + \hbar  \Delta_{T_1} ) \exp ( S / \hbar ) = 0
$$
if and only if $\Gamma( P(T_1,T_2) , S)$ satisfies the $\Delta_{T_2}$ quantum master equation
$$
(Q + \hbar \Delta_{T_2} )  \exp (\Gamma( P(T_1,T_2) , S) / \hbar  )  = 0. 
$$
\label{lemma propagator homotopy}
\end{lemma}
\begin{proof}
Indeed,
$$
[ \partial_{P(T_1,T_2)}  , Q ] = \partial_{Q P(T_1,T_2) } = \partial_{K_{T_1} } - \partial_{K_{T_2} } = \Delta_{T_2} - \Delta_{T_1} . 
$$
This implies that 
$$
[ Q , \exp ( \partial_{P(T_1,T_2)}  )  ] = \exp (\partial_{P(T_1,T_2)}  ) ( \Delta_{T_1} - \Delta_{T_2} ) . 
$$
Thus,
\begin{align*}
(Q + \hbar \Delta_{T_2 } )  \exp (\Gamma( P(T_1,T_2) , S) / \hbar  )  &=  (Q + \hbar \Delta_{T_2 } ) \exp ( \hbar \partial_{P(T_1,T_2)} ) \exp ( S / \hbar  )\\ 
&= \exp ( \hbar \partial_{P(T_1,T_2)}  ) ( Q + \hbar \Delta_{T_1 }  )  \exp ( S / \hbar  ) .
\end{align*}
The converse follows because the operator $\exp ( \hbar \partial_{P(T_1,T_2)}  )$ is invertible.

\end{proof} 

Suppose we have an effective action $S^{eff}(T_1) \in \Oo(\E, \C[[\hbar]])$ at scale $T_1$.   Then $\Gamma( P(T_1,T_2), S^{eff}(T_1))$  is the effective action at scale $T_2$.      What we see is that if the effective action $S^{eff}(T_1)$ satisfies the scale $T_1$ quantum master equation, then the effective action $S^{eff}(T_2)$ at scale $T_2$ satisfies the scale $T_2$ quantum master equation, and conversely.

\subsection{Renormalised quantum master equation}

We know, using the renormalisation techniques we have discussed so far, how to make sense of the expression $\Gamma(P(0,T) , S ) $ if $S \in \Oo_l(\E, \C[[\hbar]])$ is a local functional which is at least cubic modulo $\hbar$.  The  renormalised scale $T$ effective action is
$$\Gamma^R ( P(0,T) , S ) \defeq \lim_{\eps \to 0} \Gamma(P(\eps,T), S^R )$$ 
where $S^R = S - S^{CT}$.  

We would like to define a renormalised quantum master equation, which is a replacement for the ill-defined equation
$$
(Q + \hbar \Delta_0 ) e^{S / \hbar }= 0.
$$
When $\dim M = 0$, so $\E$ is finite dimensional, this equation is well-defined and is equivalent to the equation
$$
(Q + \hbar \Delta_T )\exp \left( \hbar^{-1} \Gamma ( P(0,T) , S ) \right)= 0.
$$
In the infinite dimensional situation, we will take this second equation as a definition, where we use the renormalised version $\Gamma^R ( P(0,T), S )$ of $\Gamma(P(0,T), S)$.  Thus, 
\begin{definition}
$S \in \Oo_l(\E, \C[[\hbar]])$ satisfies the renormalised quantum master equation if
$$
(Q + \hbar \Delta_T ) \exp \left( \hbar^{-1}  \Gamma^R (P(0,T), S ) \right) = 0 .
$$
\end{definition}
Theorems B and C show that there is a bijection between the set of local functionals $S \in \Oo_l(\E, \C[[\hbar]])$ and the set of systems of effective actions, sending $S$ to $\{\Gamma^R ( P(0,T) , S) \mid T \in \R_{> 0} \}$ .  The renormalised QME is thus just saying that the scale $T$ effective action satisfies the scale $T$ QME.

It is immediate from Lemma \ref{lemma propagator homotopy} that this condition is independent of $T$.    That is, if $S$ satisfies the renormalised QME for one value of $T$, then it does so for all other values of $T$. 

It seems to me that a successful quantisation of the quantum field theory with classical action $S_0$ consists of replacing $S_0$ by a solution $S = S_0 + \sum_{i > 0} \hbar^i S_i$ of the renormalised quantum master equation. Equivalently, a quantisation of a classical action $S_0$ is given by replacing  $S_0$ by a system of effective actions $\{S^{eff}(T) \mid T \in \R_{> 0} \}$, which satisfy the quantum master equation, and which, modulo $\hbar$, converge to $S_0$ as $T \to 0$. 

\begin{lemma}
Suppose that $S \in \Oo_l(\E, \C[[\hbar]])$ satisfies the renormalised quantum master equation.  Then 
$$
\Gamma^R ( P(0,\infty), S ) \mid_{\Ker H} 
$$
satisfies the quantum master equation on the finite dimensional vector space 
$$
\Ker H = H^\ast(\E, Q) . 
$$
\end{lemma}
\begin{proof}
This is immediate from the fact that the heat kernel $K_\infty$ is the element of $\Ker H \otimes \Ker H \subset \E \otimes \E$ which (as a kernel on $\Ker H$) represents the identity map $\Ker H \to \Ker H$. 
\end{proof} 
Recall that
$\Gamma^R ( P(0,\infty), S )(a)$ is a renormalised version of the functional integral
$$
\hbar \log \int_{x \in (\Im \GF ) _\R }  \exp\left( \tfrac{1}{2\hbar}  \ip{x,Q x}  + \tfrac{1}{\hbar} S ( x + a ) \right) .
$$
If we take $a$ to be in $\Ker H = H^\ast(\E,Q)$, then this functional integral is the Wilsonian effective action in the Batalin-Vilkovisky formalism, obtained by integrating over a Lagrangian $\Im \GF$ in the space $\Im \GF \oplus \Im Q$ of positive eigenvalues of the Hamiltonian $H$.

\section{Homotopies of solutions of the quantum master equation}

\label{section simplicial sets}
In the usual finite-dimensional Batalin-Vilkovisky formalism \cite{Sch93, AleKonSch97}, one reason the quantum master equation is so important is that it implies that the functional integral is independent of the choice of gauge fixing condition (usually, a Lagrangian or isotropic subspace in the space of fields).    

We would like to prove a version of this in the infinite dimensional situation, for the renormalised quantum master equation.  However, things are more delicate in this situation; the renormalised QME itself depends on the choice of a gauge fixing condition.  What we will show is that if $\GF(t)$ is a one-parameter family of gauge fixing conditions, then the set of homotopy classes of solutions to the renormalised QME using $\GF(0)$ is isomorphic to the corresponding set using $\GF(1)$.  This result is a corollary of a result about certain simplicial sets of gauge fixing conditions and of solutions to the renormalised QME.

As I explained earlier, we are fixing a renormalisation scheme throughout  the paper, which allows us to identify the set of local functionals $S \in \Oo_l(\E) [[\hbar]]$ with the set of systems of effective actions.  Everything we say should be in terms of the effective action $\Gamma^R(P(0,T), S)$ and not directly in terms of the local functional $S$.  As long as we stick to this, everything we do will be independent of the choice of renormalisation scheme. 
\subsection{Families and simplicial sets of gauge fixing conditions}

Let us denote the $n$-simplex by $\Delta^n$.

\begin{definition}
 A family of gauge fixing conditions parametrised by $\Delta^n$ is a first order differential  operator 
$$
\GF : \E \otimes \cinfty(\Delta^n) \to \E \otimes \cinfty (\Delta^n)
$$
which is $\cinfty(\Delta^n)$ linear, and which satisfies a family version of the axioms for a gauge fixing condition.  Namely, 
\begin{enumerate}
\item
the operator $H_0 = [Q \otimes 1, \GF]$ must be a smooth family of generalised Laplacians in the sense of \cite{BerGetVer92}, section 9.5.
\item
there is a direct sum decomposition (of $\cinfty(\Delta^n)$ modules)
$$
\E \otimes \cinfty(\Delta^n) = \Im (Q\otimes 1) \oplus \Im \GF \oplus \Ker H_0 .
$$
\end{enumerate}
\end{definition}
We will normally just write $Q$ for the operator $Q \otimes 1$  on $\E \otimes \Omega^\ast(\Delta^n)$, and similarly we will write $\d_{DR}$ for $ 1 \otimes \d_{DR}$.  The operator $\GF$ extends in a unique $\Omega^\ast(\Delta^n)$-linear way to $\E \otimes \Omega^\ast(\Delta^n)$.  Let us use the notation
$$
H = [ Q + \d_{DR} , \GF ]  : \E \otimes \Omega^\ast (\Delta^n) \to  \E \otimes \Omega^\ast (\Delta^n) . 
$$
It is easy to see that $H$ is $\Omega^\ast(\Delta^n)$-linear.  If $H_1 = [\d_{DR}, \GF]$, then $H = H_0 + H_1$ and $H_1$ is an order $1$ differential operator.  Thus, the symbol of $H$ is the same as that of $H_0$, and, by assumption, the symbol of $H_0$ is given by a smooth family of metrics on $M$, times the identity on the vector bundle $E$ on $M$.    Thus, we can think of $H$ as a smooth family of generalised Laplacians on $\E$ parametrised by the supermanifold $\Delta^n \times \R^{0,n}$.  The results of \cite{BerGetVer92}  imply that there exists a unique heat kernel for $H$, 
$$
K_t \in \E \otimes \E \otimes \Omega^\ast(\Delta^n) .
$$

Let us denote by $\mbf{GF} [n]$ the set of families of gauge fixing conditions parametrised by $\Delta^n$.  If $\phi : \Delta^n \to \Delta^m$ is a smooth map (for instance, the inclusion of a face) then there is an induced map $\phi^\ast : \mbf{GF}[m]   \to \mbf{GF}[n]$.   In this way, the sets $\mbf{GF}[n]$ form a simplicial set, which we denote $\mbf{GF}$.

\subsection{Simplicial sets of solutions to the renormalised QME}

Let $$S \in \Oo_l(\E, \Omega^\ast(\Delta^n) [[\hbar]] )$$ be a family of functionals on $\E$, which, as always, is at least cubic modulo $\hbar$.   As before, if $P \in \E \otimes \E \otimes \Omega^\ast (\Delta^n)$, one can define
$$
\Gamma ( P , S ) = \hbar \log \left( \exp ( \hbar \partial_P ) \exp ( S  / \hbar ) \right)  \in \Oo ( \E , \Omega^\ast(\Delta^n) [[\hbar]]  ) .
$$

Let us suppose, as above, that we have a smooth family of gauge fixing conditions $\GF$ parametrised by $\Delta^n$. Let $H$ and $K_t$ be as above.  This allows us to construct the propagator we use, 
$$
P(\eps,T)= \int_{\epsilon}^T  (\GF \otimes 1) K_t \d t .
$$
Thus, one can define the family version of the renormalisation group flow, which sends
$$
f \in \Oo(\E, \Omega^\ast(\Delta^n) [[\hbar]] ) \to \Gamma(P(\eps,T) ,f ) \in \Oo(\E, \Omega^\ast(\Delta^n) [[\hbar]] ).
$$
\begin{definition}
A family of systems of effective actions on $\E$ parametrised by the supermanifold $\Delta^n \times \R^{0,n}$ is given by an effective action
$$
S^{eff}(T) \in \Oo(\E, \Omega^\ast(\Delta^n) [[\hbar]] )
$$
for all $T \in \R_{> 0}$, 
such that
\begin{enumerate}
\item
The renormalisation group equation is satisfied,
$$
S^{eff}(T_2) =  \Gamma ( P(T_1,T_2) , S^{eff}(T_1)).
$$
\item
As $T \to 0$, $S^{eff}(T)$ must become local, in the following sense.  There must exist some $T$-dependent local functional 
$$\Phi \in \Oo_l( \E , \Omega^\ast(\Delta^n) \otimes \cinfty(0,\infty) \otimes \C[[\hbar]]  ) $$ 
where $T$ is the coordinate on $(0,\infty)$,
such that 
$$\lim_{T \to 0} \left( S^{eff}(T) - \Phi(T) \right) = 0.$$ 
(The $T \to 0$ limit of $S^{eff}(T)$ itself will generally not exist). 
\end{enumerate}
\end{definition}

Theorems A, B and C apply in this situation. Indeed, theorem A is stated and proved in this generality in the appendix, and the proofs of the family version of theorems  B and C is identical to the proof we gave earlier.    Thus, we can define 
$$
\Gamma^R(P(0,T), S ) = \lim_{\eps \to 0} \Gamma(P(\eps,T), S - S^{CT} ) \in \Oo ( \E , \Omega^\ast(\Delta^n) \otimes \C [[\hbar]]  ) . 
$$
This map gives us a  bijection between the set of local functionals $S \in \Oo_l ( \E, \Omega^\ast(\Delta^n) \otimes \C[[\hbar]] )$, and the set of families of systems of effective actions parametrised by $\Delta^n \times \R^{0,n}$. This bijection depends on the choice of renormalisation scheme, which we fix throughout the paper.

We say that $S$ satisfies the renormalised quantum master equation if
$$
( \d_{DR} + Q + \hbar \Delta_T) \exp ( \hbar^{-1} \Gamma^R (P(0,T), S )  ) = 0 . 
$$
This condition, as before, is independent of $T$.   
\begin{definition}
Define a simplicial set $\mbf{BV}_\E$ by saying that $\mbf{BV}_\E[n]$ is the set of triples 
\begin{enumerate}
\item
$\GF \in \mbf{GF}[n]$
\item $S \in \Oo_l(\E, \Omega^\ast(\Delta^n) [[\hbar]] )$ which is at least cubic modulo $\hbar$, and which satisfies the renormalised quantum master equation corresponding to $\GF$.
\end{enumerate}

Similarly, define a simplicial set $\mbf{BV}_{H^\ast(\E,Q) }$, by saying $\mbf{BV}_{H^\ast(\E,Q) }[n]$  is the set of functions $S \in \Oo (H^\ast(\E, Q)  ) \otimes \Omega^\ast(\Delta^n) [[ \hbar ]]$ which satisfy the quantum master equation $$ ( \d_{DR} + \hbar \Delta_{H^\ast(\E, Q) }  ) e^{S / \hbar } = 0.$$  
\end{definition}

An equivalent (and more natural) definition of the simplicial set $\mbf{BV}_\E$ would be to say that $\mbf{BV}_\E[n]$ is the set of families of systems of effective actions, parametrised by $\Delta^n \times \R^{0,n}$, such that the scale $T$ effective action $S^{eff}(T)$ satisfies the scale $T$ QME, 
$$
(Q + \d_{DR} + \hbar \Delta_T ) e^{ S^{eff}(T) / \hbar } = 0 . 
$$
Of course, the choice of renormalisation scheme (which we are fixing throughout the paper) leads to an equivalence between the two definitions.

\begin{lemma}
There is a map
$$
\mbf{BV}_\E[n] \to \mbf{BV}_{H^\ast(\E,Q)} [n]
$$
of simplicial sets, which sends
$$
S \mapsto \Gamma^R ( P(0,\infty) , S )\mid_{ H^\ast(\E, Q)  } .
$$
This extends to a map of simplicial sets
$$
\mbf{BV}_\E[n] \to \mbf{BV}_{H^\ast(\E,Q)} [n]
$$
in the natural way.
\end{lemma}

\subsection{Fibration property}
Let $\GF$ be a gauge fixing condition on $\E$.  Let $\mbf{BV}_{\E, \GF}$ denote the simplicial set of solutions to the renormalised quantum master equation with fixed gauge fixing condition $\GF$.  In other words, $\mbf{BV}_{\E, \GF}$ is the fibre of $\mbf{BV}$ over the point $\GF \in  \mbf{GF} [0]$, under the natural map $\mbf{BV}_{\E} \to \mbf{GF} $.

 The set $\pi_0 \left( \mbf{BV}_{\E, \GF} \right)$ is the set of homotopy classes of solutions to the renormalised quantum master equation.  The fact that the map $\mbf{BV}_{\E, \GF} \to \mbf{BV}_{H^\ast(\E,Q)}$ is a map of simplicial sets tells us that there is an induced map
\begin{equation}
\pi_0 \left( \mbf{BV}_{\E, \GF} \right)  \to  \pi_ 0 \left( \mbf{BV}_{H^\ast(\E,Q)}  \right). \tag{$\dagger$}  \label{map to qme cohomology} 
\end{equation}
We would like to show that the set $\pi_0 ( \mbf{BV}_{\E, \GF}) $, and the map \eqref{map to qme cohomology}, are in some sense independent of $\GF$..  This will follow from an abstract result.

\begin{theorem}
The map
$$
\mbf{BV}_\E \to \mbf{GF} 
$$
is a fibration of simplicial sets.  
\label{theorem fibration}
\end{theorem}
We will prove this later.

\begin{corollary}
Let $\GF(t)$ is a  smooth families of gauge fixing conditions  parametrised by $[0,1]$.  Then the sets $\pi_0 ( \mbf{BV}_{\E, \GF(0)} ) $, $\pi_0 ( \mbf{BV}_{\E, \GF(1) } ) $ are canonically isomorphic. Further, the diagram
$$
\xymatrix{ \pi_0 \left( \mbf{BV}_{\E, \GF(0)} \right)  \ar[r] \ar[d] &  \pi_0 \left( \mbf{BV}_{\E, \GF(1)}  \right) \ar[dl] \\
\pi_ 0 \left(\mbf{BV}_{H^\ast(\E, Q )    } \right) } 
$$
commutes.  
\end{corollary}
\begin{proof}
This follows immediately from the path and homotopy lifting properties for fibrations of simplicial sets.  
\end{proof}

\section{A local obstruction to solving the renormalised quantum master equation} 
We will prove Theorem \ref{theorem fibration} using obstruction theory.  In this section we will construct  the required obstructions to solving the renormalised QME. 

Let us fix a smooth family $\GF$ of gauge fixing conditions  parametrised by $\Delta^n$. Let $S \in \Oo_l(\E, \Omega^\ast(\Delta^n) \otimes \C[[\hbar]])$ be, as always, at least cubic modulo $\hbar$.   One way to write the renormalised quantum master equation for $S$ is the equation
$$
\lim_{\eps \to 0}    (Q +  \d_{DR} + \Delta_T ) \exp ( \hbar  \partial_{P(\eps,T)} ) \exp ( S^R / \hbar )   = 0. 
$$

One can re-express this using the  following identity: 
\begin{multline}
\Gamma \left( P(\eps,T)+ \delta \Delta_\eps ,  S^R + \delta Q S^R + \delta \d_{DR} S^R \right)  = \\ \hbar \log \left( (1 + \delta Q + \delta \d_{DR} + \delta  \Delta_T ) \exp ( \hbar  \partial_{P(\eps,T)} ) \exp ( S^R / \hbar )  \right)  \label{equation Q inside gamma}
\end{multline}
where $\delta$ is an odd parameter.  This identity is a corollary of the fact that 
$$
[\partial_{P(\eps,T)} , Q + \d_{DR}] = \Delta_T -\Delta_\eps. 
$$
Thus, we  have an equivalent expression of the renormalised quantum master equation; the renormalised QME holds if and only if 
\begin{equation*}
\lim_{\eps \to 0} \frac{\d}{\d \delta} \left<   \Gamma \left( P(\eps,T)+ \delta \Delta_\eps,  S^R + \delta Q S^R + \delta \d_{DR} S^R \right)  \right> = 0 .
 \tag{$\dagger$} \label{rqme} 
\end{equation*}

Recall that our choice  of renormalisation scheme gives us a subspace
$$
\mscr{A}_{< 0} \subset \mscr{A}
$$
which is complementary to the subspace $\mscr{A}_{\ge 0}$ .
Let
$$
\mscr{A}_{\le 0} = \mscr{A}_{< 0} \oplus \C \subset \mscr{A}.
$$
If we denote  by $\mscr{A}_{> 0}$ the space of functions whose $\eps \to 0$ limit exists and is equal to zero, we have a direct sum decomposition
$$
\mscr{A} = \mscr{A}_{\le 0} \oplus \mscr{A}_{> 0} . 
$$
\begin{proposition}
\label{proposition obstruction}
\begin{enumerate}
\item
Let $S \in \Oo_l(\E, \Omega^\ast(\Delta^n) \otimes \C[[\hbar]] )$ be an even element.  Then
there exists a unique odd element 
$$O(S) \in \Oo_l(\E,  \mscr{A}_{\le 0}\otimes \Omega^\ast(\Delta^n) \otimes \C[[\hbar]] )$$ such that 
$$\lim_{\eps \to 0} \frac{\d}{\d \delta} \left<   \Gamma \left( P(\eps,T)+ \delta \Delta_\eps,  S^R + \delta Q S^R + \delta \d_{DR} S^R + \delta O(S)   \right)  \right>    $$ 
exists, and has value $0$. 
\item
$O_{(I,K)}(S)$ depends only on $S_{\le (I,K)}$, and 
$$
O_{(I,K)}(S_{<(I,K)} + \hbar^I S_{(I,K)} ) =   \hbar^I Q S_{(I,K)} + \hbar^I \d_{DR} S_{(I,K)} + O_{(I,K)}(S_{<(I,K)} ). 
$$
\item
Further,  if $k > 0$, the germ of $O_{i,k}(S)$ at $x \in M$ depends only on the germ of the data $(\E, \ip{\ , \ }, Q , \GF , \{ S_{i,k} \mid k > 0 \} )$ near $x\times \Delta^n$.  
\end{enumerate}

\end{proposition}
Thus, $S$ satisfies the renormalised quantum master equation if and only if $O(S) = 0$.
\begin{proof}
The proof is very similar to that of theorem B.  

Suppose we've constructed $O_{(i,k)}(S) \in \Oo_l(\E, \Omega^\ast(\Delta^n) \otimes \mscr{A}_{\le 0})$ for all $(i,k) < (I,K)$, which are independent of $T$, and which are such that 
$$\frac{\d}{\d \delta} \left<   \Gamma_{< (I,K)} \left( P(\eps,T)+ \delta \Delta_\eps,  S^R + \delta Q S^R + \delta \d_{DR} S^R + \delta O_{<(I,K)}(S)   \right)  \right>  $$
is regular at $\eps = 0$ with value $0$.  

Define
\begin{multline*}
O_{(I,K)}(S,T) = -\text{projection onto $\mscr{A}_{\le 0}$ of } \\   \frac{\d}{\d \delta} \left<   \Gamma_{(I,K)} \left( P(\eps,T)+ \delta \Delta_\eps,  S^R + \delta Q S^R + \delta \d_{DR} S^R + \delta O_{<(I,K)}(S)   \right)  \right>.  
\end{multline*}
The second clause of the proposition is immediate from this definition.  What  remains to be shown is that  $O_{(I,K)}(S,T)$ is independent of $T$, and so (by the same argument as in the proof of theorem B) local.

It follows from the definition of $O_{(I,K)}(S,T)$ that
$$
\lim_{\eps \to 0}\Gamma_{\le (I,K)} \left( P(\eps,T)+ \delta \Delta_\eps,   S^R 
+ \delta Q S^R + \delta \d_{DR} S^R  +
 \delta O_{< (I,K)}(S) + \hbar^I \delta O_{(I,K)}(S,T)   \right)  
$$
exists, and this limit is independent of $\delta$. Thus, the same holds for 
\begin{multline*}
\Gamma_{\le(I,K)} \left( P(T,T')  ,  \Gamma_{\le(I,K)}\left( P(\eps,T)+ \delta \Delta_\eps, S^R + \delta Q S^R  \delta \d_{DR} S^R 
\right. \right.
\\ 
 \left. \left.
+ \delta O_{< (I,K)}(S) +  \hbar^I \delta O_{(I,K)}(S,T) \right) \right) . 
\end{multline*}

Observe that
\begin{multline*}
\Gamma_{\le(I,K)} \left( P(T,T')  ,  \Gamma_{\le(I,K)}\left( P(\eps,T)+ \delta \Delta_\eps, S^R + \delta Q S^R + \delta \d_{DR} S^R  \right. \right. \\
\shoveright {\left. \left.
 + \delta O_{< (I,K)}(S) +  \hbar^I  \delta O_{(I,K)}(S,T) \right) \right)   }
\\ =   \Gamma_{\le(I,K)}\left(P(\eps,T') + \delta \Delta_\eps, S^R + \delta Q S^R 
 + \delta \d_{DR} S^R 
 \right. \\
\shoveright{ \left.
 + \delta O_{< (I,K)}(S)) +  \hbar^I \delta O_{(I,K)}(S,T) \right) }
\\ =  \Gamma_{\le(I,K)}\left(P(\eps,T') + \delta \Delta_\eps, S^R + \delta Q S^R + \delta \d_{DR} S^R  
 \right. \\
\left.
+ \delta O_{< (I,K)} (S)  \right) + \hbar^I \delta O_{ (I,K)}(S,T).
\end{multline*}
As we have seen,  the first line in this expression has an $\eps \to 0$ limit which is independent of $\delta$.  The same is therefore true for the last line; it follows that $O_{(I,K)}(S,T)  = O_{(I,K)} (S,T')$, and so $O_{(I,K)}(S) = O_{(I,K)}(S,T)$ is local as desired.  

\end{proof} 

We would like to use an inductive, term-by-term method to construct solutions to the renormalised QME.    The following lemma allows us to do so.
\begin{lemma}
Let $S \in \Oo_l ( \E, \Omega^\ast(\Delta^n) \otimes \C[[\hbar]] )$ be a local action functional such that $O_{< (I,K)} (S) = 0$.  (In other words, $S$ satisfies the QME up to order $(I,K)$).  Then:
\begin{enumerate}
\item
$$\lim_{\eps \to 0} O_{(I,K)}(S)$$ exists. This implies that $O_{(I,K)}(S)$ is independent of $\eps$, and so is an element
$$
O_{(I,K)}(S) \in  \Oo_l(\E, \Omega^\ast(\Delta^n)  ).
$$
\item
$$(Q + \d_{DR}) O_{(I,K)}(S)  = 0. $$ 
\end{enumerate}
\label{closed local obstruction}
\end{lemma} 
\begin{remark}
This  lemma implies that if $S$ solves the $<(I,K)$ renormalised QME, that is, $O_{<(I,K)}(S)= 0$, then lifting $S$ to a solution to the $\le(I,K)$ renormalised QME amounts to finding some $S'\in \Oo_l(\E, \Omega^\ast(\Delta^n))$, homogeneous of degree $K$ as a function on $\E$, with  $(Q+ \d_{DR}) S' = -  O_{(I,K)}(S)$.  The desired solution of the renormalised QME will be $S + \hbar^I S'$. 
\end{remark}

\begin{proof}
By definition,
\begin{multline*}
O_{(I,K)}(S) =   -\text{projection onto $\mscr{A}_{\le 0}$ of }   \\ \frac{\d}{\d \delta} \left<   \Gamma_{(I,K)} \left( P(\eps,T)+ \Delta_\eps,  S^R + \delta Q S^R + \delta \d_{DR} S^R + \delta O_{<(I,K)}(S)   \right)  \right>.  
\end{multline*}
Equation \eqref{equation Q inside gamma} says that
\begin{multline*}
  \Gamma \left( P(\eps,T)+ \Delta_\eps,  S^R   + \delta Q  S^R  + \delta \d_{DR} S^R  \right)  \\ 
  = \hbar \log \left(  (1 + \delta Q + \delta \d_{DR} + \delta  \Delta_T ) \exp ( \hbar  \partial_{P(\eps,T)} ) \exp ( S^R / \hbar ) \right) .
\end{multline*}
Since $O_{<(I,K)}(S) = 0$, it follows that 
\begin{multline*}
O_{(I,K)}(S) =   -\text{projection onto $\mscr{A}_{\le 0}$ of }  \\ \hbar \log \left(  (1 + \delta Q + \delta \d_{DR} + \delta  \Delta_T ) \exp ( \hbar  \partial_{P(\eps,T)} ) \exp ( S^R / \hbar  ) \right) .
\end{multline*}
It follows that the $\eps \to 0$ limit of $O_{(I,K)}(S)$ exists.  This implies that $O_{(I,K)}(S)$ is constant as a function of $\eps$. 

The equation $O_{<(I,K)}(S) = 0$ implies that 
\begin{multline*}
\hbar^I O_{I,K} (S) = - \text{leading term of }  \\ \lim_{\eps \to 0} \frac{ \d}{\d \delta}  \left< \hbar \log \left(  (1 + \delta Q + \delta \d_{DR}  +\delta T \Delta_T ) \exp ( \hbar  \partial_{P(\eps,T)} ) \exp ( S^R / \hbar ) \right) \right>
\end{multline*}
where ``leading term'' means the first non-zero term in the expansion labelled by $(i,k)$, where the $(i,k)$ term, as always, is $\hbar^i$ times a homogeneous function of degree $k$ on $\E$.  We use the  lexicographical ordering on the labels $(i,k)$.

It follows that
\begin{multline*}
\hbar^I O_{I,K}(S) = - \text{leading term of }  \\ ( \delta Q + \delta \d_{DR} +  T \Delta_T )  \lim_{\eps \to 0}  \left< \exp ( \hbar  \partial_{P(\eps,T)} ) \exp ( S^R / \hbar  ) \right>
\end{multline*}
so that
\begin{align*}
\hbar^I  (Q +  \d_{DR})  O_{I,K}(S) =& -\text{leading term of } \\  & (\delta Q + \delta \d_{DR} + T \Delta_T )^2 \lim_{\eps \to 0}  \left< \exp ( \hbar  \partial_{P(\eps,T)} ) \exp ( S^R / \hbar  ) \right> 
\\ 
  =& 0  
\end{align*}
as desired.

\end{proof}

\section{Proof of Theorem \ref{theorem fibration}}
Now we can prove Theorem \ref{theorem fibration}, that the map $\mbf{BV}_\E \to \mbf{GF}$ is a fibration of simplicial sets.
\subsection{Preliminary definitions}
Sending $\Delta^n$ to  $\Omega^\ast(\Delta^n)$ defines a simplicial algebra, which we denote $\Omega^\ast_\Delta$.  If $X$ is a simplicial set, let $\Omega^\ast(X)$ denote the space of maps $X \to \Omega^\ast_\Delta$ of simplicial sets.  Concretely, an element of $\Omega^\ast(X)$ is something which assigns to an $n$-simplex of $X$ an element of $\Omega^\ast(\Delta^n)$, in a way compatible with the natural maps between simplices.

The spaces 
$$
\Oo_l(\E, \Omega^\ast(\Delta^n) )
$$
assemble into a simplicial object in the category of chain complexes, which we denote $\Oo_l(\E, \Omega^\ast_{\Delta}) $.  For any simplicial set $X$, let $\Oo_l(\E, \Omega^\ast(X) )$ denote the space of maps $X \to \Oo_l(\E, \Omega^\ast_{\Delta}) $.    More concretely, we can view an element of $\Oo_l(\E, \Omega^\ast(X) )$ as being an element of $\Oo_l(\E, \Omega^\ast(\Delta^n) )$ for each $n$-simplex of $X$, compatible with the natural maps between simplices.

Earlier we defined a map 
\begin{align*}
\Oo_l(\E, \Omega^\ast(\Delta^n) \otimes \C [[\hbar]] )  &\to \Oo_l(\E, \mscr{A}_{\le 0} \otimes \Omega^\ast(\Delta^n)\otimes \C [[\hbar]] ) \\
S &\mapsto O(S) .
\end{align*}
This gives a map of simplicial sets
$$
\Oo_l(\E, \Omega^\ast_{\Delta} \otimes \C [[\hbar]])  \to \Oo_l(\E, \mscr{A}_{\le 0}\otimes \Omega^\ast_{\Delta}\otimes \C [[\hbar]] ).
$$
For each $X$, we get an obstruction  map
$$
O : \Oo_l(\E, \Omega^\ast(X)\otimes \C [[\hbar]] )  \to \Oo_l(\E, \mscr{A}_{\le 0} \otimes  \Omega^\ast(X)\otimes \C [[\hbar]]).
$$
$S \in \Oo_l(\E, \Omega^\ast(X) \otimes \C[[\hbar]] )$ satisfies the renormalised QME if and only if $O(S) = 0$.

$\mbf{BV}_{\E}$ is the simplicial set whose $n$-simplices are elements $S \in \Oo_l(\E, \Omega^\ast(\Delta^n)  ) $ such that $O(S) = 0$.  

\begin{definition}
Let  $\mbf{BV}_{\E, \le(I,K)}$ be the simplicial set whose $n$-simplices are given by pairs
$$
(\GF ,S ) 
$$ 
where $\GF$ are smooth families of gauge fixing conditions, and 
$$S \in \Oo_l(\E, \Omega^\ast(\Delta^n)  )$$ 
is such that $S_{(i,k)} = 0$ if $(i,k) > (I,K)$ and $O_{\le (I,K)} (S) = 0$.  

Similarly, define $\mbf{BV}_{\E, <(I,K)}$ to be the simplicial set whose $n$-simplices are pairs $(\GF,  S)$ where $\GF$ are as before, and 
$S \in \Oo_l(\E, \Omega^\ast(\Delta^n)  ) $ is such that $S_{(i,k)} = 0$ if $(i,k) \ge (I,K)$ and $O_{< (I,K)} (S) = 0$.  
\end{definition}  
Suppose $X$ is equipped with a map $X \to \mbf{GF}$.  Then a lift of this to a map  $X \to \mbf{BV}_{\E, \le(I,K)}$ is given by $S \in \Oo_l(\E, \Omega^\ast(X) )$ such that $S_{>(I,K)} = 0$ and $O_{\le(I,K)} (S) = 0$.   A similar description holds for maps to $\mbf{BV}_{\E, <(I,K)}$.  

The following lemma is an immediate corollary of Lemma \ref{closed local obstruction}. 
\begin{lemma}
Let $f : X \to \mbf{BV}_{\E, <(I,K)}$ be a map.  Then the obstruction $O_{(I,K)}(X)$ to lifting $f$ to a map $X \to \mbf{BV}_{\E, \le (I,K) }$ is independent of the parameter $\eps$, and so is an element of  $\Oo_l(\E, \Omega^\ast(X) )$. This is closed, $(Q + \d_{DR}) O_{(I,K)}(X) = 0$.  

A lift of $f$ to a map $ X \to \mbf{BV}_{\E, \le (I,K)}$ is given by an element $S' \in \Oo_l(\E, \Omega^\ast(X) )$, homogeneous of degree $K$ as a function on $\E$, with 
$$
(Q + \d_{DR}) S' = - O_{(I,K)}(X) . 
$$
\end{lemma}

\subsection{Completion of proof} 
Note that $\mbf{BV}_{\E, \le (0,0) }$ is $\mbf{GF}$.  Thus, to prove Theorem \ref{theorem fibration},  it suffices to show that
\begin{proposition}
The map
$$\mbf{BV}_{\E, \le (I,K) } \to \mbf{BV}_{\E, < (I,K) }$$
is a fibration.  
\label{proposition fibration}
\end{proposition}
\begin{proof}

Let $\Lambda^n \subset \Delta^n$ denote an $n$-horn, obtained by removing any face from $\partial \Delta^n$.  Suppose we have a map $\Lambda^n \to \mbf{BV}_{\E, \le (I,K) }$, and an extension of this to a map $\Delta^n \to \mbf{BV}_{\E, <(I,K) }$.   We need to show that there exists a map $\Delta^n \to \mbf{BV}_{\E, \le (I,K) }$ such that the diagram
$$
\xymatrix{
\Lambda^n \ar[r] \ar[d] & \mbf{BV}_{\E, \le (I,K) } \ar[d]  \\
\Delta^n \ar@{-->}[ur] \ar[r] & \mbf{BV}_{\E,  < (I,K) }
}
$$
commutes.  There is an obstruction 
$$
O_{I,K} ( \Delta^n ) \in \Oo_l(\E,  \Omega^\ast(\Delta^n ) ) .
$$
to the existence of this lift. 

There is a map $\Lambda^n \to \mbf{BV}_{\E, <(I,K) }$.  The obstruction to the lift of this to a map $\Lambda^n \to \mbf{BV}_{\E, \le (I,K) }$ is an element
$$
O_{I,K}( \Lambda^n ) \in \Oo_l ( \E, \Omega^\ast(\Lambda^n) ) .
$$
Because we are given a lift $\Lambda^n \to \mbf{BV}_{\E, \le (I,K) }$, we are given an element
$$
L_{I,K}( \Lambda^n ) \in \Oo_l ( \E, \Omega^\ast(\Lambda^n) ) .
$$
such that $(Q + \d_{DR}) L_{I,K} (\Lambda^n) = O_{I,K} (\Lambda^n)$.

Under the natural map
$$
\Oo_l ( \E, \Omega^\ast(\Delta^n) )  \to \Oo_l ( \E, \Omega^\ast(\Lambda^n) )
$$
the restriction of $O_{I,K}(\Delta^n)$ is $O_{I,K}(\Lambda^n)$. 

The spaces $\Oo_l ( \E, \Omega^\ast(\Delta^n))$ form a simplicial abelian group.  All simplicial groups are Kan complexes; see  \cite{May92}.  This implies that this restriction map is surjective.  It is also a quasi-isomorphism.  (This argument was explained to me by Ezra Getzler).

 It follows immediately that there exists an $L_{I,K}(\Delta^n) \in \Oo_l ( \E, \Omega^\ast(\Delta^n) )$ which satisfies $(Q + \d_{DR}) L_{I,K}(\Delta^n) = O_{I,K}(\Delta^n)$, and which restricts to $L_{I,K}(\Lambda^n)$.  The choice of such an $L_{I,K}(\Delta^n)$ is precisely the choice of a lift of $\Delta^n \to \mbf{BV}_{\E,  < (I,K) }$ to $\mbf{BV}_{\E,  \le (I,K) }$.  The fact that $L_{I,K}(\Delta^n)$ restricts to $L_{I,K}(\Lambda^n)$ implies that the above diagram commutes.

\end{proof}

\section{A local to global principle} 

\label{section local global}
\subsection{Simplicial presheaves and the \v{C}ech complex}

Let $\op{sSets}$ denote the category of simplicial sets.   Recall that a simplicial presheaf on a topological space is simply a presheaf of simplicial sets.   Let $\mc U$ be an open cover of $M$, so that $\mc U$ is a topological space with a surjective open map $\mc U \to M$ which, on every connected component of $\mc U$, is an open embedding.   Let 
$$
\mc U_n = \mc U \times_M  \cdots \times_M \mc U 
$$
where $\mc U$ appears $n+1$ times on the right hand side. As usual, we can associate to $\mc U$ a simplicial space, whose space of $n$-simplices is $\mc U_n$.

Every connected component of $\mc U_n$ is an open subset of $X$.  Let $G$ be a simplicial presheaf on $M$, and let 
$$G(\mc U_n) = \prod_{V \text{ a connected component of } \mc U_n} G(V) . $$
The simplicial sets $G(\mc U_n)$ form a cosimplicial simplicial set. In other words, if $\Delta$ denotes the category of totally ordered finite sets, there is a functor 
\begin{align*}
\check{G}_\mc U : \Delta &\to \op{sSets} \\
\check{G} _\mc U[n] &= G(\mc U_n) 
\end{align*}
where $[n] \in \op{Ob} \Delta$ is the set $\{0,1,\ldots,n\}$. 

Next we will define a simplicial set $\check{C}( \mc U, G )$ which is the \v{C}ech complex for the open cover $\mc U$, with coefficients in $G$.  

For any simplicial sets $X,Y$, let us denote by
$$
\Hom^{\Delta}( X, Y)
$$
the simplicial set of maps $X \to Y$,  so that 
$$
\Hom^{\Delta}( X, Y) [n] = \Hom( X \times \Delta^n, Y) 
$$

For each face map $\phi : \Delta_{n-1} \to \Delta^n$, there are maps
$$
\Hom^\Delta (\Delta^n , G( \mc U_n ) ) \xto{\phi^\ast} \Hom^\Delta (\Delta_{n-1} , G( \mc U_n ) ) \xleftarrow{\phi_\ast} \Hom^\Delta (\Delta_{n-1}  , G( \mc U_{n-1} ) ).
$$
The first arrow is simply composition with $\phi \times \op{Id}$, and the second arrow comes from the cosimplicial set structure on the simplicial sets $G(\mc U_n)$.  

Let 
$$
\check{C}( \mc U, G ) \subset  \prod_n \Hom^\Delta (\Delta^n , G( \mc U_n ) )
$$
be the sub-simplicial set of those elements which are compatible with the face maps.  More precisely, an element $f \in \check{C}( \mc U, G ) [m] $ is a sequence of elements    
$$f_n \in \Hom (\Delta^n \times  \Delta^m , G( \mc U_n ) )$$ for each $n \ge 0$, such that  for each face map $\phi : \Delta_{n-1} \to \Delta^n$, 
$$
\phi^\ast f_n = \phi_\ast f_{n-1} \in  \Hom ( \Delta_{n-1} \times \Delta^m , G( \mc U_n ) ).
$$
Note that we \emph{do not} require compatibility with degeneracy maps.  The simplicial set $\check{C}(\mc U, G)$ is a version of the total complex of the cosimplicial simplicial set defined by $G(\mc U_n)$ (see \cite{BouKan72}).

Similarly, let 
$$
\check{C}_{\le N}( \mc U, G ) \subset \prod_{n = 0}^N \Hom^\Delta (\Delta^n , G( \mc U_n ) )
$$
be the sub-simplicial set of those elements compatible with face maps $\Delta_{n-1} \to \Delta_{n}$, for all $0 \le   n \le N$.  Note that 
$$
\check{C}( \mc U, G ) = \liminv_N \check{C}_{\le N}( \mc U, G )
$$

The simplicial set $\check{C}( \mc U, G )$ is the non-linear \v{C}ech complex for the open cover $\mc U$, with values in the simplicial presheaf $G$.  As with ordinary \v{C}ech cohomology, we will define the simplicial set of derived global sections of $G$ as a limit over all open covers.  That is, let 
$$
\R \Gamma (M , G ) = \op{colim}_{\mc U} \check{C} ( \mc U,G)  .  
$$
\begin{definition}
A simplicial presheaf $G$ is \emph{locally fibrant} if for there exists an open cover $\{U_i\}$ of $M$, such that for all $i$ and all open sets $W \subset U_i$, $G(W)$ is a Kan complex (i.e.\ a fibrant simplicial set). 
\end{definition}

We will now see that our definition of $\R \Gamma(M,G)$ is well-behaved as long as $G$ is locally fibrant.
\begin{lemma}
If $G$ is locally fibrant, then the map
$$
\check{C}_{\le N}( \mc U, G ) \to \check{C}_{\le N-1}( \mc U, G )
$$
is a fibration for sufficiently fine open covers $\mc U$.
\end{lemma}

\begin{proof}
We  can assume that $G(W)$ is fibrant for all open sets $W$ contained inside some element of the cover $\mc U$. 

Suppose we have a map $\Delta^n \to \check{C}_{\le N-1}( \mc U, G )$, and a map $\Lambda^n \to \check{C}_{\le N}( \mc U, G )$ lifting the composition $\Lambda^n \to \Delta^n \to \check{C}_{\le N-1}( \mc U, G )$.  We need to construct a lift $\Delta^n \to \check{C}_{\le N}( \mc U, G )$.

Such a lift is given by a map $\Delta^n \to \Hom^\Delta(\Delta_N , G ( \mc U_N ) )$, in other words, a map $\Delta_N \times \Delta^n \to G( \mc U_N)$.    This map is known on $\partial \Delta_N \times \Delta^n$ and on $\Delta_N \times \Lambda^n$.  Indeed, the definition of $\check{C}_{\le N}( \mc U, G )$ implies that for each face $\Delta_{N-1} \to \Delta_N$, this map must restrict to the given  map $\Delta^n \to \Hom^\Delta( \Delta_{N-1}, G( \mc U_{N-1})$.  On $\Delta_N \times \Lambda^n$, this map must restrict to the given map $\Lambda^n \to \Hom^\Delta(\Delta_N , G ( \mc U_N ) )$.  

Thus, we are left with finding an extension of a map 
$$
\left( \partial \Delta_N \times \Delta^n \right) \amalg_{\partial \Delta_N \times \Lambda^n} \Delta_N \times \Lambda^n \to G ( \mc U_N ) 
$$
to a map 
$$
\Delta_N \times \Delta^n \to G ( \mc U_N ) . 
$$
The map
$$
\left( \partial \Delta_N \times \Delta^n \right) \amalg_{\partial \Delta_N \times \Lambda^n} \Delta_N \times \Lambda^n  \to \Delta_N \times \Delta^n 
$$
is a cofibration and a weak equivalence.    Since $G ( \mc U_N ) $ is fibrant, the desired extension exists.  
\end{proof} 

\begin{lemma}
Let $G, G'$ be a locally fibrant simplicial presheaves on $M$, and let $\mc U$ be a sufficiently fine open cover of $M$.  Let $G \to G'$ be a map of simplicial presheaves such that the maps $G(\mc U_n) \to G'(\mc U_n)$ are all weak equivalence. 

Then the map
$$
 \check{C} ( \mc U,G)   \to  \check{C} ( \mc U,G')    
$$
is a weak equivalence. 

\label{lemma local weak equivalence}
\end{lemma}
\begin{remark}
This lemma is standard in the theory of cosimplicial spaces, see \cite{BouKan72}.
\end{remark}

\begin{proof}
The long exact sequence for the homotopy groups of a fibration implies that the maps $\check{C}_{\le N} ( \mc U,G)   \to  \check{C}_{\le N} ( \mc U,G')$ are weak equivalences for all $N$.  The result follows easily from this. 
\end{proof} 
\begin{lemma}
Let $G, G'$ be locally fibrant simplicial presheaves on $M$.  Let $G \to G'$ be a map of simplicial presheaves such that, for all sufficiently small open balls $B$ in $M$, the map $G (B) \to G'(B)$ is a weak equivalence.

Then the map $\R \Gamma(M, G) \to \R \Gamma(M,G') $ is a weak equivalence.
\end{lemma}
\begin{proof}
Recall that a \emph{good cover} of $M$ is a cover all of whose iterated intersections are open balls.  Every cover of $M$ has a refinement by a good cover.  

If $\mc U \to M$ is a sufficiently fine good cover, then $\check{C}(\mc U, G ) \to \check{C}(\mc U, G')$ is a weak equivalence, by the previous lemma. 

The colimit over open covers $\mc U$ we use to define $\R \Gamma(M,G)$ is a filtered colimit.  Homotopy groups commute with filtered colimits.    It follows immediately that the map $\R \Gamma(M,G) \to \R \Gamma(M,G')$ is a weak equivalence.  
\end{proof} 

\begin{corollary}
Let $G$ be a simplicial presheaf such that for all sufficiently small open balls $U$ in $M$, $G(U)$ is contractible.  Then $\R \Gamma(M,G)$ is contractible.
\label{corollary locally contractible}
\end{corollary}
\begin{proof}
The hypothesis, together with Lemma \ref{lemma local weak equivalence}, implies that the map $\R \Gamma(M,G) \to \R \Gamma(M,\ast)$ is a weak equivalence, where $\ast$ is the constant simplicial presheaf on $M$ with value a point.  For any open cover $\mc U$ of $M$, $\check{C}(\mc U, \ast) = \ast$, so $\R \Gamma(\mc U, \ast) = \ast$.
\end{proof} 

\subsection{Simplicial presheaf of gauge fixing conditions}
In order to construct a simplicial presheaf of gauge fixing conditions, we need some extra data, which is typically present in applications.    

The extra data consists of:
\begin{enumerate}
\item
A smooth, locally trivial fibre bundle $F$ over $M$, whose fibres are diffeomorphic to $\R^k$ for some $k$.
\item
A map $D : F \to \op{Diff}^{\le 1} (E)$ of fibre bundles over $M$.
\item
If $f : U \to  F\mid_U$ is any local section, defined over some open set $U$ in $M$, then 
$$
D(f) : \Gamma_c(U, E) \to \Gamma_c(U,E)
$$
must be an odd operator of square zero, which is self adjoint with respect to the symplectic pairing, and which is such that $[Q, D(f)]$ is a generalised Laplacian.  Here $\Gamma_c$ denotes sections with compact support.
\item
If $f : M \to F$ is a global section of $F$, then $D(f)$ satisfies the axioms of a gauge fixing condition: in addition to the properties already mentioned, this means that there is a direct sum decomposition
$$
\E = \Gamma(M,E) = \Im Q \oplus \Ker [ D(f), Q ] \oplus \Im D(f) . 
$$
\end{enumerate}

It is a little \emph{ad hoc} to require the existence of such a bundle $F$.  It would be better to construct a bundle of local gauge fixing conditions directly from the data $(\E,Q, \ip{\ , \ })$.    However, I couldn't see a clean way of doing this.   In all examples, one is naturally given an $F$ satisfying these properties.     

For example, in Chern-Simons theory, $F$ will be the bundle of metrics on the tangent bundle of $M$.  This satisfies these conditions.  Once we fix some metric, we can identify the fibre of the bundle of metrics on the tangent bundle is $GL(n) / SO(n)$, which is diffeomorphic to $\R^{n(n+1)/2}$.  

Associated to the bundle $F$ we have a simplicial presheaf $\mbf{F}$  on $M$, such that 
$$\mbf{F}(U)[n] = \text{smooth sections of }  F \vert_U \times \Delta^n \to U \times \Delta^n.$$
There is a map
$$
\Gamma(M,\mbf{F} )  \to \mbf{GF} .
$$
The simplicial presheaf $\mbf{F}$ has the following properties.
\begin{lemma}
\label{lemma properties sheaf gf}
\begin{enumerate}
\item
$\mbf{F}$ is locally fibrant. 
\item
Both $\Gamma(M, \mbf{F})$ and $\R \Gamma(M, \mbf{F})$ are contractible.
\item
The sheaves of sets $U \mapsto \mbf{F}(U)[n]$ are soft. 
\end{enumerate}
\end{lemma}

Recall that a soft sheaf $G$ is one such that for all closed sets $C \subset M$, the map $G(M) \to G(C)$ is surjective. 

\begin{proof}
I will give a proof of the first statement, that $\mbf{F}$ is locally fibrant.  The other assertions are straightforward to prove. It suffices to show that for all open subsets $U \subset M$ on which the bundle $F$ is trivial, the simplicial set $\mbf{F}(U)$ is fibrant.

Trivialising $F$ on $U$ gives a diffeomorphism $F \mid_U \iso U \times \R^{k}$.   Thus, the simplicial presheaf $\mbf{F}$ is isomorphic to the simplicial presheaf whose $n$-simplices, on an open set $U$, are smooth maps $U \times \Delta^n \to \R^k$.  Therefore the simplicial set $\mbf{F}(U)$ admits the structure of simplicial abelian groups.  All simplicial abelian groups are Kan complexes, thus $\mbf{F}(U)$ is a Kan complex.

\end{proof}

\subsection{The sheaf of local functionals}
We want to talk about sheaves of solutions to the renormalised quantum master equation.  In order to do this, we need to first define sheaves of local functionals.  Throughout this section, we will work with functions \emph{modulo constants}.  

\begin{definition}
Let $\Oo(E)$ denote the presheaf on $M$ whose value on an open set $U$ is 
$$
\Gamma(U, \Oo(E)) = \prod_{n \ge 1} \Hom ( \Gamma_c(U,E)^{\otimes n}, \C ) _{S_n}. 
$$
Here, we work with functions modulo constants, so we take $n \ge 1$. $\Hom$ denotes continuous linear maps,  the subscript $S_n$ denotes coinvariants, and the subscript $c$ denotes sections with compact support.
\end{definition}
Thus,
$$
\Gamma(M, \Oo(E)) = \Oo(\E)  / \C . 
$$
In a similar way, if $X$ is an auxiliary manifold, there is a presheaf $\Oo(E, \cinfty(X) )$  on $M$,  such that
$$
\Gamma(U, \Oo(E,  \cinfty(X)) = \prod_{n \ge 1} \Hom ( \Gamma_c(U,E)^{\otimes n}, \cinfty(X) ) _{S_n}. 
$$

Recall that  polydifferential operators $\E^{\otimes n} \to \dens(M)$ can be identified with global sections of the infinite rank vector bundle
$$
\Diff(E,\C)^{\otimes n} \otimes \dens(M)
$$
on $M$. Here $\C$ denotes the trivial bundle, $\Diff$ is the sheaf of differential operators between two vector bundles, and tensor products are fibrewise tensor products of vector bundles.  
\begin{definition}
Let $\Oo_l (E)\subset \Oo(E)$ be the sub-presheaf given locally by the image of 
$$
\prod_{n \ge 1} \Gamma(U,  \Diff(E,\C)^{\otimes n} \otimes \dens_M ) \xto{\int}  \prod_{n \ge 1} \Hom ( \Gamma_c(U,E)^{\otimes n}, \cinfty(X) ) _{S_n} .
$$
\end{definition}
Thus, 
$$
\Gamma(M , \Oo_l (E) ) = \Oo_l(\E) / \C .
$$
In a similar way, we can define the sheaf $\Oo_l (E ,  \cinfty(X) )$ of local functionals with values in $\cinfty(X)$, or in $\Omega^\ast(X)$, if $X$ is an auxiliary manifold (with corners).
\begin{lemma}
The presheaf $\Oo_l (E, \Omega^\ast(X))$  is a soft sheaf. 
\end{lemma}

\subsection{Simplicial presheaf of solutions to the renormalised QME}

Proposition \ref{proposition obstruction} shows that if we have $S \in \Oo_l(\E, \Omega^\ast(\Delta^n)\otimes \C [[\hbar]])$, and a smooth family of gauge fixing conditions $\GF$, parametrised by $\Delta^n$, then there are obstructions 
$$O_{(I,K)}(S, \GF )\in \Oo_l(\E, \mscr{A}_{\le 0} \otimes \Omega^\ast(\Delta^n))$$
 to $S$ satisfying the renormalised quantum master equation. These obstructions are homogeneous of order $K$; thus, if $K > 0$, they are elements of $\Gamma(M, \Oo_l (E, \mscr{A}_{\le 0 } \otimes \Omega^\ast(\Delta^n) ))$.  It is easy to see that the $O_{(I,K)}$, when $K > 0$, don't depend on the constant term of $S$.  Thus, the obstructions give us maps 
$$
O_{(I,K)}^{global} : \Gamma(M, \Oo_l(E, \Omega^\ast(\Delta^n )   \otimes \C [[\hbar]] ) \times \mbf{F} [n] ) \to \Gamma(M, \Oo_l(E, \mscr{A}_{\le 0}\otimes \Omega^\ast(\Delta^n) ) ).
$$ 
\begin{lemma}
For all $I \ge 0$, $K > 0$ and $n \ge 0$, there are unique maps of simplicial presheaves 
$$O_{(I,K)} :  \Oo_l(E, \Omega^\ast(\Delta^n )   \otimes \C [[\hbar]] )  \times \mbf{F}  [n] \ \to \Oo_l(E, \mscr{A}_{\le 0} \otimes \Omega^\ast(\Delta^n) )  $$ such that the diagram
$$
\xymatrix{  \Gamma(M,  \Oo_l(E, \Omega^\ast(\Delta^n )   \otimes \C [[\hbar]] ) \times \mbf{F} [n]  )    \ar[r]^>>>>>{O_{(I,K)}^{global}} \ar[d] &   \Gamma(M, \Oo_l(E, \mscr{A}_{\le 0} \otimes \Omega^\ast(\Delta^n) ) )  \ar[d] \\ 
\Gamma(U,  \Oo_l(E, \Omega^\ast(\Delta^n )   \otimes \C [[\hbar]] )   \times \mbf{F} [n])    \ar[r]^>>>>>{O_{(I,K)}(U)}  &   \Gamma(U, \Oo_l(E, \mscr{A}_{\le 0}\otimes \Omega^\ast(\Delta^n) ) )  }.
$$
commutes, for all open sets $U \subset M$. 
\end{lemma}
\begin{proof}
The sheaf $ \Oo_l(E, \Omega^\ast(\Delta^n )   \otimes \C [[\hbar]] )  \times \mbf{F}  [n]$ is soft.   This implies that for all open sets $U \subset M$, and all sections $(S,f) \in \Gamma(U,  \Oo_l(E, \Omega^\ast(\Delta^n )   \otimes \C [[\hbar]] )  \times \mbf{F}  [n] \ ) $, there exists an open cover of $U$ by sets $\{V_i\}$, and global sections $(S_i,f_i )$, such that $(S_i,f_i) \mid_{V_i} = (S,f) \mid_{V_i}$.   

This implies uniqueness.   For existence, we need to show that $O_{(I,K)}^{global}(S_i,f_i) \mid_{V_i}$ only depends on $(S_i,f_i) \mid {V_i}$. This follows from the third clause of Proposition \ref{proposition obstruction}. 

\end{proof}

\begin{lemma}

Let $S \in \Gamma(U,   \Oo_l(E, \Omega^\ast(\Delta^n )   \otimes \C [[\hbar]] ) \times \mbf{F}  [n])$ be such that $O_{(i,k)}(S)=  0$ for all $(i,k) < (I,K)$.  Then \begin{enumerate}
\item
$O_{(I,K)}(S)$ is independent of $\eps$, so it is a section of $\Oo_l(E, \Omega^\ast(\Delta^n ) )$. 
\item
$$
(Q + \d_{DR} )  O_{(I,K)}(S) = 0. 
$$
\item
If $S'\in \Gamma(U, \Oo_l(E, \Omega^\ast(\Delta^n )  )$ is homogeneous of degree $K$, then 
$$
O_{(I,K)} ( S + \hbar^I S' ) = O_{(I,K)} (S) + (Q + \d_{DR} ) S'.
$$
\end{enumerate}
\end{lemma}
\begin{proof}
The proof is identical to that of Lemma \ref{closed local obstruction}. 
\end{proof} 

\begin{definition}
Let $\mbf{BV}$ be the simplicial presheaf on $M$ such that   $\Gamma(U, \mbf{BV}[n] )$ is the set of $$(S,f) \in \Gamma(U,  \Oo_l(E, \Omega^\ast(\Delta^n )   \otimes \C [[\hbar]] )  \times \mbf{F}[n])$$  such that $$O_{(I,K)}(S,f) = 0$$ for all $I \ge 0,K > 0$.

Let $\mbf{BV}_{\le (I,K)}$  be the simplicial presheaf such that $\Gamma(U, \mbf{BV}[n] )$ is the set  
$$(S,f) \in \Gamma(U,  \Oo_l(E, \Omega^\ast(\Delta^n )   \otimes \C [[\hbar]] )  \times \mbf{F}[n]  )$$ 
such that
\begin{enumerate}
\item
$S_{(i,k)} = 0$ if $(i,k) > (I,K)$, where, as usual, $S_{(i,k)}$ is the part homogeneous of degree $(i,k)$ as a function of $(\hbar, e \in \Gamma_c(U,E) )$.   
\item
$O_{(i,k)}(S,f) = 0$ for $(i,k) \le (I,K)$.
\end{enumerate}

Let $\mbf{BV}_{< (I,K)}$  be the simplicial presheaf such that $\Gamma(U, \mbf{BV}[n] )$ is the set  
$$(S,f) \in \Gamma(U, \Oo_l(E, \Omega^\ast(\Delta^n )   \otimes \C [[\hbar]] )  \times \mbf{F}[n] )$$ 
such that
\begin{enumerate}
\item
$S_{(i,k)} = 0$ if $(i,k) \ge (I,K)$.
\item
$O_{(i,k)}(S,f) = 0$ for $(i,k) < (I,K)$.
\end{enumerate}

\end{definition}

\begin{lemma}
$\mbf{BV}$ is a locally fibrant simplicial presheaf.
\end{lemma}
\begin{proof}
It suffices to show that  $\Gamma ( U , \mbf{BV}) $ is a fibrant simplicial set, for all open subsets $U$ of $M$ on which the bundle $F$ is trivial. Since we know $\Gamma(U, \mbf{F})$ is fibrant, it suffices to show that the map $\Gamma(U, \mbf{BV} ) \to \Gamma(U, \mbf{F})$ is a fibration. The proof is identical to that of Theorem \ref{theorem fibration}.
\end{proof} 

This implies that our definition of $\R \Gamma (M ,\mbf{BV})$ is well-behaved.

\begin{definition}
The operator $Q$ on the sheaf of sections of $E$ is called \emph{triangular} if 
$E$ admits a decomposition as a direct sum of vector bundles  
$$E = E_0 \oplus E_1 \oplus \cdots \oplus E_n$$ such that the operator $Q$ takes sections of $E_i$ to sections of $\oplus_{j > i} E_j$.
\end{definition}

\begin{thmD}
Suppose $Q$ is triangular.
Then the natural map
$$
\Gamma(M, \mbf{BV} ) \to \R \Gamma(M, \mbf{BV} ) 
$$
is a weak equivalence. 
\end{thmD}
This theorem I will call the ``local to global principle''.  This result allows one to construct global solutions to the quantum master equation from local local solutions. The local solutions will typically be constructed in explicit model cases, like $\R^n$ with the flat metric, as we will see below for the case of Chern-Simons theory.

\subsection{Proof of theorem D} 
There is a simplicial presheaf on $M$ whose presheaf of $n$-simplices is
$$
\Gamma(U, \Oo_l(E, \Omega^\ast(\Delta^n ) ) .
$$
Let us denote this simplicial presheaf by $\Oo_l^{\Delta}(E)$.  Note that
$$
\mbf{BV} \subset \Oo_l^{\Delta}(E) \times \mbf{F} 
$$
is a sub-simplicial presheaf. 

\begin{lemma}
Let $X$ be a simplicial set, and let $S_X : X \to \R \Gamma(M, \mbf{BV}_{< (I,K)})$ be a map. Then there is an obstruction 
$$
O(S_X) \in \op{Hom} ( X, \R \Gamma(M, \Oo_l^{\Delta} (E) ) )
$$
which is closed, $(Q + \d_{DR} )O(S_X) = 0$, and homogeneous of degree $K$.  The simplicial set of lifts $X \to \R \Gamma(M, \mbf{BV}_{\le (I,K)})$ is isomorphic to the presheaf of maps $S_X' : X \to \R \Gamma(M, \Oo_l^{\Delta}(E))$, homogeneous of degree $K$, satisfying
$$
(Q + \d_{DR} ) S'_X = O(X) . 
$$
\end{lemma}
\begin{proof}
Fix an open cover $\mc U$ of $M$.  First, we will prove the corresponding statement for lifts of maps
$$
S_X : X \to \check{C}(\mc U, \mbf{BV}_{<(I,K)}).
$$
Such a map is given by maps
$$
S_X^n : X \times \Delta^n \to \Gamma(\mc U_n,  \mbf{BV}_{<(I,K)} )
$$
for all $n$, which are compatible with the face maps $\Delta^m \to \Delta^n$, $\mc U_m \to \mc U_n$, for $m < n$. 

The obstructions to lifting $S_X^n$ to $\mbf{BV}_{\le (I,K)}$ are elements
$$
O(S_X, \mc U_n ) : X \times \Delta^n \to \Gamma(\mc U_n, \Oo_l^\Delta(E) ).
$$
Naturality of these obstructions implies that they fit together into an element
$$
O(S_X, \mc U) : X \to \check{C}(\mc U, \Oo_l^\Delta(E) ).
$$

There is a bijection between lifts of each $S_X^n$ and ways of making the obstruction $O(S_X, \mc U_n)$ exact.  Thus, there is a corresponding bijection between lifts of $S_X$ and ways of making the total obstruction $O(S_X, \mc U)$ exact. Of course, the map $S' : X \to \check{C}(\mc U, \Oo_l^\Delta(E) )$ making $O(S_X, \mc U)$ exact must be homogeneous of degree $K$.

The obstruction $O(S_X, \mc U)$ is natural with respect to restrictions to finer open covers $\mc V$.  Thus, taking the direct limit, we get an element
$$
O(S_X) = \limdir_{\mc U} O(S_X, \mc U) \in \limdir_{\mc U} \check{C}(\mc U, \Oo_l^\Delta(E) ) = \R \Gamma ( M , \Oo_l^\Delta(E) ).
$$
Again, naturality of the obstruction implies that lifts $X \to \R \Gamma(M, \mbf{BF}_{\le (I,K)})$ correspond to ways of making $O(S_X)$ exact.
\end{proof} 
\begin{proposition}
The map
$$
\R \Gamma (M, \mbf{BV}_{\le (I,K)} ) \to \R \Gamma ( M, \mbf{BV}_{< (I,K)} )
$$
is a fibration.
\end{proposition}
\begin{proof}
The proof is exactly the same as that of Proposition \ref{proposition fibration},
using  the obstruction to lifting maps  to $ \R \Gamma ( M, \mbf{BV}_{< (I,K)} )$ to maps to  $\R \Gamma (M, \mbf{BV}_{\le (I,K)} )$, \end{proof} 

Now we will prove the theorem by induction. The starting point in the induction is the statement that the map
$$
\Gamma(M, \mbf{BV}_{(0,0)} ) \to \R \Gamma (M, \mbf{BV}_{(0,0)})
$$
is a weak equivalence. Note that $\mbf{BV}_{(0,0)} = \mbf{F}$.  Lemma \ref{lemma properties sheaf gf} implies that the map $\Gamma(M, \mbf{F}) \to \R \Gamma(M, \mbf{F} )$ is a weak equivalence. 

We will assume by induction that the map
$$
\Gamma(M, \mbf{BV}_{< (I,K)} )\to \R \Gamma (M, \mbf{BV}_{<(I,K)} )
$$
is a weak equivalence.

Consider the diagram
$$
\xymatrix{  \Gamma (M, \mbf{BV}_{\le (I,K)} ) \ar[r] \ar[d]^{\pi} &\R \Gamma ( M, \mbf{BV}_{\le (I,K)} )  \ar[d]^{\R \pi} \\ 
 \Gamma (M, \mbf{BV}_{< (I,K)} ) \ar[r] &  \R \Gamma ( M, \mbf{BV}_{< (I,K)} ) .  } 
$$
By Proposition \ref{proposition fibration}, 	the map $\pi$ is a fibration, and we have just seen that $\R \pi$ is also a fibration.  By induction, the bottom horizontal arrow is a weak equivalence.    By considering the long exact sequence of homotopy groups of a fibration, it suffices to show that the fibres of $\pi$ and of $\R \pi$ are weakly homotopy equivalent. 

Thus, fix a base point
$$
x_0 \in  \Gamma (M, \mbf{BV}_{< (I,K)} )[0] .
$$
It remains to show that the map
$$
\pi^{-1}(x_0) \to \R \pi^{-1} (x_0)
$$
is a weak equivalence.

There is an obstruction
$$
O_{(I,K)}(x_0) \in \Gamma(M, \Oo_l(E) )
$$
to lifting $x_0$ to $\Gamma (M, \mbf{BV}_{\le (I,K)} ) $.  This is a closed element, $Q O_{(I,K)}(x_0) = 0$. 

We can identify $\pi^{-1}(x_0)$ with the simplicial set whose $n$-simplices are elements of 
$$\phi \in \Gamma(M, \Oo_l^{\Delta}(E) ) [n] = \Gamma(M, \Oo_l(E, \Omega^\ast(\Delta^n) ) )$$
which are homogeneous of degree $K$ as functionals on $\Gamma(M,E)$, and which satisfy
$$
(Q + \d_{DR}) \phi = O_{(I,K)}(x_0).
$$ 

Similarly, we can identify $\R \pi^{-1}(x_0)$ with the simplicial set whose $n$-simplices are
$$
\phi \in  \R \Gamma(M, \Oo_l^{\Delta}(E)  ) [n]
$$
which are homogeneous of degree $K$ and satisfy $(Q + \d_{DR}) \phi = O_{(I,K)}(x_0)$.

\begin{lemma}
The natural map of chain complexes
$$
(\Gamma(M, \Oo_l(E) ) , Q)  \to  ( \R \Gamma(M, \Oo_l^{\Delta}(E) )   [0] , Q + \d _{DR} )  
$$
is a quasi-isomorphism. 
\label{lemma linear cech}
\end{lemma}

\begin{proof}
It suffices to show that for all open covers $\mc U$ of $M$, the natural map 
$$
( \Gamma(M, \Oo_l(E) ) , Q) \to ( \check{C}(\mc U, \Oo_l^{\Delta} (E) ) [0] , Q + \d_{DR} )
$$
is an isomorphism.

 The triangular nature of the differential $Q$ on $E$ implies that the complex of sheaves $\Oo_l(E)$ on $M$ admits a grading,  labelled by non-negative integers, such that $Q$ is of strictly negative degree.   This grading is inherited by all auxiliary complexes, such as $\check{C}(\mc U, \Oo_l^{\Delta} (E) ) [0]$ and $\Gamma(M, \Oo_l(E) ) $.  It follows from a spectral sequence argument that it suffices to show that the map of complexes
$$
( \Gamma(M, \Oo_l(E) ) , 0) \to ( \check{C}(\mc U, \Oo_l^{\Delta} (E) ) [0] ,  \d_{DR} )
$$
is a quasi-isomorphism, where the complex on the left hand side now has zero differential.

The complex $\check{C}(\mc U, \Oo_l^{\Delta} (E) ) [0]$, with the differential $\d_{DR}$, is quasi-isomorphic to the ordinary  \v{C}ech complex with coefficients in sheaf $\Oo_l(E)$ (with zero differential).  This sheaf is soft. Thus, its  \v Cech cohomology is the same as its global sections, which is the statement that the map above is a quasi-isomorphism.

\end{proof}

\begin{corollary}
$\pi^{-1}(x_0)$ is empty if and only if $\R \pi^{-1}(x_0)$ is empty.
\end{corollary}
\begin{proof}
Indeed, $\pi^{-1}(x_0)$ is empty if and only if the cohomology class
$$
[\phi ] \in H^\ast (   \Gamma(M, \Oo_l(E) ), Q ) 
$$
is non-zero, and similarly for $\R \pi^{-1}(x_0)$. 
\end{proof}

Let us use the temporary notation $\Oo_{l,K}(E) \subset \Oo_{l}(E)$ for the sub-sheaf of functionals which are homogeneous of degree $K$. Also let $\Oo_{l,K}^\Delta(E) \subset \Oo_{l}^\Delta(E)$  be the sub-simplicial presheaf of $\Oo_{l}^\Delta(E)$ consisting of those functionals  which are homogeneous of degree $K$. Thus, a section of the presheaf of $n$-simplices of $\Oo_{l,K}^\Delta(E)$ is given by an element of $\Gamma(U, \Oo_{l,K}(E, \Omega^\ast(\Delta^n)) )$.   $\Oo_{l,K}^{\Delta}(E)$ is a simplicial presheaf of $\Z/2$-graded complexes, with differential $Q + \d_{DR}$.

Let us suppose that $\pi^{-1}(x_0)$ is non-empty.  (Otherwise, we are done, as $\R \pi^{-1}(x_0)$ will also be empty).  Then by picking a base point in $\pi^{-1}(x_0)$, we can identify $\pi^{-1}(x_0)$ with the simplicial set whose $n$-simplices are closed, even elements of $\Gamma(M, \Oo_{l,K}^\Delta(E) )[n]$.   Let us denote this simplicial set by $\Gamma(M, \Oo_{l,K}^\Delta(E) )^{closed, even}$.

By taking the corresponding base point in $\R \pi^{-1}(x_0)$, we can identify $\R \pi^{-1}(x_0)$ with the simplicial set whose $n$-simplices are closed even elements of $\R \Gamma(M, \Oo_{l,K}^\Delta(E) )[n]$.  Let us denote this simplicial set by $\R \Gamma(M, \Oo_{l,K}^\Delta(E) )^{closed, even}$.

Thus, to complete the proof of the theorem, it suffices to show that the map
$$
\Gamma(M, \Oo_{l,K}^\Delta(E) )^{closed, even} \to \R \Gamma(M, \Oo_{l,K}^\Delta(E) )^{closed, even}
$$
is a weak equivalence.  This will follow from the following lemma, combined with Lemma \ref{lemma linear cech}.
\begin{lemma}
The homotopy groups of the simplicial sets $\Gamma(M, \Oo_{l,K}^\Delta(E) )^{closed, even}$ and $\R \Gamma(M, \Oo_{l,K}^\Delta(E) )^{closed, even}$ are given by
\begin{align*}
\pi_n \left(  \Gamma(M, \Oo_{l,K}^\Delta(E) )^{closed, even} \right) &=
H^{\abs{n}}   (\Gamma(M, \Oo_{l,K} ), Q ) \\
\pi_n \left(  \R \Gamma(M, \Oo_{l,K}^\Delta(E) )^{closed, even} \right)  &= H^{\abs{n}}   (\R \Gamma(M, \Oo_{l,K}^{\Delta} )[0], Q + \d_{DR} ).
\end{align*}
On the right hand side we are using $\Z/2$-graded cohomology groups. These equations are true for all $n \ge 0$. 
\end{lemma}
This is just a $\Z/2$-graded version of the Dold-Kan correspondence. 
\begin{proof}

The $n$-simplices of $\Gamma(M, \Oo_{l,K}^\Delta(E) )^{closed, even}$ are closed and even elements of the complex  $\Omega^\ast(\Delta^n)\otimes \Gamma(M, \Oo_{l,k}(E) )$, where we are using a completed tensor product.  Similarly, the $n$-simplices of $\R \Gamma(M, \Oo_{l,K}^\Delta(E) )^{closed, even}$ are given by closed even elements of $\Omega^\ast(\Delta^n) \otimes \R \Gamma ( M , \Oo_{(l,K) } ^\Delta ) [0]$, where we are using, again, an appropriate completed tensor product.

Let $(V, \d_V)$ be one of these two $\Z/2$-graded complexes, so
$$
(V, \d_V) = \begin{cases}
( \Gamma(M, \Oo_{l,k}(E), Q )  &  \\
( \R \Gamma ( M , \Oo_{(l,K) } ^\Delta ) [0], Q + \d_{DR} ) &
\end{cases}
$$
Let $V^\Delta$ be the simplicial set whose $n$-simplices are closed even elements of $V \otimes \Omega^\ast(\Delta^n)$, where we use an appropriate completed tensor product.  We need to calculate the homotopy groups of $V^\Delta$.

As $V^\Delta$ is a Kan complex, we can calculate these homotopy groups as homotopy classes of maps from spheres to $V^\Delta$.  Let us take $0$ as a base point in $V^\Delta$.  Let $S^n$ be a simplicial set representing the $n$ sphere, and let $p \in S^n$ be a  base point.    A base point preserving map $f : S^n \to V^\Delta$ is given by an even element 
$$
\omega_f \in V \otimes \Omega^\ast(S^n , p ) .
$$
which is closed, $(\d_V + \d_{DR} ) \omega = 0$.

Thus, $\omega_f$ gives a class 
$$[\omega_f] \in H^{0 }  (V \otimes \Omega^\ast(S^n , p )) = H^{\abs{n} } (V) .$$
(Here we are using $\Z/2$-graded cohomology groups).

We need to show that if the map $f : S^n \to V^\Delta$ changes by a homotopy, then the class $[\omega_f]$ remains unchanged, and conversely, if we change $\omega_f$ to something cohomologous, then the corresponding map $S^n \to V^\Delta$ changes by a homotopy.

If $f,g : S^n \to V^\Delta$ are base point preserving maps, a base point preserving homotopy between them is given by an element $\eta \in V \otimes \Omega^\ast ( S^n \times [0,1], p \times [0,1] )$, which is even and closed, and which restricts to $\omega_f$ and $\omega_g$ at $0$ and $1$.

Since 
$$H^\ast (V \otimes \Omega^\ast ( S^n \times [0,1], p \times [0,1] ) ) = H^{\abs{n}}(V) $$
the existence of such a homotopy implies that $[\omega_f] = [ \omega_g]$.

Conversely, suppose $\phi \in V \otimes \Omega^\ast(S^n , p )$ is an odd element such that 
$$
(\d_V + \d_{DR}) \phi = \omega_f - \omega_g . 
$$
Then if we let 
$$
\eta = \omega_f + t (\omega_g - \omega_f ) + \phi \d t \in V \otimes \Omega^\ast ( S^n \times [0,1], p \times [0,1] )
$$
then $\eta$ is an even element satisfying $(\d_V + \d_{DR}) \eta = 0$, so $\eta$ defines a homotopy between $f$ and $g$. 

\end{proof} 
This completes the proof of theorem D. 

\section{Local computations in Chern-Simons theory }
We will apply the theory developed in the previous section to construct a canonical up to homotopy solution of the renormalised quantum master equation in the case of Chern-Simons theory.   

In this section, therefore, let $M$ be a manifold of dimension 
$$\dim M = n \ge 2.$$ 
Let $\mf{g}$ be a locally trivial sheaf of Lie algebras \footnote{We could also use  $L_\infty$ algebras, and the construction works in more or less the same way. }  over $\C$  on $M$, with a non-degenerate invariant symmetric pairing of the opposite parity to $n$, with values in the constant sheaf $\C$.

Let
$$
\E  = \Omega^\ast(M , \mf g )[1] . 
$$
Let $E$ denote the corresponding vector bundle on $M$, so $\Gamma(M,E)  = \E$.    Let $S_{CS} \in \Oo_l(\E)$ denote the Chern-Simons action.  

\begin{theorem}
There is a canonical, up to a contractible choice, element 
$$S \in \Oo_l(\E , \C[[\hbar]] )/ \C[[\hbar]]$$ 
which solves the renormalised quantum master equation modulo constants, and is of the form
$$
S = S_{CS } + \hbar S^{(1)} + \hbar^2 S^{(2)} +\cdots
$$
\label{theorem chern simons}
\end{theorem}
The proof was sketched in the introduction; the idea is to show that locally, in a flat metric, $S_{CS}$ satisfies the renormalised QME.  These local solutions together will give a global solution in a curved metric. 

\subsection{Simplicial presheaves of metrics}
The bundle on $M$ playing the role of $F$ in section \ref{section local global} is the bundle $\op{Met} \to M$, whose fibres are metrics on the tangent space of $M$.  As in section \ref{section local global}, we have a simplicial presheaf associated to $\op{Met}$,  defined by saying $\Gamma(U, \mbf{Met} [d])$ is the set of smooth family of metrics on $U$ parametrised by $\Delta^d$.   

Let $\mbf{FMet} \subset \mbf{Met}$ be the sub-simplicial presheaf given by families of flat metrics on $U$.  
\begin{lemma}
If $U$ is a ball in $M$, then $\Gamma(U, \mbf{FMet} )$ is contractible. \label{lemma flat metrics locally contractible}
\end{lemma}
\begin{proof}
Pick an isomorphism between $U$ and $\{x \in \R^{n} \mid \norm{x} < 1 \}$.  The standard flat metric on the ball in $\R^{n}$ gives a metric on $U$, which we call $g^0$. 

Let $x_1,\ldots, x_{n}$ denote our coordinates on $U$.  If $g$ is a smooth family of flat metrics on $U$, parametrised by $\Delta^d$, then, in these coordinates, we can write $g$ as 
$$
\sum_{1 \le i,j \le n} g_{ij}  ( x , \sigma ) \d x_i \otimes \d x_j .
$$
Here $g_{ij}(x,\sigma)$ is smooth in the $\sigma$ variables, and smooth in the $x$ variables.
To show that the simplicial set of such $g$ is contractible, it suffices to construct a simplicial map 
$$\Gamma(U, \mbf{FMet} ) \times \Delta_1 \to \Gamma(U, \mbf{FMet} ),$$ which is a simplicial homotopy between the identity on $\Gamma(U, \mbf{FMet} )$ and the projection onto the flat metric $g^0 = \sum \d x_i \otimes \d x_i$.   A $d$-simplex in $\Gamma(U, \mbf{FMet} ) \times \Delta_1$ is of course a $d$-simplex in $\Gamma(U, \mbf{FMet} )$  and a $d$-simplex on $\Delta_1$.  A $d$-simplex in $\Delta_1$ is given by a non-decreasing map of sets $\{0,1,\ldots,d\} \to \{0,1\}$.   Such a $d$-simplex yields an affine linear map $\Delta^d \to \Delta_1$.  

Let $A \subset \Gamma(U, \mbf{FMet} )$ denote the simplicial set of constant metrics, i.e. those for which $g_{ij}(x,\sigma)$ is independent of $x \in U$.  Of course, this is just the simplicial set whose $0$-simplices are positive-definite inner products on the vector space $\R^{n}$, and whose $d$-simplices are smooth families of such.

First we construct a deformation retraction of $\Gamma(U, \mbf{FMet} )$ onto $A$. This is given by a map $\Phi : \Gamma(U, \mbf{FMet} ) \times \Delta_1 \to \Gamma(U, \mbf{FMet} )$, defined as follows.  Suppose we have a $d$-simplex in $\Gamma(U, \mbf{FMet} ) \times \Delta_1$, corresponding to $g \in \Gamma(U, \mbf{FMet} )[d]$ and an affine linear map $p : \Delta^d \to \Delta_1$. If, as above, $g = \sum g_{ij}(x,\sigma) \d x_i \otimes \d x_j$,   we define $\Phi(g,p)$ by
$$
\Phi(g,p ) = \sum g_{ij} ( p(\sigma) x, \sigma  ) \d x_i \otimes \d x_j . 
$$
$\Phi$ is easily seen to give the required homotopy; when we restrict to $\Gamma(U, \mbf{FMet} ) \times \{1\} \subset \Gamma(U, \mbf{FMet} ) \times \Delta_1$, $\Phi$ is the identity, and when we restrict to $\Gamma(U, \mbf{FMet} )\times \{0\}$, $\Phi$ is a projection onto $A$.

Next we need to show that the simplicial set $A$ is contractible.  Note that the space of positive-definite inner products on $\R^n$ is $\op{GL}(n,\R) / \op{O}(n.\R)$.  A simplex in $A$ is a continuous map $\Delta^d \to \op{GL}(n,\R) / \op{O}(n,\R)$, which is smooth on the top simplices of some barycentric subdivision of $\Delta^d$.  Thus, to show $A$ is contractible, it suffices to construct a  smooth homotopy equivalence between $\op{GL}(n,\R) / \op{O}(n,\R)$ and a point.  This is easy to do using (for instance) the Gram-Schmidt procedure.

\end{proof} 

The simplicial presheaf $\mbf{FMet}$ is not necessarily locally fibrant.  Thus, our definition of $\R \Gamma$ is not well-behaved when applied to $\mbf{FMet}$.  Therefore, we will instead consider $\R \Gamma$ applied to a modification of $\mbf{FMet}$. 

Let $\op{Ex}^\infty$ denote Kan's fibrant replacement functor for simplicial sets.  If $X$ is a simplicial set, then $\op{Ex}^\infty(X)$ is a fibrant simplicial set equipped with a map to $X$ which is a weak equivalence.  Define $\op{Ex}^\infty \mbf{FMet}$ by saying that 
$$\Gamma (U, \op{Ex}^\infty \mbf{FMet } ) = \op{Ex}^\infty \Gamma ( U, \mbf{FMet }).$$
Since $\op{Ex}^\infty \mbf{FMet}$ is locally fibrant, $\R \Gamma ( M, \op{Ex}^\infty \mbf{FMet})$ is well-behaved. 
This is the correct version of derived global sections of $\mbf{FMet}$.  

\begin{corollary}
$\R \Gamma (M , \op{Ex}^\infty \mbf{FMet} ) $ is contractible.
\end{corollary}
This is immediate from Corollary \ref{corollary locally contractible} and Lemma \ref{lemma flat metrics locally contractible}.

We apply our local to global principle as follows.  We would like to construct a homotopy point of $\Gamma(M, \mbf{BV})$.  Any map $\mbf{FMet} \to \mbf{BV}$ of simplicial presheaves would yield such a homotopy point; indeed, such a map induces a map
$$
\R \Gamma(M, \op{Ex}^\infty \mbf{FMet}) \to \R \Gamma(M,  \op{Ex}^\infty \mbf{BV}) .
$$
We know that $\R \Gamma(M, \op{Ex}^\infty \mbf{FMet}) $ is contractible. Also, $\Gamma(M,\mbf{BV})$ is weakly equivalent to $\R \Gamma(M, \mbf{BV})$, and therefore to $\R \Gamma ( M ,\op{Ex}^\infty \mbf{BF} )$.     Thus, a map $\mbf{FMet} \to \mbf{BV}$ yields a point (up to contractible choice) of $\Gamma(M, \mbf{BV})$.  

\begin{theorem}
Let $U \subset M$ be an open subset, and let
$$
g \in \Gamma(U, \mbf{FMet}[d] ) 
$$
be a smooth family of flat metrics on $U$ parametrised by the $d$-simplex $\Delta^d$.  Then the usual Chern-Simons action $S_{CS}$ satisfies the renormalised quantum master equation on $U$, and so defines an element
$$
S_{CS} \in \Gamma(U , \mbf{BV}[d] ) .
$$
\label{theorem chern simons local calculation}
\end{theorem}
This theorem gives a map 
$$
\mbf{FMet} \to \mbf{BV}
$$
 of simplicial presheaves.  Thus, once we have proved this result, we will have proved Theorem \ref{theorem chern simons}. 
\subsection{Proof of Theorem \ref{theorem chern simons local calculation}}
We want to show that, locally in a flat metric, all the obstructions $O_{(I,K)}(S_{CS})$ vanish, for $I \ge 0, K > 0$. 

Since everything is local, it suffices to work with small open sets $U$ in $M$, with a smooth family of flat metrics parametrised by $\Delta^d$.

Let us trivialise our flat bundle of Lie algebras $\mf g$ on $U$.  The space of fields is the $\Omega^\ast(\Delta^d)$-module
$$
\Omega^\ast(U) \otimes \Omega^\ast(\Delta^d) \otimes \mf g
$$
with the differential $\d_{DR}^U + \d_{DR}^{\Delta^d}$, and the gauge fixing condition $\d^\ast$.  The operator $\d^\ast$ depends on the coordinates in $\Delta^d$.   The Hamiltonian $H$ is
$$
H = [\d_{DR}^U + \d_{DR}^{\Delta^d}, \d^\ast ] .
$$

Recall that $n = \dim M$. Let us give $\R^n$ the standard flat metric.  It is easy to check that one can find an open subset $V \subset \R^n\times \Delta^d$, and an isomorphism $U \times \Delta^d \to V$, such that the diagram
$$
\xymatrix{ U \times \Delta^d \ar[r] \ar[d] & V \ar[dl] \\ 
\Delta^d &   }
$$
commutes, and which is compatible with the metrics along the fibres of the projection maps to $\Delta^d$.

The isomorphism chosen above gives an isomorphism
$$
\Omega^\ast(U) \otimes \Omega^\ast(\Delta^d) \otimes \mf g \iso \Omega^\ast(V)
$$
of differential graded $\Omega^\ast(\Delta^d)$-modules.  The differential on both sides is simply the de Rham differential on $U \times \Delta^d$, or on $V$. Because the metric along the fibres of the map $U\times \Delta^d \to \Delta^d$ corresponds to that along the fibres of the map $V \to \Delta^d$, this isomorphism also takes the operator $\d^\ast$ on the left hand side to the corresponding one on the right hand side.

The propagator is constructed from the de Rham differential on $U \times \Delta^d$, and the operator $\d^\ast$. Thus,  it suffices to show that the obstructions $O_{(I,K)}(S_{CS})$ vanish when we work with $V$ instead of $U \times \Delta^d$. Since the obstruction is local, it suffices to work in $\R^n \times \Delta^d$.

We have reduced to the case of $\R^n$ with the constant family of flat metrics.   Thus, it suffices to show the following.
\begin{theorem}
Let $S_{CS}$ be the usual Chern-Simons action on $\R^n$ with the flat metric.  Then
\begin{enumerate}
\item
There are no counter-terms.
\item
$S_{CS}$ satisfies the renormalised quantum master equation.
\end{enumerate}
\end{theorem}
This result implies, in particular, that the construction is independent of the choice of renormalisation scheme.  The choice of renormalisation scheme is only involved in constructing the bijection between the set of local functionals $S\in \Oo_l(\E, \C[[\hbar]])$ and the set of systems of effective actions satisfying the renormalisation group equation. However, because there are no counter-terms, the system of effective actions $\{\Gamma ( P (0,T) , S_{CS} )\}$ associated to the Chern-Simons action $S_{CS}$ is independent of the choice of renormalisation scheme.

The proof of this theorem will follow closely the results of Kontsevich \cite{Kon94, Kon03a}.  In particular, we make use of the compactifications of configuration spaces used in these papers, and we rely on a certain vanishing theorem proved by Kontsevich in \cite{Kon94} when $n \ge 3$, and in \cite{Kon03a} when $n = 2$.

Let $\op{Conf}_m(\R^n)$ denote the space of $m$ distinct ordered points in $\R^n$.    There is a partial compactification $\br{\op{Conf}}_m(\R^n)$ of $\op{Conf}_m(\R^n)$ with the following properties.  
\begin{enumerate}
\item
$\br{\op{Conf}}_2(\R^n)$ is the real blow-up along the diagonal of $\R^n \times \R^n$.
\item
$\br{\op{Conf}}_m(\R^n)$ admits a proper map to $\R^{nm}$ such that the diagram
$$
\xymatrix{ \op{Conf}_m(\R^n) \ar[r] \ar[d] & \br{\op{Conf}}_m(\R^n)  \ar[dl] \\ 
\R^{nm} &  }
$$
commutes.
\item
For any  $i,j$ with $1 \le i,j \le m$ and $i \neq j$,  there is a projection map
$$
\pi_{i,j} : \br{\op{Conf}}_m(\R^n) \to \R^{2n}
$$
given by forgetting all factors except $i$ and $j$. This map lifts to $\br{\op{Conf}}_2(\R^n)$.
\end{enumerate}
The partial compactification we use was constructed in \cite{Kon94,Kon03a}. This is the real version of the Fulton-MacPherson compactification of configuration spaces of algebraic varieties \cite{FulMac94}.

Let
$$
\E = \Omega^\ast(\R^n) \otimes \mf g [1]
$$
be the space of fields, and let
$$
\E_c = \Omega^\ast_c(\R^n) \otimes \mf g [1]
$$
be the space of compactly supported fields.

Let
$$
P(\eps,T ) = \int_{\eps}^T \d^\ast K_t \d t \in \E \otimes \E .
$$
\begin{lemma}
Up to an overall constant,
$$
P(\eps,T) = \left ( \int_{ \norm{x - y}^2 / 4 T }^{\norm{x - y}^2 / 4 \eps }  u^{n/2 - 1 } e^{- u }  \d u \right)   \pi^\ast \op{Vol}_{S^{n-1}} \otimes  I_{\mf g} 
$$
where
\begin{enumerate}
\item
$$
\pi : \op{Conf}_2(\R^n) \to S^{n-1}
$$
is the projection 
$$(x,y) \to (x - y)  / \norm{x - y} .$$
\item
$ \op{Vol}_{S^{n-1}}$ is the standard volume form on $S^{n-1}$, given by the formula
$$
 \op{Vol}_{S^{n-1}} = \norm{z}^{-n} \sum_{i = 1}^{n} (-1)^{i}  z_i  \d z_1 \cdots \what{ \d z_i } \cdots  \d z_n .
$$
\item
$$
I_{\mf g} \in \mf g \otimes \mf g
$$
is the tensor dual to the pairing on $\mf g$. 
\end{enumerate}
\end{lemma}
\begin{proof}
The proof is an explicit calculation.  The sign conventions we are using imply that the heat kernel is
$$
K_t =  C t^{-n/2 }  e^{- \norm{x-y}^2 / 4t } \d (x_1 - y_1) \wedge \cdots \wedge \d (x_n - y_n) \otimes I_{\mf g}  
$$
where $C$ is a certain normalising constant. 
It follows that 
\begin{multline*}
\d^\ast K_t  =  C' t^{-n/2 - 1 } e^{- \norm{x-y}^2 / 4t } \sum_{i = 1}^{n} (-1)^{i}  (x_i - y_i) \d (x_1 - y_1)  \cdots  \what{ \d (x_i -  y_i) } \\  \cdots \d (x_n - y_n )  \otimes  I_{\mf g} 
\end{multline*}
where $C'$ is a constant.
Now, 
$$
\int_\eps^T  t^{-n/2 - 1 } e^{- \norm{x-y}^2 / 4t } \d t = C'' \norm{x - y}^{-n}  \int_{ \norm{x - y}^2 / 4 T }^{\norm{x - y}^2 / 4 \eps }  u^{n/2 - 1 } e^{- u }  \d u 
$$
for some constant $C''$.   This follows from the change of variables $u = \norm{x-y}^2 / 4t$.   
\end{proof} 
\begin{corollary}
The form $P(0,T) \in \Omega^\ast(\op{Conf}_2(\R^n) ) \otimes \mf g^{\otimes 2}$ extends to a smooth form in $ \Omega^\ast(\br{\op{Conf}}_2(\R^n)) \otimes \mf g^{\otimes 2} $. 
\end{corollary}
\begin{proof}
The map $\pi : \op{Conf}_2(\R^n) \to S^{n-1}$ extends to a map $\br{\op{Conf}}_2(\R^n) \to S^{n-1}$, which implies that the form $ \pi^\ast \op{Vol}_{S^{n-1}}$ extends to a smooth form on $\br{\op{Conf}}_2(\R^n)$.  It remains to show that the $\eps \to 0$ limit of $\int_{ \norm{x - y}^2 / 4 T }^{\norm{x - y}^2 / 4 \eps }  u^{n/2 - 1 } e^{- u }  \d u$ is a smooth function on $\br{\op{Conf}}_2(\R^n)$.   Note that this $\eps \to 0$ limit is the incomplete gamma function
$$
\int_{ \norm{x - y}^2 / 4 T }^{\infty }  u^{n/2 - 1 } e^{- u }  \d u  = \Gamma (n/2, \norm{x - y}^2 / 4 T ) .
$$
 The function $t \to \Gamma(n/2, t^2)$ is smooth, and the function $\norm{x - y} : \br{\op{Conf}}_2(\R^n) \to \R$ is smooth, which implies that $\Gamma( n/2, \norm{x- y}^2 / 4T)$ is smooth.  
\end{proof} 

\begin{lemma}
The differential of $P(0,T)$ on $\br{\op{Conf}}_2 (\R^n)$ is $-K_T$.
\end{lemma}
\begin{proof}
If we think of $P(0,T)$ as a de Rham current on $\R^n \times \R^n$, we know its differential is the delta current on the diagonal, minus $K_T$.  Since $P(0,T)$ is a smooth form on $\br{\op{Conf}}_2(\R^n)$, its differential must be smooth, and is determined by its restriction to $\op{Conf}_2(\R^n)$, where it is equal to $-K_T$.   
\end{proof} 

As before, we can attempt to define
$$
\Gamma( P(\eps,T), S_{CS} ) 
$$
as an element of $\Oo(\E_c, \C[[\hbar]])$, the space of functionals on the compactly supported fields.  It is not completely obvious that this is well defined. However,
\begin{proposition}
Fix $0 < T < \infty$ and $i \ge 0, k > 0$. 
Then $\Gamma_{(i,k)}(P(\eps,T), S_{CS})$  is a well-defined continuous linear functional $\E_c^{\otimes k} \to \C$. Further, the limit $\lim_{\eps \to 0} \Gamma_{(i,k)}(P(\eps,T), S_{CS})$ exists.
\label{prop cs convergence}
\end{proposition}
This proposition implies that all the counter-terms vanish. 
\begin{proof}
Let $\alpha \in \E_c^{\otimes k}$.
Standard Feynman graph techniques allow one to write
$$
\Gamma_{(i,k)}(P(\eps,T) , S_{CS} )(\alpha) = \sum_{\gamma} \frac{1}{\Aut(\gamma)} w_\gamma(\eps, T, \alpha)
$$
where the sum is over connected trivalent graphs $\gamma$, whose first Betti number is $i$, with $k$ external edges.   The weight $w_\gamma$ attached to each graph is a certain integral;  we will show that these integrals converge absolutely, even when we set $\eps =0$. 

Let $\gamma$ be such a trivalent graph, and let $V(\gamma), E(\gamma)$ denote the sets of vertices and edges of $\gamma$. Vertices are all trivalent; the end points of the external edges are not considered vertices.  Also, the external edges are not elements of $E(\gamma)$.

We will define a differential form $\omega_\gamma(\eps, T, \alpha)$ on the space $\br{\op{Conf}}_{V(\gamma)}(\R^n)$.  I'll ignore signs, as we only want to show convergence. (Technically, the form $\omega_\gamma$ is associated to a trivalent graph with a certain orientation).  

Only graphs with no loops (i.e.\ edges joining a vertex to itself) appear in the sum.  The weights of graphs with loops vanish, because the propagator $P(\eps,T)$  (which is form on $\R^n \times \R^n$)  is zero when restricted to the diagonal.  

For each edge $e \in E(\gamma)$, let
$$
\pi_e : \br{\op{Conf}}_{V(\gamma)}(\R^n) \to \br{\op{Conf}}_2(\R^n)
$$
be the projection correspond to the two vertices attached to $e$.  

Let 
$$\phi : \br{\op{Conf}}_{V(\gamma)}(\R^n) \to \R^{n k}$$
be the projection corresponding to the $k$ vertices of $\gamma$ which are attached to external edges. 

The form attached to the graph is
$$
\omega_\gamma (\eps,T , \alpha) = \otimes_{v \in V(\gamma)} \Tr^{\mf g}_v \left(  \phi^\ast \alpha \wedge_{e \in E(\gamma)} \pi_e^\ast P(\eps,T) \right).
$$
Let me explain this notation.  The expression inside the bracket is an element  
$$ \phi^\ast \alpha \wedge_{e \in E(\gamma)} \pi_e^\ast P(\eps,T) \in \Omega^\ast(\br{\op{Conf}}_{V(\gamma)} (\R^n) ) \otimes \mf g^{\otimes H(\gamma)}, $$ where $H(\gamma)$ is the set of half-edges (or germs of edges) of $\gamma$. Half-edges include external edges. Let $H(v)$ denote the set of $3$ half-edges at a vertex $v$.   For each $v$,  there is a trace map 
$$
\Tr_v^{\mf g} : \mf g^{\otimes H(v) } \to \C
$$
defined by $\ip{X, [Y,Z]}$.  Thus, 
$$
\otimes_{v \in V(\gamma)} \Tr_v^{\mf g} : \mf g^{\otimes H(\gamma)} \to \C 
$$
so that
$$
\omega_\gamma (\eps,T , \alpha) \in \Omega^\ast(\br{\op{Conf}}_{V(\gamma)} (\R^n) ) .
$$

The weight attached to $\gamma$ is
$$
w_\gamma(\eps,T, \alpha) = \int_{\br{\op{Conf}}_{V(\gamma)}} \omega_\gamma(\eps,T,\alpha).
$$
We have to show that this integral converges absolutely, for all $\eps \ge 0$ and all $0 < T < \infty$.  

The problems that could occur would be  if some of the points in the configuration space went to $\infty$. However, we are putting a compactly supported form at the vertices which are attached to external edges, and (since $k > 0$) there is at least one such vertex.     If the two points in $\R^{n}$ attached to the end points of an edge $e$ become very far apart, then $\pi_e^\ast P(\eps,T)$ decays exponentially.   If the point in $\R^n$ attached to one of the vertices goes to $\infty$, and another point is constrained to lie in a compact set because of the compactly supported form, then we must have an edge $e$ whose end points are far apart.  Thus, if any of the points in $\R^n$ attached to vertices tend to $\infty$, the integrand decays exponentially. This implies that the integral converges absolutely.
\end{proof} 
This proposition implies that
$$
\Gamma(P(0,T), S_{CS} ) \in \Oo( \E_c, \C[[\hbar]] ) / \C[[\hbar]] 
$$
is well-defined. Note that we are working modulo the ring of constants. 

To prove our result, it remains to show that the appropriate quantum master equation is satisfied. 
\begin{proposition}
$$
(\d_{DR}  +\hbar \Delta_T ) \exp ( \Gamma(P(0,T), S_{CS} ) / \hbar )  = 0
$$
modulo constants.
\end{proposition}
\begin{proof}
We can re-express the quantum master equation as 
$$
\d_{DR} \Gamma(P(0,T), S_{CS}  ) +  \{\Gamma(P(0,T), S_{CS}  ), \Gamma(P(0,T), S_{CS}  )\}_T   +\hbar \Delta_T \Gamma(P(0,T), S_{CS}  ) = 0
$$
where $\{\quad\}_T$ is a certain bracket on the space of functionals.   Up to sign, if $f, g \in \Oo(\E)$ are functionals, and $K_T = \sum \phi' \otimes \phi''$, then
$$
\{f,g\}_T = \sum \frac{\partial f }{\partial \phi'} \frac{\partial g }{\partial \phi''} .
$$
Working modulo constants amounts to ignoring the terms $\hbar \Delta_T \Gamma_{(i,2)}(P(0,T), S_{CS} )$ for any $i$,  and $\{ \Gamma_{(i,1)}(P(0,T), S_{CS} ) , \Gamma_{(j,1)}(P(0,T), S_{CS} ) \}_T$ for any $i,j$. 

In the proof of Proposition \ref{prop cs convergence} we showed how to express $\Gamma_{(i,k)}(P(0,T) , S_{CS} ) (\alpha)$ as a sum over trivalent graphs, with $k$ external edges, and first Betti number $i$.  The weight attached to each trivalent graph $\gamma$  is
$$
w_\gamma(0,T,\alpha) = \int_{\br{\op{Conf}}_{V(\gamma)}} \otimes_{v \in V(\gamma)} \Tr^{\mf g}_v \left(  \phi^\ast \alpha \bigwedge_{e \in E(\gamma)} \pi_e^\ast P(0,T) \right)
$$
Here $\alpha \in \E_c^{\otimes k}$ is the input, which we put at the external edges.

By definition, 
$$\d_{DR} \Gamma_{(i,k)}(P(0,T), S_{CS})(\alpha) =  \Gamma_{(i,k)}(P(0,T), S_{CS})(\d_{DR} \alpha). $$
Stokes' theorem implies that
\begin{multline*}
w_\gamma(0,T, \d_{DR} \alpha) = - \int_{\partial \br{\op{Conf}}_{V(\gamma)}} \otimes_{v \in V(\gamma)} \Tr^{\mf g}_v \left(  \phi^\ast  \alpha \bigwedge_{e \in E(\gamma)} \pi_e^\ast P(0,T) \right) \\
+ \sum_{e \in E(\gamma)} \pm  \int_{ \br{\op{Conf}}_{V(\gamma)}} \otimes_{v \in V(\gamma)} \Tr^{\mf g}_v \left(  \phi^\ast  \alpha \wedge \pi_e^\ast (K_T) \bigwedge_{e' \neq e} \pi_e^\ast P(0,T) \right).
\end{multline*}
In the second line, we are using the fact that $\d_{DR} P(0,T) = -K_T$. 

The terms in the sum
$$
\sum_{e \in E(\gamma)} \pm  \int_{ \br{\op{Conf}}_{V(\gamma)}} \otimes_{v \in V(\gamma)} \Tr^{\mf g}_v \left(  \phi^\ast  \alpha \wedge \pi_e^\ast (K_T) \bigwedge_{e' \neq e} \pi_e^\ast P(0,T) \right)
$$
where $e$ is a separating edge cancel with terms in the graphical expansion of $$ \{\Gamma(P(0,T), S_{CS}  ), \Gamma(P(0,T), S_{CS}  )\}_T .$$  The terms where $e$ is a non-separating edge cancel with terms in the graphical expansion of $\hbar \Delta_T \Gamma(P(0,T), S_{CS}  )$.  

Thus, it remains to show that 
$$
\sum_{\gamma} \pm \int_{\partial \br{\op{Conf}}_{V(\gamma)}} \otimes_{v \in V(\gamma)} \Tr^{\mf g}_v \left(  \phi^\ast  \alpha \bigwedge_{e \in E(\gamma)} \pi_e^\ast P(0,T) \right)  = 0 . 
$$
This has been proved by Kontsevich in \cite{Kon94} when $n \ge 3$, and in  \cite{Kon03a} when $n = 2$.   Indeed, lemma 2.1 of \cite{Kon94} implies that when $n \ge 3$, the only boundary strata of $\br{\op{Conf}}_{V(\gamma)}$ which could contribute are those where precisely two points collide.  The first lemma in section 6 of \cite{Kon03a} proves the same statement when $n = 2$.  The boundary strata where only two points collide are taken care of by the Jacobi identity.
\end{proof} 

\section*{Appendix}

This appendix contains the proof of a generalised version of theorem A. 
\subsection{Statement of results}
Let $E$ be a super vector bundle on a compact manifold $M$, whose space of global sections is denoted by $\E$. As always, suppose $\E$ has an odd symplectic structure.  
Let
$$
H_0 : \E \otimes \cinfty(\Delta^d) \to \E \otimes \cinfty(\Delta^d) 
$$
be a smooth family of generalised Laplacians, parametrised by $\Delta^d$. Thus, $H_0$ is a $\cinfty(\Delta^d)$ linear map, which is an order two differential operator with respect to $M$.  The symbol of $H_0$ is a section
$$
\sigma(H_0 ) : \Gamma( T^\ast M \times \Delta^d, \End E) .
$$
The statement that $H_0$ is a family of generalised Laplacians says that $\sigma(H_0)$ is a smooth family of metrics on $T^\ast M$, parametrised by $\Delta^d$, times the identity in $\End E$.

Let 
$$
H_1 : \E \otimes \Omega^\ast(\Delta^d) \to \E \otimes \Omega^\ast(\Delta^d)
$$
be an even $\Omega^\ast(\Delta^d)$ linear map, which is a first order differential operator with respect to $M$.  Then
$$
H = H_0 + H_1
$$
is a smooth family of generalised Laplacians parametrised by the super-manifold $\Delta^d \times \R^{0,d}$. 

The results of \cite{BerGetVer92}, appendix to chapter 9, imply that there is a unique heat kernel
$$K_t \in \E \otimes \E \otimes \Omega^\ast(\Delta^d)$$ for the operator $H$.

Let
$$
D : \E \otimes \Omega^\ast(\Delta^d) \to \E \otimes \Omega^\ast(\Delta^d)
$$
be any odd $\Omega^\ast(\Delta^d)$ linear differential operator, commuting with $H$.

As before, let
$$
P(\eps,T) = \int_\eps^T (D \otimes 1) K_t \d t \in \E \otimes  \E \otimes \Omega^\ast(\Delta^d) 
$$
be the propagator.
\begin{theorem}
Let 
$$S \in \Oo_l(\E,  \Omega^\ast(\Delta^d) \otimes \mscr{A} \otimes \C[[\hbar]] )$$ be a function which, modulo $\hbar$, is at least cubic.    More explicitly, $S$ has a Taylor series expansion 
$$
S = \sum_{i,k \ge 0} \hbar^i S_{i,k}  \phi_{i,k}(\eps)
$$
where $S_{i,k} : \E \to \Omega^\ast(\Delta^d)$ is a local functional, homogeneous of degree $k$, and $\phi_{i,k}(\eps) \in \mscr{A}$.  

Then we can form
$$
\Gamma\left( P (\epsilon,T) , S  \right)  = \sum_{i\ge 0,k \ge 0} \hbar^i \Gamma_{i,k} (P (\epsilon, T),  S)
$$
as before.  

\begin{enumerate}
\item
There exist functions $f_r \in \mscr{A}$ and $\Phi_r  \in \Oo(\E, \Omega^\ast(\Delta^d) \otimes \cinfty(0,\infty))$, for   $r \in \Z_{\ge 0}$, such that there is a small $\eps$ asymptotic expansion 
$$
\Gamma_{i,k}(P(\epsilon,T),S ) (e)    \simeq \sum f_r(\eps) \Phi_r  (e,T)
$$
for all $e \in \E$.
\item
Each $\Phi_r(e,T)$ has a small $T$ asymptotic expansion 
$$
\Phi_r(e,T) \simeq \sum \Phi_{r,s} (e) g_s(T)
$$
where the functions $\Phi_{r,s}$ are local functionals of $e$, that is, $$\Phi_{r,s} \in \Oo_l(\E,\Omega^\ast(\Delta^d)).$$ The $g_s(T)$ are certain smooth functions of $T \in (0,\infty)$.  
\item
We can view the coefficients $\Phi_{r,s}$ as linear maps $\E^{\otimes k} \to \Omega^\ast(\Delta^d)$.  Thus, we can speak of the germ of $\Phi_{r,s}$ near a point $x$.  This germ only depends on the germ of the operators $H,D, S$ near $x \times \Delta^d$.  
\end{enumerate}
\end{theorem}
\begin{remark}
At one stage in the paper, we need a slight generalisation of this result, which involves a propagator $\delta K_\eps$ where $\delta$ is an odd parameter of square zero.  The proof given below incorporates this case also. 
\end{remark}

\begin{remark}
Recall that  $\mscr{A} \subset \cinfty((0,\infty))$ is the sub-algebra spanned over $\C$ by  functions of $\epsilon$ of the form
$$
f(\epsilon )  = \int_{U \subset (\epsilon,1)^n } \frac{  F(t_1,\ldots,t_n)^{1/2} } { G(t_1,\ldots, t_n)^{1/2} } \d t_1 \cdots \d t_n
$$
and functions of the form
$$
f(\epsilon )  = \int_{U \subset (\epsilon,1)^{n-1} } \frac{  F(t_1, \ldots,t_n = \eps)^{1/2} } { G(t_1,\ldots, t_n = \eps)^{1/2} } \d t_1 \cdots \d t_{n-1}
$$
where 
\begin{enumerate}
\item
$F, G \in \Z_{\ge 0} [t_1,\ldots, t_n] \setminus \{0\}$; $n$ can take on any value. 
\item
the region of integration $U$ is an open subset cut out by finitely many equations of the form $t_i^l > t_j$, for $l \in \Z$.  
\end{enumerate}

\end{remark}

Let me explain more precisely what I mean by small $\eps$ asymptotic expansion. Let $\norm{\cdot}_{l}^{\Delta^d}$  be the $C^l$ norm on the space $\Omega^\ast(\Delta^d)$, so that $\norm{\omega}_{l}^{\Delta^d}$ is the supremum over $\Delta^d$ of the sum of all order $\le l$ derivatives of $\omega$. 

Let us consider $  \Gamma_{i,k}(P(\eps,T),S ) $ as a linear map $\E^{\otimes k} \to \Omega^\ast(\Delta^d) \otimes \cinfty( \{ 0 < \eps <  T \}  ) $.   The precise statement is that for all $R , l\in \Z_{\ge 0}$ and compact subsets $K \subset (0,\infty)$, there exists $m \in \Z_{\ge 0}$ such that 
$$
\sup_{T \in K}   \norm{   \Gamma_{i,k}(P(\eps,T),S ) (\alpha)  -  \sum_{r = 0}^R f_r(\eps) \Phi_r (T,\alpha)   }_{l}^{\Delta^d}  < \eps^{R+1} \norm{\alpha}_m
$$
for all $T$ sufficiently small.  Here $\norm{\alpha}_m$ denotes the $C^m$ norm on the space $\E^{\otimes k}$. 

The small $T$ asymptotic expansion in part (2) has a similar definition.

\subsection{Expressions in terms of integrals attached to graphs}

Standard Feynman graph techniques allow one to represent each $\Gamma_{i,k}(P(\epsilon,T),S ) $ as a finite sum 
$$
\Gamma_{i,k}(P(\epsilon,T),S )  (\alpha)  = \sum_\gamma \frac{1}{\op{Aut} \gamma} w_\gamma (\alpha)
$$
where $\alpha \in \E^{\otimes k}$,  and the sum is over certain graphs $\gamma$.   The graphs $\gamma$ that appear in the sum are connected graphs, with $k$ legs (or external lines).  Each vertex is assigned a degree in $\Z_{\ge 0}$, corresponding to the power of $\hbar$ attached to the vertex.    Degree zero vertices are at least trivalent, and the sum of the degrees of the vertices, plus the first Betti number of the graph, must be equal to $i$.   There are a finite number of such graphs.  To each graph is attached a certain integral, whose value is $w_\gamma(\alpha)$.   

I won't describe in detail the formula for the particular graph integrals appearing in the expansion of $\Gamma_{(i,k)}(P(\eps,T),S)$.  Instead, I will describe a very general class of graph integrals, which include those appearing in the sum above.   The theorem will  be proved for this general class of graph integrals.

\subsection{Asymptotics of graph integrals}

\begin{definition}
A labelled graph is a connected graph $\gamma$,  with some number of legs (or external edges). 

For each vertex $v$ of $\gamma$, let $H(v)$ denote the set of half edges (or germs of edges) emanating from $v$; a leg attached at $v$ counts as a half edge.

Also, $\gamma$ has an ordering on the sets of vertices, edges, legs and on the set of half edges  attached to each vertex.  

Each vertex $v$ of $\gamma$ is labelled by an element
$$S_v \in  \op{PolyDiff}( \E^{\otimes H(v)}, \op{Dens}(M) ) \otimes \Omega^\ast(\Delta^d).$$
$S_v$ is a smooth family of polydifferential operators parametrised by $\Delta^d$.  Let $O(v)$ denote the order of $S_v$.
\end{definition}
Let $L(\gamma)$ denote the set of legs of $\gamma$.
Fix $\alpha \in \E^{\otimes L(\gamma)}$, and fix $t_e \in (0,\infty)$ for each $e \in E(\gamma)$.  Define a function $f_\gamma (t_e, \alpha)$ as follows.

Let $H(\gamma)$ denote the set of half edges of $\gamma$, so $H(\gamma) = \cup_{v \in V(\gamma)} H(v)$.  By putting $K_{t_e}$ at each edge of $\gamma$, and $\alpha$ at the legs, we get an element
$$
\alpha \otimes_{e \in E(\gamma)} K_{t_e}   \in \E^{\otimes H(\gamma)} .
$$
On the other hand, the polydifferential operators $S_v$ at the vertices define a map
$$
\int_{M^{V(\gamma) } } \otimes S_v :   \E^{\otimes H(\gamma) }   \to \op{Dens}(M)^{\otimes V(\gamma) } \otimes \Omega^\ast(\Delta^d) \xto{\int_{M^{V(\gamma) } } }      \Omega^\ast(\Delta^d).
$$
Let 
$$
f_\gamma(t_e,\alpha) = \int_{M^{V(\gamma) } } \otimes S_v \left( \alpha \otimes_{e \in E(\gamma)} K_{t_e}  \right)  \in \Omega^\ast(\Delta^d).
$$

The graph integrals $w_\gamma$ which appear in the expansion of $\Gamma_{i,k}(P(\eps,T), S)$   can be realised as finite sums of functions of the form 
$$f_\gamma (t_e, \alpha, S_v )$$
for certain local functionals $S_v$.   Terms in the sum can be multiplied by elements of $ \mscr{A}$.

\begin{theorem}
\label{graph integral expansion}
The integral
$$
F_\gamma(\epsilon, T , \alpha )  = \int_{t_e \in (\epsilon,T)^{E(\gamma)} } f_\gamma (t_e, \alpha  )  \prod \d t_e
$$
has an asymptotic expansion as $\epsilon \to 0$ of the form
$$
F_\gamma(\epsilon, T , \alpha  ) \simeq \sum f_r(\eps) \Psi_r( T,\alpha)
$$
where $f_r \in \mscr{A}$, and the $\Psi_r$ are continuous linear  maps 
$$
\E^{\otimes L(\gamma ) } \to  \cinfty((0,\infty) ) \otimes \Omega^\ast(\Delta^d) 
$$
where $T$ is the coordinate on $(0,\infty)$.

Further, each $\Psi_r(T,\alpha)$ has a small $T$ asymptotic expansion
$$
\Psi_r(T,\alpha) \simeq \sum g_r(T) \int_M \Phi_{r,k} (\alpha)
$$
where 
$$
\Phi_{r,k} \in \op{PolyDiff}( \E^{\otimes L(\gamma)}, \dens(M) ) \otimes \Omega^\ast(\Delta^d). 
$$
and $g_r$ are smooth functions of $T \in (0,\infty)$. 
\end{theorem}
The phrase ``asymptotic expansion'' is to be interpreted in the same sense as before. 

\begin{remark}

There is a variant of this result, which incorporates a propagator of the form $P(\eps,T) + \delta K_\eps$, where $\delta$ is an odd parameter of square zero.  In this case, the graphs we use may have one special edge, on which we put $K_\eps$ instead of $K_{t_e}$.  Thus, instead of integrating over the parameter $t_e$ for this edge, we specialise to $\eps$.  
\end{remark}
\subsection{Asymptotics of heat kernels}
The proof is based on the asymptotic expansion of the heat kernels of generalised Laplacians, proved in \cite{BerGetVer92}.  The operator
$$
H  : \E \otimes \Omega^\ast(\Delta^d) \to \E \otimes \Omega^\ast(\Delta^d)
$$
is a smooth family of generalised Laplacians parametrised by the supermanifold with corners $\Delta^d \times \R^{0,d}$.     A generalised Laplacian on the vector bundle $E$ on $M$ (whose global sections is $\E$) is specified by a metric on $M$, a connection on $E$ and a potential $F \in \Gamma(M, \End(E))$.    A smooth family of generalised Laplacians is a family where this data varies smoothly; this is equivalent to saying that the operator varies smoothly. We are dealing with a smooth family parametrised by the super-manifold $\Delta^d \times\R^{0,d}$. The metric on $M$ will be independent of the odd coordinates $\R^{0,d}$, but the parameters for the connection and the potential will both involve Grassmann variables.  (I'm very grateful to Ezra Getzler for explaining this point of view to me.)

Let $x \in M$.   Let $U \subset M \times \Delta^d$ denote the open subset of points $(y,\sigma)$ where $d_\sigma(x,y) < \eps$.  Normal coordinates on $U$ gives an isomorphism 
$$
U \iso B_\eps^n \times \Delta^d
$$
where $B_\eps^n$ is the ball of radius $\eps$ in $\R^n$, and $n = \dim M$.  

Thus, we get an isomorphism 
$$
\cinfty( U \times \R^{0,d} ) = \cinfty(U) \otimes \C [ \d t_1 ,\ldots, \d t_d ]  \iso \cinfty(B_\eps^n) \otimes \Omega^\ast(\Delta^d)
$$ 
of $\Omega^\ast(\Delta^d)$ algebras.    

The vector bundle $E$ becomes a vector bundle (still called $E$) on $B_{\eps}^{n}$, and we find
$$
\Gamma( U, E) \otimes \C[\d t_1, \ldots, \d t_d ] \iso \Gamma(B_\eps^n, E ) \otimes \Omega^\ast(\Delta^d) .   
$$
The following is a variant of a result proved in \cite{BerGetVer92}, following \cite{MinPle49, BerGauMaz71,McKSin67}. 
\begin{theorem}
There exists a small $t$ asymptotic expansion of $K_t$ which, in these coordinates, is of the form
$$
K_t \simeq  t^{-\op{dim} M  / 2} e^{-\norm{x-y}^2/t }  \sum_{i \ge 0} t^i \phi_i
$$
where $x,y$ denote coordinates on the two copies of $B_\eps^n$,  and
$$
\phi_i \in  \Gamma(B_\eps^n,E) \otimes \Gamma(B_\eps^n,E) \otimes \Omega^\ast(\Delta^d).
$$
If we denote
$$
K_t^N = t^{-\op{dim} M  / 2}  e^{-\norm{x-y}^2/t }  \sum_{i =0}^N t^i \phi_i
$$
then for all $l \in \Z_{\ge 0}$
$$
\norm{ K_t - K_t^N }_l  = O(t^{N - \dim M / 2 - l}).
$$
Here $\norm{\cdot}_l$ denotes the $C^l$ norm on the space $\cinfty(B_\eps^n) \otimes \cinfty(B_\eps^n) \otimes \Omega^\ast(\Delta^d)$.  
\end{theorem}
This precise statement is not proved in \cite{BerGetVer92}, as they do not use Grassmann parameters.  But, as Ezra Getzler explained to me, the proof in \cite{BerGetVer92} goes through \emph{mutatis mutandis}. 

\subsection{Proof of Theorem \ref{graph integral expansion}}
In this proof, I will often avoid mention of the parameter space $\Delta^d \times \R^{0,d}$.  Thus, if $f$ is some expression which depends on this parameter space, I will often abuse notation and write $\abs{ f }$ for the $C^l$ norm of $f$ as a function of the $\Delta^d \times \R^{0,d}$ variables, for some $l$.  

Let us enumerate the edges of $\gamma$ as $e_1,\ldots, e_k$.  Let $t_i = t_{e_i}$, and let us consider the region where $t_1 > t_2 > \cdots > t_k$.    (Of course we have to consider all orderings of the edges of $\gamma$).

For a function $I : E(\gamma) = \{1,\ldots, k\} \to \Z_{\ge 0}$, let $\abs{I} = \sum I(i)$.    Let $t^I = \prod t_i^{I(i)}$.  Similarly, if $n \in \Z$, let $t^n = \prod t_i^n$.

Let 
$$
O(\gamma) = \sum_{v \in V(\gamma)} O(v) .
$$

For $R > 1$ let 
$$
A_{R,T} \subset (0,T)^{E(\gamma)}
$$
be the region where $t_i^R < t_{j}$ for all $i,j$.  This means that the $t_i$ are all of a similar size.
\begin{proposition}
\label{prop wicks lemma expansion}
For all $r \ge 0$, there exist
\begin{enumerate}
\item
 $F_i,G_i \in \Z_{\ge 0} [ t_1,\ldots, t_k ] \setminus \{0\} $ for  $1 \le i \le c_r$,
 \item
 polydifferential operators
$$\Psi_{i} \in \op{Poly Diff} ( \E^{\otimes L(\gamma)}, \op{Dens}(M) ) \otimes \Omega^\ast(\Delta^d)     $$ 
for $1 \le i \le c_r$,
\end{enumerate}
such that
$$
\abs{ f_\gamma(t_1,\ldots, t_k) -  \sum_{i = 1}^{c_r} \frac{ F_i (t_1,\ldots, t_k )^{1/2}} { G_i (t_1,\ldots, t_k)^{1/2} } \int_M \Psi_i  (\alpha  ) } \le \norm{\alpha}_{r + 1 - \chi(\gamma) \dim M / 2 + O(\gamma)(\abs{E(\gamma)} + 1) } t_1^{r+1} 
$$
for all $\{t_1,\ldots, t_k\} \in A_{R,T}$ with $t_1 > t_2 > \cdots  > t_k$ and  $t_1$ sufficiently small. In this expression, $\chi(\gamma)$ is the Euler characteristic of the graph.
\end{proposition}

\begin{proof}
As above, let
$$
K^N_{t}(x,y) =  t^{-\op{dim} M  / 2} \Psi(x,y) e^{-\norm{x-y}^2/t }  \sum_{i = 0}^N t^i \phi_i(x,y)
$$
be the  approximation to the heat kernel to order $t^N$.  This expression is written in normal coordinates near the diagonal in $M^2$.  The $\phi_i(x,y)$ are sections of $E \boxtimes E$ defined near the diagonal on $M^2$. $\Psi(x,y)$ is a cut-off function, which is $1$ when $\norm{x-y} < \epsilon$ and $0$ when $\norm{x-y} > 2\epsilon$.

We have the bound
$$
\norm{  K_{t}(x,y) - K^N_{t}(x,y) }_l = O ( t^{N  - \op{dim} M  /2 -  l    }  ) .
$$

The first step is to replace each $K_{t_i}$ by $K^{N}_{t_i}$ on each edge of the graph.   Thus, let $f_\gamma^N(t_i,\alpha)$ be the function constructed like $f_\gamma(t_i,\alpha)$ except using $K_{t_i}^N$ in place of $K_{t_i}$. 

Each time we replace $K_{t_i}$ by $K_{t_i}^{N}$, we get a contribution of $t_i^{N  - O(\gamma) - \op{dim} M / 2}$ from the edge $e_i$, times the $O(\gamma)$ norm of the contribution of the remaining edges, times $\norm{\alpha}_{O(\gamma)}$. 

The $O(\gamma)$ norm of the contribution of the remaining edges is $\prod_{j \neq i} t_j^{ -  O(\gamma) - \op{dim} M / 2}$.  We are thus left with the bound 
$$
\abs {  f_\gamma(t_i,\alpha ) -  f_\gamma^N (t_i,\alpha )   } <  C    t^{- O(\gamma) - \dim M / 2  }    t_1^{N} \norm{\alpha}_{O(\gamma) } 
$$
where $C$ is a constant.   (Recall our notation : $t^n$ denotes $\prod t_i^n$). 

In particular, if the $t_i$ are in $A_{R,T}$, we find that 
$$ \abs { f_\gamma(t_i, \alpha  ) -   f_\gamma^N(t_i,\alpha) }   <   \norm{\alpha}_{O(\gamma)} t_1^{N -  \abs{E(\gamma)} ( O(\gamma) + \dim M / 2 )  R + 1 }  $$
if $t_i \in A_{R,T}$ and $t_i$ are all sufficiently small.

Next, we construct a  small $t_i$ asymptotic expansion of $f_\gamma^N(t_i,\alpha)$.   Recall that $f_\gamma^N(t_i,\alpha)$ is defined as an integral over a small neighbourhood of the small diagonal in $M^{V(\gamma)}$.  Let 
$$n = \op{dim} M .$$  
By using a partition of unity we can consider $f_\gamma^N(t_i,\alpha)$ as an integral over a small neighbourhood of zero in $\R^{n V(\gamma)}$.  This allows us to express  $f_\gamma^N(t_i,\alpha)$ as a finite sum of integrals over $\R^{n V(\gamma)}$, of the following form.

For each vertex $v$ of $\gamma$, we have a coordinate map $x_v : \R^{n \abs{V(\gamma)}} \to \R^n$.  Fix any $\eps> 0$, and let $\chi : [0,\infty) \to [0,1]$ be a smooth function with $\chi(x) = 1$ if $x < \epsilon$, and $\chi(x) = 0$ if $x > 2 \epsilon$.  Let us define a cut-off function $\psi$ on $\R^{n V(\gamma)}$ by the formula
$$
\psi = \chi (  \norm{\sum x_v}^2  ) \chi \left( \sum_{v' \in V(\gamma)}  \norm{ x_{v'} -  \abs{V(\gamma)}^{-1} \sum_{v \in V(\gamma)} x_v }^2  \right). 
$$
Thus, $\psi$ is zero unless all points $x_v$ are near their centre of mass $\abs{V(\gamma)}^{-1}\sum x_v$, and $\psi$ is zero when this centre of mass is too far from the origin.  

 For each $1 \le i \le k$, let $Q_i$ be the quadratic for on $\R^{n \abs{V(\gamma)}}$ defined by
$$
Q_i ( x) = 
\begin{cases} 
0 & \text{ if the edge } e_i \text{ is a loop, i.e.\ is attached to only one vertex } \\
\norm{ x_{v_1} - x_{v_2} }^2 & \text{ if } v_1, v_2 \text{ are the vertices attached to the edge } e_i 
\end{cases}
$$

We can write $f_\gamma^N$ as a finite sum of integrals of the form
$$
\int_{\R^{n V(\gamma)} } \psi    e^{- \sum Q_i  / t_i } \sum_{ I, K} t^{ I -  \dim M / 2 - O(\gamma) } \Phi_{I,K} \partial_{K,x} \alpha.
$$
In this expression, 
\begin{itemize}
\item
The sum is over $I : E(\gamma) \to \Z_{\ge 0}$, with all $I(e) \le N + O(\gamma) + 1$, and multi-indices 
$K : V(\gamma) \times \{1,\ldots,n\} \to \Z_{\ge 0}$, with $\sum K(v,i) \le O(\gamma)$.  The notation $\partial_{K,x}$ denotes 
$$\partial_{K,x}=  \prod_{v \in V(\gamma), 1 \le i \le n}  \frac {\partial}{\partial x_{v,i}^{K(v,i)} }.$$

In this notation, we are pretending (by trivialising the vector bundle $E$ on some small open sets in $M$) that $\alpha$ is a function on $\R^{n V(\gamma)}$. 
\item
The $\Phi_{I,K}$ are smooth functions on $\R^{n V(\gamma)}$.
\end{itemize}

Next, we will use Wick's lemma to compute the asymptotics of this integral.   Let $c = \left(1/\abs{V(\gamma)} \right)\sum x_v$ be the centre of mass function $\R^{n V(\gamma)}\to \R^n$.       We can perform the integral in two steps, first integrating over the coordinates $y_v =  x_v -  c$,  and secondly by integrating over the variable $c$.  (Of course, there are $\abs{V(\gamma)} - 1$ independent $y_v$ coordinates). The quadratic form $\sum Q_i / t_i$ on $\R^{n V(\gamma)}$ is non-degenerate on the subspace of $\R^{n V(\gamma)}$ of vectors with a fixed centre of mass, for all $t_i \in (0,\infty)$.  Thus, the integral over the variables $y_v$ can be approximated with the help of Wick's lemma.

Let us order the set $V(\gamma)$ of vertices as $v_1,v_2, \ldots, v_m$.  We will use the coordinates $y_1,\ldots,y_{m-1}$, and $c$ on $\R^{n V(\gamma)}$.  Then  $f^N_\gamma$  is a finite sum of integrals of the form
$$
\int_{w \in \R^n} \chi ( \abs{w}^2 ) \int_{y_1,\ldots,y_{m-1} \in \R^n} \chi ( \sum \norm{y_v}^2 )  e^{- \sum Q_i(y) / t_i} \sum t^{I - O(\gamma)  - \dim M / 2} \Phi_{I,K} \partial_{K, w,y_i} \alpha.
$$
Here we are using the same notation as before, in these new coordinates.    

To get an approximation to the inner integral, we take a Taylor expansion of the functions $\alpha$ and $\Phi_{I,K}$ around the point $\{y_i = 0, w\}$, only expanding in the variables $y_i$.  We take the expansion to order $N$.  We find, as an approximation to the inner integral, an expression of the form
$$
\int_{y_1,\ldots,y_{m-1} \in \R^n} e^{- \sum Q_i(y) / t_i} \sum t^{I -  \dim M / 2 - O(\gamma)} y^K c_{K,I,L} \left( \prod \frac{\partial^{L_i}} {\partial^{L_i} y_i} \frac{\partial ^{L_w}} {\partial^{L_w} w}  \alpha  \right)_{y_i = 0}
$$
where the sum is over a finite number of multi-indices $I,K,L$, and $c_{K,I,L}$ are constants.

We can calculate each such integral by Wick's lemma.   The application of Wick's lemma involves inverting the quadratic form $\sum Q_i(y) / t_i$.    Let $A = A(t_i)$ denote the matrix of the quadratic form $\sum Q_i(y) / t_i$; this is a square matrix of size $(\dim M) \abs{V(\gamma)}$, whose entries are sums of $t_{i}^{-1}$.   Note that $\left( \prod_{i = 1}^{k} t_i \right)  A$ has polynomial entries.  

Let 
$$
P_\gamma (t_i) = \det \left( \left( \prod_{i = 1}^{k}t_i \right) A \right). 
$$
This is the graph polynomial associated to $\gamma$ (see \cite{BloEsnKre06}).   One important property of $P_\gamma$ is that it is a sum of monomials, each with a non-negative integer coefficient. 

We can write
$$
A^{-1} = P_\gamma^{-1} B
$$
where the entries of $B$ are polynomial in the $t_i$.   Note also that 
$$
\det A = P_\gamma t^{-(\dim M)(\abs{V(\gamma)} - 1)  / 2 } .
$$
(The expansion from Wick's lemma gives an overall factor of $(\det A)^{-1/2}$).

Thus, we find, using Wick's lemma, an approximation of the form 
$$
f^N_\gamma(t_i) \simeq P_\gamma^{-1/2} \sum_{l \ge 0, I: E(\gamma) \to \tfrac{1}{2}\Z_{\ge 0}} P_\gamma^{-l} t^{I - \dim M / 2 - O(\gamma) } \int_M \Psi_{l,I}(\alpha)
$$
where the $\Psi_{l,I}$ are polydifferential operators
$$
\Psi_{l,I} : \E^{\otimes L(\gamma)} \to \dens(M)
$$
 and the sum is finite (i.e.\ all but finitely many of the $\Psi_{l,I}$ are zero).    
 
This expansion is of the desired form; it remains to bound the error term.
\begin{lemma}
The error term in this expansion is bounded by $$t_1^{R (N+1 + (\dim M) \chi(\gamma)/2 - \abs{E(\gamma)} O(\gamma) ) } \norm{\alpha}_{N + 1 + O(\gamma)}$$ for $N$ sufficiently large and $t_1$ sufficiently small.  Here $\chi(\gamma)$ is the Euler characteristic of the graph. 
\end{lemma}
\begin{proof}

The error in this expansion arises from the error in the Taylor expansion of the functions $\alpha, \Phi_{I,K}$ around $0$, and from the fact that we are neglecting the cut-off function.  Thus, if $N$ is sufficiently large, the magnitude of the error in the expansion can be bounded by an expression of the form
$$
t^{- \dim M / 2 - O(\gamma) }\int_{w \in \R^n} \chi(w)   \int_{y_1,\ldots, y_{m-1} \in  \R^{n}  } \abs{ \prod_e e^{-\sum Q_i(y)/t_i} \sum_{K} \phi_K \partial_{K,w,y_j} \alpha  } . 
$$
Here, the sum is over multi-indices $K$ for the variables $y_i,w$, each $\phi_K$ is homogeneous of order $N+1$ as a function of the variables $y_i$,  and $\abs{K} \le N+1+O(\gamma)$, so we are differentiating $\alpha$ at most this number of times.

This integral only decreases if we  decrease each $t_i$.  Since $t_i > t_1^{R}$, we find we can bound the integral by the corresponding integral using the quadratic form $\sum Q_i(y) / t_1^{R}$. Thus, we find a bound of   
$$
t_1^{- R \abs{E(\gamma)} ( \dim M / 2 + O(\gamma) )}\int_{w \in \R^n} \chi(w)   \int_{y_1,\ldots, y_{m-1} \in  \R^{n}  } \abs{ \prod_e e^{-\sum Q_i(y)/t_1^R} \sum_{K} \phi_K \partial_{K,w,y_j} \alpha  } . 
$$
Wick's lemma bounds this integral by
$$
\op{det}(\sum Q_i(y)/ t_1^{R})^{-1/2} t_1^{R (N+1 - \abs{E(\gamma)}(\dim M /  2 + O(\gamma)) )} \norm{\alpha}_{N+1 + O(\gamma) }
$$
Also, 
$$
\op{det}(\sum Q_i(y)/ t_1^{R})^{-1/2} = C t_1^{  R\op{dim} M  \abs{V(\gamma)} /2    }
$$
where $C$ is a constant, to yield the desired bound. 
\end{proof} 

\end{proof} 

A \emph{subgraph} $\gamma'$ of a graph $\gamma$ is given by the set  of edges $E(\gamma') \subset E(\gamma)$.   The vertices of the subgraph $\gamma'$  are the ones that adjoin edges in $E(\gamma')$.     The legs of $\gamma'$ are the half-edges of $\gamma$, which adjoint vertices of $\gamma'$, but which are not part of any edge of $\gamma'$.  

Let us fix a proper subgraph $\gamma'$, and let us enumerate the edges of $\gamma$ as $e_1,\ldots, e_k$,  where $e_1,\ldots, e_l \in E(\gamma')$ and $e_{l+1},\ldots, e_k \in E(\gamma) \setminus E(\gamma')$.  Let $t_i = t_{e_i}$.

Let 
$$ \phi_{ \gamma,\gamma' } (t_{l+1}, \ldots, t_k ) =  \int_{ t_1,\ldots, t_l \in (t_{l+1}^{1/R} , T)   }  f_\gamma(t_1,\ldots, t_k) \d t_1 \cdots \d t_l  $$
\begin{lemma}
Let $R \gg 0$ be sufficiently large. Then for all $r > 0$,  there exist $m_r \in \Z_{\ge 0}$,  a finite number of $g_i \in \mscr{A}$, $F_i, G_i \in \Z_{\ge 0 }[t_{l+1},\ldots,t_k]\setminus\{0\}$, and continuous linear maps $\Psi_i : \E^{\otimes L(\gamma)} \to \Omega^\ast(\Delta^d) \otimes \cinfty((0,\infty))$ such that 
$$
\abs{ \phi_{\gamma,\gamma'} (t_{l+1},\ldots, t_k ) -  \sum g_i(t_{l+1} )\frac{ F_i (t_{l+1}, \ldots, t_k )^{1/2} } {G_i(t_{l+1}, \ldots, t_k)^{1/2} }  \Psi_i ( T,\alpha)     } < \norm{\alpha}_{m_r} t_{l+1}^{r+1}
$$
for all $t_{l+1} >  \cdots >  t_k > 0$ with $t_{l+1}^R < t_k$, and $t_{l+1}$ sufficiently small.

Further, the $\Psi_i$ admit small $T$ asymptotic expansions 
$$
\Psi_i(T,\alpha) \simeq \sum \eta_{i,k}(T)\int_M \Phi_{i,k}(\alpha)
$$
where 
$$\Phi_{i,k}  \in \op{Poly Diff}(\E^{\otimes L(\gamma)}, \dens(M) ) \otimes \Omega^\ast(\Delta^d) $$
and $\eta_{i,k}(T)$ are smooth functions of $T$.
\label{lemma subgraph}
\end{lemma}

\begin{proof}
We can write  
$$f_\gamma(t_1,\ldots, t_k, \alpha)  = f_{\gamma'}(t_{l+1}, \ldots, t_k, \alpha \otimes K_{t_1} \otimes \cdots \otimes K_{t_l}). $$ The right hand side of this equation denotes the graph integral for $\gamma'$ with inputs being tensor products of the heat kernels $K_{t_i}$, for $1 \le i \le l$, and $\alpha$.

The starting point of the proof is Proposition \ref{prop wicks lemma expansion} applied to each connected component of the graph $\gamma'$.  This proposition implies that we can approximate the contribution of each such connected component by a local functional applied to the legs of that connected component.  If  we approximate the contribution of each connected component in this way, we are left with a graph integral $\Psi_{\gamma/\gamma'}(t_1,\ldots, t_l, \alpha)$ for the quotient graph $\gamma / \gamma'$, with certain local functionals at the vertices of $\gamma / \gamma'$.   

More formally, Proposition \ref{prop wicks lemma expansion} implies that  there exists a finite number of $F_i, G_i \in \Z_{\ge 0}[t_{l+1},\ldots,t_k]\setminus\{0\}$, and $\Psi^i_{\gamma / \gamma' } (t_1,\ldots, t_l , \alpha )$, which are graph integrals for $\gamma / \gamma'$,  such that 
\begin{multline*}
\abs{f_{\gamma}(t_{1}, \ldots, t_k, \alpha)  -  \sum \frac{ F_i (t_{l+1}, \ldots, t_k )^{1/2} } {G_i(t_{l+1}, \ldots, t_k)^{1/2} }    \Psi^i_{\gamma / \gamma' } (t_1,\ldots, t_l , \alpha )  } \\ < t_{l+1}^{r+1} \norm{\alpha \otimes K_{t_1} \otimes \cdots \otimes K_{t_l}}_{r + 1 -\chi(\gamma) \dim M / 2 + O(\gamma)(\abs{E(\gamma)} + 1) }  
\end{multline*}
for all $t_1,\ldots,t_k$ with $t_{l+1} > \ldots > t_k$, $t_{l+1}^R < t_k$, and $t_{l+1}$ sufficiently small.

We are only interested in the region where $t_i^R > t_{l+1}$ if $1 \le i \le l$.  Note that 
\begin{align*}
\norm{K_{t_i}}_{r + 1 -\chi(\gamma) \dim M / 2 + O(\gamma) (\abs{E(\gamma)}+1)} &= O(t_i^{(\chi(\gamma) -1)\dim M/2    - r -1 - O(\gamma)(\abs{E(\gamma)}+1) }  )\\ & = O( t_{l+1}^{R^{-1}\left(( \chi(\gamma) -1)\dim M/2    - r -1 - O(\gamma)(\abs{E(\gamma)}+1)  \right) }   )
\end{align*}
Since we can take $R$ as large as we like, this contribution is small, and we find, 
$$
\abs{f_{\gamma}(t_{1}, \ldots, t_k, \alpha)  -  \sum \frac{ F_i (t_{l+1}, \ldots, t_k )^{1/2} } {G_i(t_{l+1}, \ldots, t_k)^{1/2} }    \Psi^i_{\gamma / \gamma' } (t_1,\ldots, t_l , \alpha )  } \\ < t_{l+1}^{r ( 1 - 1 / R ) } \norm{\alpha}_{m_r}  
$$
for  some $m_r \gg 0$.  

We want to integrate over the variables $t_1,\ldots, t_l$, in the region $(t_{l+1}^{1/R} , T)$.  The integral 
$$
\int_{t_1,\ldots, t_l \in (t_{l+1}^{1/R}, T) } \Psi^i_{\gamma / \gamma' } (t_1,\ldots, t_l , \alpha )   \d t_1 \cdots \d t_l
$$
is approximated (by induction) using Theorem \ref{graph integral expansion}, applied to the graph $\gamma / \gamma'$, with $t_{l+1}^{1/R}$ playing the role of $\eps$.   It is easy to see that this yields the desired approximation of $\phi_{\gamma , \gamma ' }(t_{l+1},\ldots, t_k,\alpha)$. 
\end{proof}

Let $\gamma' \subset \gamma$ be a subgraph.  
Let $A_{R,T}^{\gamma'} \subset (0,T)^{E(\gamma)}$ be the open subset where 
\begin{align*}
t_e^R > t_{e'} & \text{ if } e \in E(\gamma) \setminus E(\gamma') \text{ and } e' \in E(\gamma') \\
t_e^R < t_{e'} & \text{ if } e,e' \in E(\gamma') .
\end{align*} 
This means that the lengths of the edges of the subgraph $\gamma'$ are all around the same size, and are all much smaller than the lengths of the other edges. 

\begin{lemma}
Fix $R \gg 0$, and $0 < T < \infty$.  Then the closures of the regions
 $A_{R^{2^k},T}^{\gamma'} $, 
where $0 \le k \le \abs{E(\gamma)}$ and $\gamma' \subset \gamma$ is non-empty, cover $(0,T)^{E(\gamma)}$. (The regions $A_{R^{2k},T}$ appear as $A_{R^{2k}, T}^{\gamma}$, where $\gamma$ is considered as a subgraph of itself). 
\end{lemma}
\begin{proof}
Let $\{t_e\} \in (0,T)^{E(\gamma)}$.  As before, label the elements of $E(\gamma)$ by $\{1,2,\ldots, k \}$,  in such a way that $t_1 \ge t_2 \ge \cdots \ge t_{k}$.  

Either $t_j^R  \ge t_{k}$ for all $j < k$, or, there is a smallest $i_1 < k$ such that $t_{i_1}^{R} \le t_{k}$.  In the first case, we are done, as then $\{t_e\} \in \br{ A_{R,T}^{\gamma'}} $ where $\gamma'$ is the subgraph with the single edge corresponding to $t_{k}$.  

Suppose the second possibility holds.  Then either for all $j < i_1$,  $t_j^{R} \ge t_{i_1}$,   or there exists a smallest $i_2 < i_1$ with $t_{i_2}^R \le t_{i_1}$.  In the first case, we're done, as we are in $\br{ A_{R,T}^{\gamma'} }$ where $\gamma'$ is the subgraph whose edges correspond to $t_{i_1}, t_{i_1+1}, \ldots, t_{k}$.  

Again, let's suppose the second possibility holds.  Then $t_{i_2}^{R^2} \le t_{i_1}^R \le t_{k}$.    Either, for all $j < i_2$, $t_j^{R^2} \ge t_{i_2}$,  and then we are in $\br{ A_{R^2,T}^{\gamma'} }$ where $\gamma'$ is the subgraph whose edges correspond to $t_{i_2}, t_{i_2+1} , \ldots, t_k$.

Otherwise, there is some smallest $i_3 < i_2$ with $t_{i_3}^{R^2} \le t_{i_2}$.  Then $t_{i_3}^{R^4} \le t_{k}$.   And so forth.

We eventually end up either finding ourselves in one of the regions $\br { A_{R^{2^k}, T}^{\gamma'} }$, for some non-empty proper subgraph $\gamma' \subset \gamma$, or we find some $i_{k+1} = 1$, so we are in  $\br{ A_{R^{2^{k}}, T}^{\gamma} } = \br{ A_{R^{2^{k}}, T} }$.  
\end{proof} 
\begin{definition}
An open subset $U$ of $(0,T)^{E(\gamma)}$ is called \emph{good} if it is a subset of some  $A_{R,T}^{\gamma'}$ which is cut out by a finite number of inequalities $t_e^{R^n} > t_{e'}$, where $l \in \Z_{\ge 0}$ and both $e,e' \in E(\gamma')$. 
\end{definition}
\begin{lemma}
The intersection of any two good regions is good. 
\end{lemma}
Now, we have seen that $(\eps,T)^{E(\gamma)}$ is covered by the closures of a finite number of good regions.  Thus, we can write 
$$F_\gamma (\eps , T, \alpha ) =  \int_{(\eps,T)^{E(\gamma)} } f_\gamma ( t_e ,\alpha ) \prod_{e \in E(\gamma)} \d t_e $$ as an alternating sum of integrals of $f_\gamma(t_e)$ over  $U \cap (\eps,T)^{E(\gamma)}$, where $U$ is a good subset of $(0,T)^{E(\gamma)}$. Thus, in order to understand the small $\eps$ asymptotic expansions of $F_\gamma(\eps, T,\alpha)$, it suffices to consider the integrals of $f_\gamma$ over such regions.

\begin{lemma}
Let us fix $R$ to be a (sufficiently large) integer.

Let $U \subset (0,T)^{E(\gamma)}$ be a good subset.  Then the integral
$$
F_{\gamma,U}(\eps,T,\alpha) = \int_{U \cap (\eps,T)^{E(\gamma)} } f_\gamma ( t_e ,\alpha ) \prod_{e \in E(\gamma)} \d t_e
$$
admits a small $\eps$ asymptotic expansion
$$F_{\gamma,U}(\eps,T,\alpha)  \simeq \sum \phi_r(\eps) \Psi_r( T,\alpha)$$
where $\phi_r \in \mscr{A}$ and $\Psi_r : \E^{\otimes L(\gamma)} \to \Omega^\ast(\Delta^d)\otimes \cinfty((0,T))$ are continuous maps.

Each $\Psi_r(T,\alpha)$ admits a small $T$ asymptotic expansion
$$
\Psi_r(T,\alpha) \simeq \sum g_k(T)\int_M \Phi_{r,k}(\alpha)
$$
where $g_k$ is a smooth function of $T$, and $\Phi_{r,k} \in \op{Poly Diff}(\E^{\otimes L(\gamma)}, \dens(M) ) \otimes \Omega^\ast(\Delta^d)$. 
\end{lemma}
\begin{proof}
This  follows from Proposition \ref{prop wicks lemma expansion} and Lemma \ref{lemma subgraph}.  

Indeed, we can assume (without loss of generality) that 
$$
U \subset  \{ t_1,\ldots, t_k \mid t_i^N > t_{l+1} \text{ for all } 1 \le i \le l  , \quad t_{l+1} >  t_{l+2} > \cdots > t_k, \quad t_{l+1}^N < t_k     \} 
$$
is an open subset cut out by a finite number of inequalities of the form $t_i^{a} > t_j$, where $a \in \Z_{> 0}$ and $l+1 \le i,j \le k$.

Then let
$$
\eta(t_{l+1}, T, \alpha)  = \int_{(t_1,\ldots, t_k) \in U} f_\gamma( t_1,\ldots, t_k , \alpha) \d t_1 \cdots \what{\d t_{l+1} } \cdots \d t_k .
$$
Here, we're integrating over all variables except $t_{l+1}$.

Then  Proposition  \ref{prop wicks lemma expansion} (if $l = 0$) or Lemma \ref{lemma subgraph} (if $l > 0$) imply that $\eta(t_{l+1}, T,\alpha)$ has a nice $t_{l+1}$ asymptotic expansion. More precisely, for all $r \in \Z_{> 0}$, there exist $m_r \in \Z_{> 0}$, and a finite number of functions $g_i(t_{l+1})$, $\psi_i(T,\alpha)$, such that
\begin{enumerate}
\item
$$\abs{ \eta(t_{l+1}, T, \alpha) -  \sum g_i(t_{l+1}) \psi_i(T,\alpha)  } < t_{l+1}^{r+1} \norm{\alpha}_{m_r}   $$
\item
the $\psi_i$ are continuous linear maps $\psi_i : \E^{\otimes L(\gamma)} \to \Omega^\ast(\Delta^d) \otimes \cinfty((0,\infty))$, which admit a small $T$ asymptotic expansion 
$$
\psi_i(T,\alpha) \sim \sum \zeta_{i,k}(T) \int_M \phi_{i,k}(\alpha)
$$ where 
$$\phi_{i,k}  \in \op{Poly Diff}(\E^{\otimes L(\gamma)}, \dens(M) ) \otimes \Omega^\ast(\Delta^d) $$
and $\zeta_{i,k}(T)$ are smooth functions of $T$.
\item
$$
\int_{\eps}^1 g_i(t_{l+1}) \d t_{l+1} \in \mscr{A}
$$
\end{enumerate}
(The third part follows from the particular form of the terms of the expansions proved in Proposition \ref{prop wicks lemma expansion} and Lemma \ref{lemma subgraph}.)

Now, it remains to integrate out the variable $t_{l+1}$.  Let
$$
\eta_r (t_{l+1}, T, \alpha) =  \sum g_i(t_{l+1}) \psi_i(T,\alpha)
$$
be the approximation to order $t_{l+1}^{r+1}$ to $\eta(t_{l+1}, T,\alpha)$.  

Then
$$
\abs{ \int_\eps^T \eta \d t_{l+1}  - \left( \int_0^T (\eta - \eta_r) \d t_{l+1}   + \int_1^T \eta_r \d t_{l+1} +    \int_\eps^1 \eta_r \d t_{l+1}  \right) } < \eps^{r+2} \norm{\alpha}_{m_r} .
$$
Thus, 
$$\int_0^T (\eta - \eta_r) \d t_{l+1}   + \int_1^T \eta_r \d t_{l+1} +    \int_\eps^1 \eta_r \d t_{l+1} $$
gives the desired small $\eps$ asymptotic expansion.  Note that the integral in the first term converges.   The first two terms are independent of $\eps$, and the third term is in $\mscr{A} \otimes \Hom ( \E^{\otimes L(\gamma)}, \Omega^\ast(\Delta^d) \otimes \cinfty((0,\infty)) )$, as desired. 

It's easy to check that the small $T$ asymptotic expansion of the approximation above is in terms of local functionals of $\alpha$, as desired.

\end{proof} 
This completes the proof of Theorem \ref{graph integral expansion}.    

The proof of the variant result, when we include a propagator of the form $P(\eps,T) + \delta K_\eps$, is identical, except that instead of integrating over the smallest variable $t_k$, we specialise $t_k = \eps$.  In our definition of the algebra $\mscr{A}$, we included functions of $\eps$ of the form
$$
f(\epsilon )  = \int_{U \subset (\epsilon,1)^{k-1} } \frac{  F(t_1,\ldots,t_k = \eps)^{1/2} } { G(t_1,\ldots, t_k = \eps)^{1/2} } \d t_1 \cdots \d t_{k-1}
$$
where $F, G \in \Z_{\ge 0}[ t_1, \ldots, t_k]$ and $U$ is cut out by equations of the form $t_i^l > t_j$, for $i,j < k$.  Functions of this kind arise when we have a $\delta K_\eps$ term in the propagator.


\def\cprime{$'$}

\end{document}